\documentclass{article}
\usepackage{pstricks, amssymb, xr, multind, multicol}
\def\diagram#1{\def\normalbaselines {\lineskip=5pt\baselineskip=0pt
\lineskiplimit=1pt}\matrix{#1}}
\def\hfl#1{\smash{\mathop{\hbox to 10mm{\rightarrowfill}}\limits^{\textstyle
#1}}}
\def\vfl#1{\llap{$\scriptstyle#1$}\left\downarrow\vbox to 6mm{}\right.}

\def\uvfl#1{\llap{$\scriptstyle#1$}\left\uparrow\vbox to 6mm{}\right.}
\def\norm#1{\parallel #1 \parallel}

\makeindex{N}
\makeindex{T}
\newcommand{\indexT}[1]{#1\index{T}{#1}}

\newtheorem{proposition}[equation]{Proposition}

\newtheorem{theorem}[equation]{Theorem}
\newtheorem{exa}[equation]{Example}
\newtheorem{ex}[equation]{Exercise}
\newtheorem{s-ex}[equation]{Side-exercise}
\newtheorem{exas}[equation]{Examples}
\newtheorem{lemma}[equation]{Lemma}
\newtheorem{sublemma}[equation]{Sublemma}
\newtheorem{remar}[equation]{Remark}
\newtheorem{remars}[equation]{Remarks}
\newtheorem{nota}[equation]{Notation}
\newtheorem{sremar}[equation]{Side-remark}
\newtheorem{definitio}[equation]{Definition}

\newenvironment{remark}{\begin{remar} \rm }{\end{remar}}

\newenvironment{example}{\begin{exa} \rm }{\end{exa}}
\newenvironment{definition}{\begin{definitio} \rm }{\end{definitio}}

\newcommand{\HH}{{\mathbb{H}}}
\newcommand{\KK}{{\bf K}}

\newcommand{\CA}{{\cal A}}

\newcommand{\CE}{{\cal E}}

\newcommand{\ZZ}{\mathbb{Z}}
\newcommand{\RR}{\mathbb{R}}
\newcommand{\QQ}{\mathbb{Q}}
\newcommand{\CC}{\mathbb{C}}
\newcommand{\NN}{\mathbb{N}}
\newcommand{\bp}{\noindent {\sc Proof: }}
\newcommand{\eop}{\nopagebreak
			\hspace*{\fill}{$\diamond$}
			\medskip}

\newcommand{\tata}{\begin{pspicture}[0.2](0,0)(.5,.4)
\pscircle(0.25,0.2){.25}
\psline{*-*}(0.05,.2)(.45,.2)
\end{pspicture}}

\begin{document}
\thispagestyle{empty}
\begin{center}
\vskip4cm

\huge{On the Kontsevich-Kuperberg-Thurston construction \\of a configuration-space invariant \\for rational homology $3$-spheres}

\vskip1cm

\large{Christine Lescop}

\vskip1cm

\begin{large}
Pr\'epublication de l'Institut Fourier n$^{o}$ 655 (2004) \\
http://www-fourier.ujf-grenoble.fr/prepublications.html \\
\end{large}
\normalsize(Second version: New title, minor modifications in the abstract and in the introduction) 
\end{center}

\vskip 4cm

\begin{abstract}
M.~Kontsevich proposed a topological construction for an invariant $Z$ of rational homology $3$-spheres using configuration space integrals. 
G.~Kuperberg and D.~Thurston proved that $Z$ is a universal real finite type invariant for integral homology spheres in the sense of Ohtsuki, Habiro and Goussarov.

We review the Kontsevich-Kuperberg-Thurston construction and we provide detailed and elementary proofs for the invariance of $Z$.
This article is the preliminary part of a work that aims to prove splitting formulae for this powerful invariant of rational homology spheres.
It contains the needed background for the proof that will appear in the second part.

\vskip1cm

{\bf A.M.S. subject classification:} 57M27 57N10 55R80 57R20

\end{abstract}

\vskip4cm

\noindent Institut Fourier \\
Laboratoire de math\'ematiques de l'Universit\'e Joseph Fourier, UMR 5582 du CNRS \\
B.P.74 \\
38402 Saint-Martin-d'H\`eres cedex (France)\\
mail: lescop@ujf-grenoble.fr

\baselineskip.5cm
\newpage
\tableofcontents

\newpage
\section{The statements}
\label{secsta}
\setcounter{equation}{0}

We review the Kontsevich-Kuperberg-Thurston construction of an invariant $Z$ of rational homology spheres in \cite{kt,ko}. (See \cite{as1,as2,bc1,bc2,cat} for another construction.) This invariant is constructed by means of configuration space integrals, it is valued in the algebra $\CA(\emptyset)$ of Jacobi diagrams.
Its main property, that was proved by Greg Kuperberg and Dylan Thurston, is that it is a universal real finite type invariant for homology spheres in the sense of 
\cite{ggp,hab,oht}. A generalization of this property is proved in \cite{sumgen}. Here, we provide detailed
and elementary proofs for the invariance of $Z$, and for the
properties of $Z$ that are needed in \cite{sumgen}.

All the main ideas here are due to Witten, Axelrod, Singer, Kontsevich, Bott, Taubes, Cattaneo, G.~Kuperberg and D.~Thurston among others. 
I thank Dylan Thurston for explaining them to me.
  
The invariant $Z$ is a powerful generalization of the Casson invariant for integral homology $3$-spheres.
In this setting, the {\em Casson invariant\/} normalized as in \cite{akmc,mar} may be described as \index{N}{lambda@$\lambda$}
$$\lambda(M)=\frac{1}{6}\int_{C_2(M)}\omega_M^3$$
for any 2-form $\omega_M$ that satisfies the hypotheses stated in Subsections~\ref{subfundform}, \ref{subapseckkt} and \ref{subpont}
on the configuration space $C_2(M)$ defined in Subsection~\ref{subconfap}.

\subsection{The configuration space $C_2(M)$}
\label{subconfap}
\index{N}{Ctwo@$C_2(M)$}

When $A$ is a subset of $B$, $(B \setminus A)$ denotes its complement in $B$, when $x \in B$, $(B \setminus \{x\})$ is also denoted by $(B \setminus x)$.

Let $X=X^d$ be a smooth $d$-dimensional (real) manifold, and let $Y^k$ be a smooth $k$-dimensional submanifold of $X$.
If $TY$ denotes (the total space of) the tangent bundle of $Y$, then
$TX/TY$ is the {\em normal bundle\/} of $Y$. When $V$ is a real vector space, the group $]0,\infty[$ acts on $V$ by multiplication, and we set
$$SV=S(V)=(V \setminus 0)/]0,\infty[.$$ 
\index{N}{SV@$S(V)$}
The {\em unit normal bundle\/} $SN_XY$ of $Y$ in $X$ is the bundle over $Y$ whose fiber over $y$ is $S(T_yX/T_yY)$.
In this article, {\em to blow-up\/} a submanifold $Y^k$ in $X^d$ amounts to replace $Y$ by its unit normal bundle.

For example, if $Y \times \RR^{d-k}$ is a tubular neighborhood of $Y=Y \times \{0\}$ in $X$, the blow-up is equivalent to the sequence of operations:
$$Y \times \RR^{d-k} \longrightarrow [(Y \times \RR^{d-k}) \setminus Y]=Y \times
]0,\infty[ \times S^{d-k-1} \longrightarrow Y \times
[0,\infty[ \times S^{d-k-1}.$$ 
In general, this provides a local definition. See Definition~\ref{defblodif}
for the general definition.
The blown-up manifold inherits a smooth structure of a manifold with corners
from the smooth structure of $X$. See Proposition~\ref{propblodifun}.
Note that the blown-up manifold has the homotopy type of $(X \setminus Y)$.

Let $M$ be a closed oriented 3-manifold. Fix $\infty \in M$. 

Let $C_1(M)$ 
\index{N}{Cone@$C_1(M)$}
denote the manifold obtained from $M$ by blowing-up $\infty$.
The boundary of $C_1(M)$ is $ST_{\infty}(M)$. It is homeomorphic to $S^2$.
Let $M^2(\infty,\infty)$ 
\index{N}{Mtwoinfty@$M^2(\infty,\infty)$}
denote the manifold obtained from $M^2$ by blowing up
$(\infty,\infty)$ that becomes $S(T_{(\infty,\infty)}M^2) \cong S^5$.
In $M^2(\infty,\infty)$, the closures of the three submanifolds
of $M^2 \setminus (\infty,\infty)$, $\infty \times (M \setminus \infty)$, 
$(M \setminus \infty) \times \infty$ and $\mbox{diag}((M \setminus \infty)^2)$
are three disjoint submanifolds of $M^2(\infty,\infty)$ which intersect
$S(T_{(\infty,\infty)}M^2)$ along $S(0 \times T_{\infty}M)$, $S( T_{\infty}M \times 0)$ and $S(\mbox{diag}((T_{\infty}M)^2))$, respectively.
The three of them are canonically diffeomorphic to $C_1(M)$. They will be denoted 
by $\infty \times C_1(M)$, 
$C_1(M) \times \infty$ and $\mbox{diag}(C_1(M)^2)$, respectively.

The normal bundle of $\mbox{diag}((M \setminus \infty)^2)$ in $(M \setminus \infty)^2$ is identified with the tangent bundle of $(M \setminus \infty)$
through $$(u,v) \in (T_xM)^2/\mbox{diag}((T_xM)^2) \mapsto (v-u) \in T_xM.$$

The {\em configuration space\/} $C_2(M)$ is 
the {\em compactification\/} of  $$\breve{C}_2(M)=(M \setminus \infty)^2 \setminus \mbox{diagonal}$$
obtained from $M^2(\infty,\infty)$
by blowing-up $\infty \times C_1(M)$, 
$C_1(M) \times \infty$ and $\mbox{diag}(C_1(M)^2)$.

\subsection{Fundamental forms on $C_2(M)$}
\label{subfundform}

For $S^3=\RR^3 \cup \infty$, we have a homotopy equivalence $p_{S^3}$ that makes
the following square commute:
$$\diagram{\breve{C}_2(S^3)&\hfl{p_{S^3}}&S^2 \cr\vfl{=}&&\uvfl{\mbox{projection}}
\cr (\RR^3)^2 \setminus \mbox{diag}&\hfl{\cong}&\RR^3 \times ]0,\infty[ \times S^2 \cr (x,y)& \mapsto & (x,\parallel y-x \parallel , \frac{y-x }{\parallel y-x \parallel})}.
$$

The following lemma is proved at the end of Subsection~\ref{subsdifblowup}.
The natural projection onto the $X$-factor of a product is denoted by $\pi_X$.

\begin{lemma}
\label{lemextproj}
The map $p_{S^3}$ smoothly extends to $C_2(S^3)$, and its extension $p_{S^3}$
satisfies:
$$p_{S^3}=\left\{\begin{array}{ll} -\pi_{S^2} \;\;& \mbox{on} \;ST_{\infty}(S^3) \times (S^3 \setminus \infty)=S^2 \times (S^3 \setminus \infty)\\ \pi_{S^2} \;\;& \mbox{on} \; (S^3 \setminus \infty) \times ST_{\infty}(S^3)=(S^3 \setminus \infty) \times S^2\\ \pi_{S^2} \;\;& \mbox{on} \;ST(\RR^3) 
{=} \RR^3 \times S^2\end{array}\right.$$
\end{lemma}

Let $B^3(r)$ \index{N}{bt@$B^3(r)$} be the ball of $\RR^3$ centered at $0$ with radius $r$.
Let $\phi$ be an orientation-preserving embedding of $(S^3 \setminus \mbox{Int}(B^3(1)))$ into $M$.
Then
 $$M=\left(\RR^3 \cup \infty \setminus \mbox{Int}(B^3(1))\right) \cup_{ ]1,3] \times S^2} B_M,$$ where $$B_M=M \setminus \phi\left(S^3 \setminus \mbox{Int}(B^3(3))\right),$$ and $ (]1,3]\times S^2 = \phi(]1,3]\times S^2))$ is a collar of $\partial B_M$ in $B_M$. Fix $(\infty \in M)=\phi(\infty)$.
This identifies $ST_{\infty}M$ to $(ST_{\infty}(S^3)=S(\RR^3)=S^2)$.

\begin{definition}
A {\em trivialisation of \/$T(M \setminus \infty)$ that is standard near $\infty$\/}
is a trivialization \index{T}{trivialisation standard near $\infty$}
$$\tau: T(M \setminus \infty) \longrightarrow (M \setminus \infty) \times \RR^3$$
of the tangent bundle  $T(M \setminus \infty)$ of $(M \setminus \infty)$ that agrees with the standard trivialization $\tau_{S^3}$ of $\RR^3$ outside $B_M(1)= B_M \setminus ( ]1,3] \times S^2 )$.
\end{definition}
Let  $\tau_M$ be such a trivialisation (that exists by Lemma~\ref{lemtrivexist}). Note that $\tau_M$ identifies $S(T(M \setminus \infty))$ to $(M \setminus \infty) \times S^2$.

\begin{remark}
In the sequel, we define an invariant under orientation-preserving diffeomorphisms for pairs $(M, \phi)$ where
$M$ is an oriented $\QQ$-sphere, and $\phi$ is an orientation-preserving embedding of $(S^3 \setminus \mbox{Int}(B^3(1)))$ into $M$. Since all such embeddings $\phi$
are isotopic in $M$, the choice of $\phi$ will not matter. This allows us to fix our decomposition
$$M=\left(\RR^3 \cup \infty \setminus \mbox{Int}(B^3(1))\right) \cup_{ ]1,3] \times S^2} B_M,$$ 
once for all. This will not be discussed anymore.
\end{remark}

Let 
$$P:C_2(M) \longrightarrow M^2$$
be the natural projection map.
The identification of $M$
and $S^3$ in a neighborhood of $\infty$ provides identifications of neighborhoods of $P^{-1}(\infty,\infty)$ in $\partial C_2(S^3)$ and in $\partial C_2(M)$.

Define $p_M(\tau_M): \partial C_2(M) \longrightarrow S^2$ 
\index{N}{pM@$p_M(\tau_M)$}
by carrying the definition of  $p_{S^3}$ in the neighborhood above and by mimicking the definition of $p_{S^3}$ elsewhere on $\partial C_2(M)$. Recall that $ST_{\infty}(M)=S^2$
$$p_M(\tau_M)=\left\{\begin{array}{ll} -\pi_{S^2} \;\;& \mbox{on} \;ST_{\infty}(M) \times (M \setminus \infty)=S^2 \times (M \setminus \infty)\\ \pi_{S^2} \;\;& \mbox{on} \; (M \setminus \infty) \times ST_{\infty}(M)=(M \setminus \infty) \times S^2\\ \pi_{S^2}(\tau_M) \;\;& \mbox{on} \;ST(M \setminus \infty) \stackrel{\tau_M}{=} (M \setminus \infty) \times S^2\end{array}\right.$$

Let $\iota$ 
\index{N}{iota@$\iota$} 
be the involution of $C_2(M)$ that extends ($(x,y) \mapsto (y,x)$)
and let $\overline{\iota}$ 
\index{N}{iotab@$\overline{\iota}$} 
be the antipode of $S^2$.\\
Let $\omega_{S^2}$ be a volume form on $S^2$ such that 
$\int_{S^2}\omega_{S^2}=1$.
We say that  $\omega_{S^2}$ is {\em antisymmetric\/} if $\overline{\iota}^{\ast}(\omega_{S^2})=-\omega_{S^2}$.

Let $\tau_M$ be a trivialisation of $T(M \setminus \infty)$ that is standard near $\infty$.
\begin{definition}
\label{defformfund}
\index{T}{form!fundamental}
\noindent A two-form $\omega_M$ on $C_2(M)$ is {\em fundamental with respect to \/$\tau_M$ and \/} $\omega_{S^2}$ if:\\
$\bullet$ its restriction to $\partial C_2(M)$ is $p_M(\tau_M)^{\ast}(\omega_{S^2})$, and,\\ 
$\bullet$ it is closed.\\
Such a two-form is  {\em antisymmetric\/} if  
\index{T}{form!antisymmetric}
$\iota^{\ast}(\omega_M)=-\omega_M$. 
\end{definition}

It will be easily shown (Lemma~\ref{lemexisfunap}) that such forms exist for any trivialization $\tau_M$ when $M$ is a $\QQ$-sphere.

\subsection{Jacobi diagrams}
Here, a {\em  \indexT{Jacobi diagram}\/} $\Gamma$ is a trivalent graph $\Gamma$ without simple loop like $\begin{pspicture}[.2](0,0)(.6,.4)
\psline{-*}(0.05,.2)(.25,.2)
\pscurve{-}(.25,.2)(.4,.05)(.55,.2)(.4,.35)(.25,.2)
\end{pspicture}$. The set of vertices of such a $\Gamma$ will be denoted by $V(\Gamma)$,
\index{N}{VGamma@$V(\Gamma)$}
its set of edges will be denoted by $E(\Gamma)$.
\index{N}{EGamma@$E(\Gamma)$}
A {\em \indexT{half-edge}\/} $c$ of $\Gamma$ is an element of
$$H(\Gamma)=\{c=(v(c);e(c)) | v(c) \in V(\Gamma); e(c) \in E(\Gamma);v(c) \in e(c)\}.$$
\index{N}{HGamma@$H(\Gamma)$}
An {\em automorphism\/} of $\Gamma$ 
\index{T}{Jacobi diagram!automorphism of} 
is a permutation $b$ of $H(\Gamma)$
such that for any $c,c^{\prime} \in H(\Gamma)$,
$$v(c)=v(c^{\prime}) \Longrightarrow v(b(c))=v(b(c^{\prime}))\;\;\mbox{and}\;\;e(c)=e(c^{\prime}) \Longrightarrow e(b(c))=e(b(c^{\prime})).$$
The number of automorphisms of $\Gamma$ will be
denoted by $\sharp \mbox{Aut}(\Gamma)$. 
\index{N}{AutGamma@$\sharp \mbox{Aut}(\Gamma)$}
For example, $ \sharp \mbox{Aut}(\tata)=12$.
{\em An orientation\/} of a vertex of such a diagram $\Gamma$ 
 is a cyclic order of the three
half-edges that meet at that vertex.
\index{T}{Jacobi diagram!orientation of}
A Jacobi diagram $\Gamma$ is {\em oriented\/} if all its vertices are oriented (equipped with an orientation).
The {\em degree} of such a diagram is 
half the number of its vertices. 

Let $\CA_n(\emptyset)$ 
\index{N}{An@$\CA_n(\emptyset)$}
denote the real vector space generated by the degree $n$ oriented Jacobi diagrams, quotiented out by the following relations AS and IHX:

$$ {\rm AS :}  \begin{pspicture}[.2](0,-.2)(.8,1)
\psset{xunit=.7cm,yunit=.7cm}
\psarc[linewidth=.5pt](.5,.5){.2}{-70}{15}
\psarc[linewidth=.5pt](.5,.5){.2}{70}{110}
\psarc[linewidth=.5pt]{->}(.5,.5){.2}{165}{250}
\psline{*-}(.5,.5)(.5,0)
\psline{-}(.1,.9)(.5,.5)
\psline{-}(.9,.9)(.5,.5)
\end{pspicture}
+
\begin{pspicture}[.2](0,-.2)(.8,1)
\psset{xunit=.7cm,yunit=.7cm}
\pscurve{-}(.9,.9)(.3,.7)(.5,.5)
\pscurve[border=2pt]{-}(.1,.9)(.7,.7)(.5,.5)
\psline{*-}(.5,.5)(.5,0)
\end{pspicture}=0,\;\;\mbox{and IHX :} 
\begin{pspicture}[.2](0,-.2)(.8,1)
\psset{xunit=.7cm,yunit=.7cm}
\psline{-*}(.1,1)(.35,.2)
\psline{*-}(.5,.5)(.5,1)
\psline{-}(.75,0)(.5,.5)
\psline{-}(.25,0)(.5,.5)
\end{pspicture}
+
\begin{pspicture}[.2](0,-.2)(.8,1)
\psset{xunit=.7cm,yunit=.7cm}
\psline{*-}(.5,.6)(.5,1)
\psline{-}(.8,0)(.5,.6)
\psline{-}(.2,0)(.5,.6)
\pscurve[border=2pt]{-*}(.1,1)(.3,.3)(.7,.2)
\end{pspicture}
+
\begin{pspicture}[.2](0,-.2)(.8,1)
\psset{xunit=.7cm,yunit=.7cm}
\psline{*-}(.5,.35)(.5,1)
\psline{-}(.75,0)(.5,.35)
\psline{-}(.25,0)(.5,.35)
\pscurve[border=2pt]{-*}(.1,1)(.2,.75)(.7,.75)(.5,.85)
\end{pspicture}
=0. 
$$
\index{N}{AS}
\index{N}{IHX}

Each of these relations relate diagrams which can be represented by planar immersions that are identical outside the part of them represented in the pictures. Here, the orientation of vertices is induced by the counterclockwise order of the half-edges. For example, 
AS identifies the sum of two diagrams which only differ by the orientation
at one vertex to zero.
$\CA_0(\emptyset)$ is equal to $\RR$ generated by the empty diagram.

\subsection{The invariants $Z_n$}
\label{subapseckkt}

\begin{definition}
Let $V$ be a finite set. An {\em orientation\/} of $V$
\index{T}{orientation of a finite set} 
is a bijection from $V$ to $\{1,2,\dots,\sharp V\}$ (or a total order on $V$) up to an even permutation. 
\end{definition}
When $M$ is an odd-dimensional oriented manifold, an orientation of $V$ provides an ordering of the factors of $M^V$ (up to an even permutation), and therefore induces an orientation of $M^V$. Thus, the datum of an orientation of $V$ is equivalent to the datum of an orientation of $M^V$.

Let $\Gamma$ be a Jacobi diagram. Let $H(\Gamma)$ be its set of half-edges.
When the edges of $\Gamma$ are oriented, the orientations of the edges
induce an orientation of $H(\Gamma)$ that is called the {\em \indexT{edge-orientation}\/} of
$H(\Gamma)$ and that is represented by a total order of $H(\Gamma)$ of the following form. Fix an arbitrary order on the set of edges, then take the two halves of the first edge ordered from origin to the end, next the two halves of the second edge, and so on.

When the set $V(\Gamma)$ of vertices of $\Gamma$ is oriented and when the vertices of $\Gamma$ are oriented (as the sets of their three half-edges),
these data induce an orientation of $H(\Gamma)$ that is called the {\em \indexT{vertex-orientation}\/} of
$H(\Gamma)$ and that is defined as follows.
Number the vertices of $\Gamma$ from $1$ to $\sharp V(\Gamma)$ by a bijection that induces 
the given orientation of $V(\Gamma)$.
The wanted order of $H(\Gamma)$ is given by taking first the half-edges of the first vertex 
with an order that agrees with the vertex-orientation, then the half-edges that contain the second vertex, and so on.
\\
Let $M$ be a $\QQ$-sphere. Set $$\breve{C}_{V(\Gamma)}(M)=(M \setminus \infty)^{V(\Gamma)} \setminus \{\mbox{all diagonals}\}.$$
\index{N}{CbreveV@$\breve{C}_{V(\Gamma)}(M)$}
The set $\breve{C}_{V(\Gamma)}(M)$ is the set of injective maps from $V(\Gamma)$ to $(M\setminus \infty)$. It is an open submanifold of $(M \setminus \infty)^{V(\Gamma)}$ that is oriented as soon as $V(\Gamma)$ is oriented.

An edge $e$ of $\Gamma$ defines a pair $P(e)$ of elements of $V(\Gamma)$.
Then the restriction of maps induces a canonical map from $\breve{C}_{V(\Gamma)}(M)$
to $\breve{C}_{P(e)}(M)$. An orientation of $e$ orders the pair $P(e)$
and produces a canonical identification of $\breve{C}_{P(e)}(M)$ with $\breve{C}_{\{1,2\}}(M) \subset C_2(M)$. (The origin of $e$ is mapped to $1$.)
For any oriented edge $e$ of $\Gamma$, the composition of these maps will be denoted by $$p_e:\breve{C}_{V(\Gamma)}(M)\longrightarrow C_2(M).$$
\index{N}{pe@$p_e$}
Let $\omega_M$ be an antisymmetric two-form that is fundamental with respect to $\tau_M$ and $\omega_{S^2}$. 

Let $\Gamma$ be an oriented Jacobi diagram.
Orient the edges of $\Gamma$, and orient $V(\Gamma)$ so that the
edge-orientation of $H(\Gamma)$ coincides with the vertex-orientation
of $H(\Gamma)$.
Set $$I_{\Gamma}(\omega_M)=\int_{\breve{C}_{V(\Gamma)}(M)}\bigwedge_{e \in E(\Gamma)}p_e^{\ast}(\omega_M).$$
\index{N}{IGammaomegaM@$I_{\Gamma}(\omega_M)$}
This integral is convergent thanks to Proposition~\ref{propconfunc} below. It is easy to see that its sign only depends on the vertex-orientation of $\Gamma$ up to an
even number of changes. In particular, the product $I_{\Gamma}(\omega_M)[\Gamma]$ only depends on the (unoriented) Jacobi diagram $\Gamma$.

\begin{proposition}[\cite{kt}]
\label{propthkktun}
Let $M$ be a $\QQ$-sphere. Let $\omega_M$ be an antisymmetric two-form that is fundamental with respect to a trivialization $\tau_M$ standard near $\infty$ and to a form $\omega_{S^2}$ such that $\int_{S^2}\omega_{S^2}=1$. Then with the notation above \index{N}{ZnMtauM@$Z_n(M;\tau_M)$}
$$Z_n(M;\tau_M)=\sum_{\Gamma \;\mbox{\small Jacobi diagram with $2n$ vertices}}\frac{I_{\Gamma}(\omega_M)}{\sharp \mbox{Aut}(\Gamma)}[\Gamma] \in \CA_n(\emptyset)$$
only depends on the oriented diffeomorphism type of $M$ and on the homotopy class of $\tau_M$.
(Here the sum runs over Jacobi diagrams without vertex-orientations.)
\end{proposition}

In the next subsection, we shall see that any $\ZZ$-sphere $M$ has a preferred homotopy class $[\tau^0_M]$ of trivialisations that are standard near $\infty$, and this
will allow us to define the invariants of $\ZZ$-spheres by
$$Z_n(M)=Z_n(M;\tau^0_M).$$
\index{N}{ZnM@$Z_n(M)$}
In general, we shall need a correction term, called the {\em framing correction\/} that is described in Subsection~\ref{subfra}.

\subsection{Homotopy classes of trivialisations of $\QQ$-spheres.}
\label{subpont}

Recall that $GL^+(\RR^3)$ is homotopy equivalent to the pathwise connected group $SO(3)$, that $\pi_1(SO(3))\cong \ZZ/2\ZZ$, $\pi_2(SO(3)) \cong 0$ and
$\pi_3(SO(3)) \cong \ZZ[\rho]$ where the generator $ [\rho]$ of $\pi_3(SO(3))$
is represented by the following covering map
$$\rho:S^3 \longrightarrow SO(3).$$
\index{N}{rho@$\rho$}
See $S^3$ as the unit sphere of the quaternionic field $(\HH=\RR \oplus \RR i \oplus \RR j \oplus \RR k)$. 
\index{N}{H@$\HH$}
Then, for any
element $\gamma$ of $S^3$, $\rho(\gamma)$ is the restriction of the conjugacy $(x \mapsto \gamma x\gamma^{-1})$ to the euclidean space $\RR^3$ of the pure quaternions.

Boundaries of oriented manifolds are oriented 
with the outward normal first convention. Unit spheres of oriented euclidean
vector spaces are oriented as the boundaries of unit balls. In particular,
the sphere $S^3$ is the oriented boundary of the unit ball of $\HH$.
The group $SO(3)$ is locally oriented as $S^2 \times S^1$ (oriented rotation axis in $S^2$, rotation angle with respect to the previous axis) (outside its center). With these orientations, $\mbox{deg}(\rho)=2$.
\index{N}{rho@$\rho$} 

\begin{lemma}
\label{lemtrivexist}
The trivialisation $\tau_M$ defined on $]1,3] \times S^2$
extends to $B_M$.
\end{lemma}
\bp Choose a cell decomposition of $B_M$ with respect to its boundary.
Since $GL^+(\RR^3)$ is pathwise connected, we may extend the trivialisation to the one-skeleton of $B_M$. If there were an obstruction in 
$$H^2(B_M,\partial B_M=S^2;\ZZ/2\ZZ)=H^2(B_M;\ZZ/2\ZZ)$$ to extend $\tau_M$
on the two-skeleton of $B_M$, there would exist a surface $\Sigma$ immersed in $M$ such that the pull-back of $TM$ under this immersion is not trivialisable on $\Sigma$. But this
pull-back is isomorphic to the sum of the tangent space $T\Sigma$ of $\Sigma$, and the unique one-dimensional fibered bundle $\eta$ over $\Sigma$ that makes $T\Sigma \oplus \eta$ orientable. Therefore, the pull-back of $TM$ under this immersion is isomorphic to the pull-back of $T\RR^3$ under any immersion of $\Sigma$ into $\RR^3$ and is trivialisable on $\Sigma$. Thus, $\tau_M$ extends to the two-skeleton of $\Sigma$. Since $\pi_2(GL^+(\RR^3)) = 0$, $\tau_M$ also extends to the three-skeleton. 
\eop

The above proof also shows that any oriented 3-manifold is parallelisable.

Recall that the signature of a $4$-manifold is the signature of the intersection form on its $H_2$. Also recall that any closed oriented three-manifold bounds a compact oriented $4$-dimensional manifold whose signature may be arbitrarily changed by connected
sums with copies of $\CC P^2$ or $-\CC P^2$.
Let $W=W^4$ be a signature $0$ cobordism between $B^3(3)$ and $B_M$, that is a compact oriented $4$-dimensional manifold with corners such that
$$\partial W= B_M \cup (-[0,1] \times S^2) \cup -B^3(3)$$
where $\partial B_M= \partial B^3(3)=S^2$.

$$
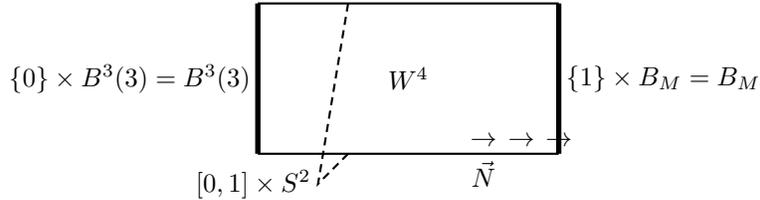
\begin{pspicture}[.4](-3,-.7)(5.5,2.5)
\psline(0,0)(4,0)
\psline(4,2)(0,2)
\psline[linewidth=2pt](0,0)(0,2)
\psline[linewidth=2pt](4,0)(4,2)
\rput(2,1){$W^4$}
\rput[r](-.1,1){$\{0\} \times B^3(3)=B^3(3) $}
\rput[l](4.1,1){$\{1\} \times B_M=B_M$}
\psline[linestyle=dashed,dash=3pt 2pt](1.2,0)(.8,-.4)(1.2,2)
\rput[r](.7,-.4){$[0,1] \times S^2$}
\rput[b](4,.1){$ \rightarrow $}
\rput[b](3.5,.1){$ \rightarrow $}
\rput[b](3,.1){$ \rightarrow $}
\rput[t](3,-.1){$ \vec{N}$}
\end{pspicture}$$

Let $\tau_M$ be a trivialisation of $(M \setminus \infty)$ that is standard near $\infty$. 
Define the {\em  \indexT{Pontryagin number}\/} of $\tau_M$ $$p_1(\tau_M) \in \ZZ$$ \index{N}{pone@$p_1$} as follows.

Consider the complex $4$-bundle $TW \otimes \CC$ over $W$. 
Near $\partial W$, $W$ may be identified to an open subspace of one of the products $[0,1] \times B^3(3)$ or $[0,1] \times B_M$. 
Let $\vec{N}$ be the tangent vector to $[0,1] \times \{\mbox{pt}\}$ (under these identifications), and let 
$\tau(\tau_M)$ denote the trivialization of $TW \otimes \CC$ over $\partial W$ that is obtained by stabilizing either $\tau_{S^3}$ or $\tau_M$ into $\vec{N} \oplus \tau_M$ or  $\vec{N} \oplus \tau_{S^3}$. Then the obstruction to extend this trivialization to $W$ is the relative first \indexT{Pontryagin class} $$p_1(W;\tau(\tau_M))=p_1(\tau_M)[W,\partial W] \in H^4(W,\partial W;\ZZ)=\ZZ[W,\partial W]$$ of the trivialisation. 

Now, we specify our sign conventions for this Pontryagin class. They are the same as in \cite{milnorsta}. 
In particular, $p_1$ is the opposite of the second Chern class $c_2$ of the complexified tangent bundle. See \cite[p. 174]{milnorsta}.
More precisely, equip $M$ with a riemannian metric that coincides with the standard metric
of $\RR^3$ outside $B^3(1)$, and equip $W$ with a riemannian metric that coincides with the orthogonal product metric of one of the products $[0,1] \times B^3(3)$ or $[0,1] \times B_M$ near $\partial W$. 
Equip $TW \otimes \CC$ with the associated hermitian structure. The determinant bundle of $TW$
is trivial because $W$ is oriented and $det(TW \otimes \CC)$ is also trivial.
We only consider the trivialisations that
are unitary with respect to the hermitian structure of $TW \otimes \CC$ and the standard hermitian form of $\CC^4$, and that are special with respect to
the trivialisation of $det(TW \otimes \CC)$.
Since $\pi_i(SU(4))=\{0\}$ when $i<3$, the trivialisation $\tau(\tau_M)$
extends to a special unitary trivialisation $\tau$ outside the interior of a $4$-ball $B^4$
and defines 
$$\tau: (TW \otimes \CC)_{|S^3} \longrightarrow S^3 \times \CC^4$$
over the boundary $S^3=\partial B^4$ of this $4$-ball $B^4$.
Over this $4$-ball $B^4$, the bundle is trivial and admits a trivialisation
$$\tau_B: (TW \otimes \CC)_{|B^4} \longrightarrow B^4 \times \CC^4.$$
Then $\tau_B \circ \tau^{-1}(v \in S^3, w \in \CC^4)=(v, \phi(v)(w))$
where $\phi(v) \in SU(4)$. 
Let $i^2(m^{\CC}_r)$ 
\index{N}{itwo@$i^2(m^{\CC}_r)$}
be the following map
$$\begin{array}{llll} i^2(m^{\CC}_r): &(S^3 \subset \CC^2) & \longrightarrow & SU(4)\\
& (z_1,z_2) & \mapsto & \left[\begin{array}{cccc} 1&0&0&0\\0&1&0&0\\0&0&z_1&-\overline{z}_2\\0&0&z_2& \overline{z}_1\end{array} \right] \end{array}.$$
When $(e_1,e_2,e_3,e_4)$ is the standard basis of $\CC^4$, the columns of the matrix contain the coordinates of the images of the $e_i$ with respect to $(e_1,e_2,e_3,e_4)$.
Then the homotopy class $[i^2(m^{\CC}_r)]$ of $i^2(m^{\CC}_r)$ generates $\pi_3(SU(4))=\ZZ[i^2(m^{\CC}_r)]$
and the homotopy class of $\phi: S^3 \longrightarrow  SU(4)$ satisfies
$$[\phi]=-p_1(\tau_M)[i^2(m^{\CC}_r)] \in \pi_3(SU(4)).$$

\begin{proposition}
\label{proppont}
The first \indexT{Pontryagin number} $p_1(\tau_M)$ \index{N}{pone@$p_1$} is well-defined by the above 
conditions. (It is independent of the choices that were made.)
It only depends on the homotopy class of the trivialisation $\tau_M$
among the trivialisations that are standard near $\infty$.

For any closed 3-manifold $M$, for any 
trivialisation $\tau_M$ of $T(M \setminus \infty)$ that is standard near $\infty$, and for any $$g: (B_M,  ]1,3] \times S^2) \longrightarrow (SO(3),1),$$
let $\mbox{deg}(g)$ denote the degree of $g$ and let 
$$\begin{array}{llll} 
\psi(g): &B_M \times \RR^3 &\longrightarrow  &B_M \times \RR^3\\
&(x,y) & \mapsto &(x,g(x)(y))\end{array}$$
\index{N}{psig@$\psi(g)$} 
then 
$$p_1(\psi(g) \circ \tau_M)-p_1(\tau_M)=-2\mbox{deg}(g).$$ 
If $M$ is a given $\ZZ$-sphere, then $p_1$ defines a bijection from the set of homotopy classes of trivialisations of $M$ that are standard near $\infty$ to $4 \ZZ$.
\end{proposition}

This proposition will be proved in Subsection~\ref{subproofpont}.
Of course, for a given $\ZZ$-sphere, our preferred class of trivialisations will be $p_1^{-1}(0)$. By definition, the standard trivialisation of $\RR^3$ is in this class when $M=S^3$.

\subsection{The framing correction}
\label{subfra}

Let $X$ be a 3-dimensional vector space. Let $V$ be a finite set. Then
$\breve{S}_V(X)$ 
\index{N}{SbreveV@$\breve{S}_V(X)$}
denotes the set of injective maps from $V$ to $X$ up to translations and dilations. It is an open subset of the smooth manifold $S(X^V/
\mbox{diag}(X^V))$. Set $\breve{S}_n(X)=\breve{S}_{\{1,2,\dots, n\}}(X)$.
\index{N}{Sbreven@$\breve{S}_n(X)$}

When $V$ and $X$ are oriented, $X^V$ and $(\mbox{diag}(X^V) \cong X)$ are oriented,
then the quotient $X^V/\mbox{diag}(X^V)$ is oriented so that $X^V$ has the
(fiber $\mbox{diag}(X^V)$ $\oplus$ quotient $X^V/\mbox{diag}(X^V)$) orientation. 
When $W$ is a vector space, $S(W)$ is oriented as the boundary of a unit ball
of $W$ equipped with an arbitrary norm, that is so that the multiplication from $]0,\infty[ \times S(W)$ to $W$
preserves the orientation. This orients $S(X^V/
\mbox{diag}(X^V))$ and hence $\breve{S}_V(X)$. 

For an $\RR^3$ vector bundle, $p:E \longrightarrow B$,  $\breve{S}_{V}(E)$ denotes the fibered space over $B$ where the fiber over $(g \in B)$ is $\breve{S}_{V}(p^{-1}(g))$.
When $B$ is an oriented manifold, $\breve{S}_{V}(E)$ is next oriented with the (base $B$ $\oplus$ fiber) orientation.

Let $p:E_1 \longrightarrow S^4$ be the $\RR^3$ vector bundle over $S^4=B^4 \cup_{S^3} (-B^4)$ whose total space is \index{N}{Eone@$E_1$}
$$E_1= B^4 \times \RR^3 \cup_{S^3 \times \RR^3} (-B^4) \times \RR^3$$ \index{N}{rho@$\rho$}
where the two parts are glued by identifying $(g,x) \in S^3 \times \RR^3$ of the first factor to $$\left((g,\rho(g)(x)) \in  (-B^4) \times \RR^3\right).$$
The vector bundle $E_1$ is equipped with the involutive bundle isomorphism  $\iota$
\index{N}{iota@$\iota$} 
over $\mbox{Id}_{S^4}$ that is the multiplication by $(-1)$ over each fiber. 


In particular, $\breve{S}_2(E_1)$ is a (compact) $S^2$-bundle that is denoted by $S_2(E_1)$. \index{N}{StwoEone@$S_2(E_1)$}
Let $\omega_T$ \index{N}{omegaT@$\omega_T$} be a closed $2$-form on $S_2(E_1)$ that represents the Thom class of this $S^2$-bundle such that $\overline{\iota}^{\ast}(\omega_T)=-\omega_T$. ($[\omega_T]$ is dual to 
a $4$-dimensional manifold that intersects the "left-hand side part" $B^4 \times S^2$ of $S_2(E_1)$ as $(B^4 \times \{\mbox{point}\})$.) 

Let $\Gamma$ be an oriented Jacobi diagram. 
Each edge $e$ of $\Gamma$ again defines a pair $P(e) \subset V(\Gamma)$
that induces  a projection 
$$\breve{S}_{V(\Gamma)}(E_1) \longrightarrow S_{P(e)}(E_1)$$
by restriction on the fibers. 
An orientation of the edge $e$ induces an order on $P(e)$ that identifies
$S_{P(e)}(E_1)$ to $S_2(E_1)$, and this again defines
$$p_e: \breve{S}_{V(\Gamma)}(E_1) \longrightarrow S_2(E_1).$$
\index{N}{pe@$p_e$}
Orient the vertices, the edges of $\Gamma$, and orient $V(\Gamma)$ so that the
edge-orientation of $H(\Gamma)$ coincides with the vertex-orientation
of $H(\Gamma)$.
Set \index{N}{IGammaomegaT@$I_{\Gamma}(\omega_T)$} \index{N}{Eone@$E_1$}
$$I_{\Gamma}(\omega_T)[\Gamma]=\int_{\breve{S}_{V(\Gamma)}(E_1)}\bigwedge_{e \; \mbox{\small edge of}\; \Gamma}p_e^{\ast}(\omega_T)[\Gamma],$$
and define
$$\xi_n=\sum_{\Gamma \;\mbox{\small connected Jacobi diagram with $2n$ vertices}}\frac{I_{\Gamma}(\omega_T)}{\sharp \mbox{Aut}(\Gamma)}[\Gamma] \in \CA_n.$$ \index{N}{ksin@$\xi_n$}
Define $$\CA(\emptyset)=\prod_{n \in \NN}\CA_n(\emptyset)$$ 
\index{N}{Aemptyset@$\CA(\emptyset)$}as the topological product 
of the vector spaces $\CA_n(\emptyset)$.
Set \index{N}{ZMtauM@$Z(M;\tau_M)$} $$Z(M;\tau_M)=(Z_n(M;\tau_M))_{n\in \NN} \in \CA(\emptyset)$$
where $Z_0(M;\tau_M)=1[\emptyset]$.
Similarly, with $\xi_0=0$, 
$$\xi=(\xi_n)_{n\in \NN}.$$\index{N}{ksi@$\xi$}
Equip $\CA(\emptyset)$ with the continuous product that maps 
two (classes of) graphs to (the class of) their disjoint union.
This product turns $\CA(\emptyset)$ into a commutative algebra.

\begin{theorem}[\cite{kt}]
\label{thkktfra}
The obtained $\xi_n$ \index{N}{ksin@$\xi_n$}\index{N}{omegaT@$\omega_T$} does not depend on the closed form $\omega_T$ that represents the Thom class of $S_2(E_1)$ such that $\iota^{\ast}(\omega_T)=-\omega_T$.
For any $\QQ$-sphere $M$, set
$$Z(M)=Z(M;\tau_M) \exp(\frac{p_1(\tau_M)}{4}\xi).$$  \index{N}{ksi@$\xi$}
\index{N}{ZM@$Z(M)$}
Then $Z$ is a topological invariant of $M$.
\end{theorem}

Note that  Theorem~\ref{thkktfra} obviously implies Proposition~\ref{propthkktun}. Thus, we are left with the proof of Theorem~\ref{thkktfra}.

Equip $\CA(\emptyset)$ with the involution that maps $(x_n \in \CA_n(\emptyset))_{n \in \NN}$ to $\overline{(x_n)_{n \in \NN}}=((-1)^nx_n)_{n \in \NN}$. 
Then, we have the following proposition.

\begin{proposition}
\label{proporrev}
For any integer $k$, $\xi_{2k}=0$.  \index{N}{ksi@$\xi$} \\
For any $\QQ$-sphere $M$, let $(-M)$ denotes the manifold obtained from $M$ by reversing its orientation, then \index{N}{ZM@$Z(M)$}
$$Z(-M)=\overline{Z(M)}.$$
\end{proposition}
\bp
The involution $\iota$ still makes sense on $S_{2n}(E_1)$, it
reverses the orientation and it commutes with the projections $p_e$. Therefore, $$\begin{array}{ll}I_{\Gamma}(\omega_T)[\Gamma]&=-\int_{\breve{S}_{V(\Gamma)}(E_1)}
\iota^{\ast}\left(\bigwedge_{e \; \mbox{\small edge of}\; \Gamma}p_e^{\ast}(\omega_T)\right)[\Gamma]\\
&=-\int_{\breve{S}_{V(\Gamma)}(E_1)}\left(\bigwedge_{e \; \mbox{\small edge of}\; \Gamma}p_e^{\ast}(-\omega_T)\right)[\Gamma]\\
&=(1)^{\sharp E(\Gamma)+1}I_{\Gamma}(\omega_T)[\Gamma].\end{array}$$
Since $3\sharp V(\Gamma)=2\sharp E(\Gamma)=12k$, when the degree of $\Gamma$ is $2k$, we conclude that $\xi_{2k}=0$.

Consider a trivialisation 
$\tau_M:T(M \setminus \infty) \rightarrow (M \setminus \infty) \times \RR^3$ 
of $M$ that is standard near $\infty$. Its composition $\tau_{-M}$ by $(\mbox{Id}_{M \setminus \infty} \times (-1)\mbox{Id}_{\RR^3})$ is a trivialisation of $T(-M \setminus \infty)$ that is standard near $\infty$, with respect to the composition of the previous embedding of $(S^3 \setminus B^3(1))$
into $(M \setminus B_M(1))$ by the multiplication by $(-1)$.
On $\partial C_2(M)$, $p_{-M}(\tau_{-M})=\overline{\iota} \circ p_M(\tau_M)$. Therefore if $\omega(\tau_M)$ is an antisymmetric form that is fundamental with respect to $\tau_{M}$, $\iota^{\ast}(\omega(\tau_M))=-\omega(\tau_M)$ is an antisymmetric form that is fundamental with respect to $\tau_{-M}$.
Since changing the orientation of $M$, does not change the orientation of $C_{2n}(M)$, we see as before that $Z_{n}(-M;\tau_{-M})=Z_{n}(-M;-\omega(\tau_M))=(-1)^nZ_n(M;\tau_M)$.
Therefore, we are left with the proof that $p_1(\tau_{-M})=-p_1(\tau_M)$.
In order to prove it, note that changing the orientation of the cobordism $W$
between $B^3$ and  $B_M$
tranforms it into a cobordism between $-B^3 \cong B^3$ and $(B_{-M}=-B_M)$.
Furthermore, changing the trivialisation by preserving its first vector and reversing the other ones
equips $B^3$ with its standard trivialisation.
The latter trivialisation extends to the complement of a $4$-ball $B^4$ as the composition of the previous one by the above symmetry. Therefore the induced change of basis on $\partial B^4$ is conjugate through this symmetry of the connected group $U(4)$, and hence homotopic. Since the orientation of $\partial B^4$ is the opposite to the one used in the computation of 
$p_1(\tau_M)$, $p_1(\tau_{-M})=-p_1(\tau_M)$. \eop

We shall also prove that $\xi_1=-\frac{1}{12}[\tata]$ in Proposition~\ref{propxiun}.

Dylan Thurston and Greg Kuperberg also proved that $Z$ is a universal finite type invariant of integral homology $3$-spheres, that $Z$ is multiplicative
under the connected sum of $3$-manifolds, and that \index{N}{lambda@$\lambda$}
$$Z_1(M)=\frac{\lambda(M)}2[\tata]$$
for any integral homology sphere $M$ where $\lambda$ denotes the {\em Casson invariant\/} normalized as in \cite{akmc,mar}.

The article \cite{sumgen} contains splitting formulae for $Z$ that generalize
the formulae used in the Thurston and Kuperberg proof of $Z$'s universality.
It also contains a proof that $Z_1(M)=\frac{\lambda_W(M)}4[\tata]$
\index{N}{lambdaW@$\lambda_W$}
for any rational homology sphere $M$ where $\lambda_W$ denotes the Walker
extension of the Casson invariant normalized as in \cite{wal}.
Since the current article has been written in order to provide the detailed background 
for \cite{sumgen}, the Thurston and Kuperberg proof of $Z$'s universality will not be discussed here. The multiplicativity of $Z$ under connected sum that is not 
needed in \cite{sumgen} is not proved here either. This article is only a partial
detailed presentation of the properties of $Z$ that were discovered by Dylan Thurston and Greg Kuperberg in \cite{kt}, or by Maxim Kontsevich.

\newpage
\section{Proof of Theorem~\ref{thkktfra}}
\setcounter{equation}{0}
\label{secprothkkt}

\subsection{More on the topology of $C_2(M)$.}
\index{N}{Ctwo@$C_2(M)$}

Since the map $p_{S^3}:C_2(S^3) \longrightarrow S^2$ is a homotopy equivalence, $C_2(S^3)$ has the homotopy type of $S^2$.

In general, $C_2(M)$ has the homotopy type of $\left[(M \setminus \infty)^2 \setminus \mbox{diagonal}\right]$. Indeed, it has the homotopy type of 
$$M^2(\infty,\infty) \setminus \left( (\infty \times C_1(M)) \cup 
(C_1(M) \times \infty) \cup \mbox{diag}(C_1(M))\right)$$ that has the homotopy
type of $\left[(M \setminus \infty)^2 \setminus \mbox{diagonal}\right]$.
Therefore, we have the following lemma.

\begin{lemma}
\label{lemhctwo}
Let $\Lambda= \ZZ$ or $\QQ$.
If  $M$ is a $\Lambda$-sphere, then
$$H_{\ast}(C_2(M);\Lambda)=H_{\ast}(S^2;\Lambda).$$
and if \/$[S(T_xM)]$ denotes the homology class of a fiber of $(ST(M \setminus \infty) \subset C_2(M))$, then $H_{2}(C_2(M);\Lambda)=\Lambda[S(T_xM)]$.
\end{lemma}
\bp
In this proof, the homology coefficients are in $\Lambda$. Since $(M \setminus \infty)$ has the homology of a point, the K\"unneth Formula implies that
 $(M \setminus \infty)^2$ has the homology of a point. Now, by excision,
$$H_{\ast}((M \setminus \infty)^2,(M \setminus \infty)^2 \setminus \mbox{diag}) \cong H_{\ast}((M \setminus \infty) \times \RR^3,(M \setminus \infty)\times (\RR^3 \setminus 0))$$
$$ \cong H_{\ast}( \RR^3, S^2) \cong \left\{\begin{array}{ll} \Lambda \;\;&\;\mbox{if} \;\ast =3,\\
0\;&\;\mbox{otherwise.} \end{array} \right.$$
Of course, $(M \setminus \infty) \times \RR^3$ denotes a tubular neighborhood of the diagonal in $(M \setminus \infty)^2$. Note that such a neighborhood can be easily obtained by integrating the vector fields given by a trivialisation of $T(M \setminus \infty)$ standard near $\infty$. With $(m, \lambda (v \in S^2))$, associate $(m,\gamma_{\lambda}(m,v))$ where $\gamma_{0}(m,v)=m$ and $\frac{\partial}{\partial t}(\gamma_{t}(m,v))(t_0)=\tau_M^{-1}((\gamma_{t_0}(m,v),v))$. When $\varepsilon$ is a small enough positive number, this defines an embedding of $(M \setminus \infty) \times (\{ x \in \RR^3; \norm{x}< \varepsilon\} \cong \RR^3)$.

Using the long exact sequence associated to the pair $((M \setminus \infty)^2,(M \setminus \infty)^2 \setminus \mbox{diag})$, we get that
$$H_{\ast}(C_2(M))=H_{\ast}(S^2)$$
and that $H_{2}(C_2(M);\Lambda)=\Lambda[S(T_xM)]$.
\eop

Therefore, there is a preferred generator $L_M$ 
\index{N}{LM@$L_M$} 
of $H^2(C_2(M);\QQ)$ such that when $B$ is a 3-ball embedded in $M$ equipped with the orientation of $M$ and when $x$ is a point in the interior of $B$, the evaluation of $L_M$ on the homology class of $(\{x\} \times \partial B) \subset C_2(M)$ is one.

If $(K_1 \sqcup K_2) \subset M \setminus \infty$ is a two-component link of $M$, then the evaluation of $L_M$ on the homology class of the torus 
$(K_1 \times K_2 \subset C_2(M))$ is the {\em \indexT{linking number}\/} 
of $K_1$ and $K_2$ in $M$ that is denoted by $\ell(K_1,K_2)$. Here, it will be our definition for the linking number.

Let us now prove the existence of fundamental forms. 

For this, we first recall the following standard consequence of the definition of the De Rham cohomology. 

\begin{lemma}
\label{lemdr}
Let $A$ be a compact submanifold of a compact manifold $B$, let $\omega_A$ be a closed $n$-form on $A$, and let $i:A \longrightarrow B$ denote the inclusion.
Then the three following assertions are equivalent:
\begin{enumerate}
\item The form $\omega_A$ extends to $B$ as a closed $n$-form. 
\item The cohomology class of $\omega_A$ belongs to $i^{\ast}(H^n(B;\RR))$.
\item The integral of $\omega_A$ vanishes on $\mbox{Ker}(i_{\ast}:H_n(A;\RR)
\longrightarrow H_n(B;\RR))$.
\end{enumerate}
\end{lemma}

\begin{lemma}
\label{lemhombordc}
The restriction map 
$$H^2(C_2(M)) \longrightarrow H^2(\partial C_2(M))$$ is an isomorphism.
\end{lemma}
\bp Write the exact sequence (with real coefficients)
$$H^2(C_2(M),\partial C_2(M)) \cong H_4(C_2(M))= 0 \longrightarrow $$ 
$$ \longrightarrow H^2(C_2(M)) \longrightarrow H^2(\partial C_2(M)) \longrightarrow $$
$$ \longrightarrow H^3(C_2(M),\partial C_2(M)) \cong H_3(C_2(M))= 0.$$
\eop

\begin{lemma}
\label{lemexisfunap}
Any closed two-form $\omega_{M}$ on $\partial C_2(M)$ extends on $C_2(M)$ as a closed two-form. If $\omega_{M}$ is antisymmetric with respect to the involution $\iota$ on $C_2(M)$, then we can demand that the extension is antisymmetric, too.
\end{lemma}
\bp 
The first assertion is a direct consequence of the two previous lemmas.
When $\omega_{M}$ is antisymmetric, let $\omega$ denote one of its closed extensions, then the average 
$$\tilde{\omega}=\frac{\omega - \iota^{\ast}(\omega)}{2}$$
is an extension of $\omega_{M}$ that is closed and antisymmetric.
\eop

In particular fundamental forms exist. Note that the cohomology class of a fundamental form is $L_M$, since its integral along the generator $S(T_xM)$ of $H_2(C_2(M))$ is one.

\subsection{Needed statements about configuration spaces}
\label{substaconf}
Compactifications of $\breve{C}_{V(\Gamma)}(M)$ are useful to study the behaviour of our integrals
$$I_{\Gamma}(\omega_M)=\int_{\breve{C}_{V(\Gamma)}(M)}\bigwedge_{e \in E(\Gamma)}p_e^{\ast}(\omega_M)$$
near $\infty$, and their dependence on the choice of $\omega_M$.
Indeed, in order to prove the convergence, it is sufficient to find a smooth
compactification (that will have corners) where the form $\bigwedge p_e^{\ast}(\omega_M)$ smoothly extends.
The variation of this integral when adding an exact form $d\eta$, will be the integral of $\eta$ on the codimension one faces of the boundary that needs to be precisely identified. 
Therefore, the proof of Theorem~\ref{thkktfra} will require a deeper knowledge of configuration spaces.
We give all the needed statements in this subsection. All of them will be proved in Section~\ref{seccomp}.

Recall that a map from $[0,\infty[^d \times \RR^{n-d}$ to $\RR^k$ is $C^{\infty}$ or smooth at $0$
if it can be extended to a $C^{\infty}$ map in a neighborhood of $0$ in $\RR^{n}$. A {\em smooth manifold with corners\/} is a manifold where every point has a neighborhood that is diffeomorphic to a neighborhood of $0$ in 
$[0,\infty[^d \times \RR^{n-d}$. 
The {\em codimension $d$ faces\/} of a smooth manifold $C$ with corners are the connected components of the set of points that are mapped to $0$ under a diffeomorphism from one of their neighborhoods to a neighborhood of $0$ in 
$[0,\infty[^d \times \RR^{n-d}$. The union of the codimension $0$ faces of such a $C$ is called the {\em interior\/} of $C$.

Let $V$ denote a finite set, let $M$ be a closed oriented three-manifold and let $X$ be a 3-dimensional vector space.

We shall study the open submanifold 
$$\breve{C}_V(M)= (M \setminus \infty)^V \setminus\mbox{all diagonals}$$
of $(M \setminus \infty)^V$. It will be seen as the space of injective maps from $V$ to $(M \setminus \infty)$. 

We shall also study the open submanifold $\breve{S}_V(X)$ of the smooth manifold $S(X^V/
\mbox{diag}(X^V))$
made of injective maps from $V$ to $X$ up to translations and dilations. 

These manifolds are our {\em configuration spaces.\/} $\breve{S}_n(X)=\breve{S}_{\{1,2,\dots, n\}}(X)$
\index{N}{Sbreven@$\breve{S}_n(X)$}
and $\breve{C}_n(M)=\breve{C}_{\{1,2,\dots, n\}}(M)$.
\index{N}{Cbreven@$\breve{C}_n(M)$}
Note that $\breve{S}_2(X)$ may be seen as the set of maps from $\{2\}$ to $X \setminus 0$ (when choosing to map $\{1\}$ to $0$). This provides a diffeomorphism
from $\breve{S}_2(X)$ to $S(X)$ that is diffeomorphic to $S^2$.

For any subset $B$ of $V$, the restriction of maps provides well-defined
projections $p_B$ from $\breve{C}_V(M)$ to $\breve{C}_B(M)$, and from $\breve{S}_V(X)$
to $\breve{S}_B(X)$. A total order on $B$ identifies $B$ to $\{1, \dots, , \sharp B\}$ and therefore identifies $\breve{C}_B(M)$ to $\breve{C}_{\sharp B}(M)$
and $\breve{S}_B(X)$ to $\breve{S}_{\sharp B}(X)$.
In particular, by composition, any ordered pair $e$ of $V$ induces canonical
maps
$$p_e:\breve{C}_{V}(M) \longrightarrow \breve{C}_{2}(M)$$
and 
$$p_e:\breve{S}_{V}(X) \longrightarrow \breve{S}_{2}(X).$$

We are going to define suitable compactifications for these spaces. Namely, we shall prove the following propositions.

\begin{proposition}
\label{propconfunc}
There exists a well-defined smooth compact manifold with corners $C_{V}(M)$ \index{N}{CV@$C_V(M)$}
whose interior is  canonically diffeomorphic to $\breve{C}_{V}(M)$
 such that
\begin{itemize}
\item $C_{\{1\}}(M)$ and $C_{\{1,2\}}(M)$ coincide with the compactifications $C_1(M)$ and $C_2(M)$
\index{N}{Ctwo@$C_2(M)$}
defined in Section~\ref{subconfap}.
\item For any ordered pair $e$ of $V$,
the projection 
$$p_e:\breve{C}_{V}(M) \longrightarrow {C}_{2}(M) $$
smoothly extends to $C_{V}(M)$.
\end{itemize}
\end{proposition}

\begin{proposition}
\label{propconfuns}
There exists a well-defined smooth compact manifold with corners $S_{V}(X)$
\index{N}{SVX@$S_{V}(X)$} 
whose interior is canonically diffeomorphic to $\breve{S}_{V}(X)$
such that, for any ordered pair $e$ of $V$,
the projection 
$$p_e:\breve{S}_{V}(X) \longrightarrow {S}_{2}(X) $$ smoothly extends to $S_{V}(X)$.
\end{proposition}

For our purposes, it will be important to know the codimension one faces 
of these compactifications. For $C_{V}(M)$, they will be the configuration spaces $F(\infty;B)$ and $F(B)$ defined below, for some subsets $B$ of $V$, where $F(\infty;B)$ will contain
limit configurations that map $B$ to $\infty$, and $F(B)$  will contain
limit configurations that map $B$ to a point of $(M \setminus \infty)$.

Let $B$ be a non-empty subset of $V$.
Let $S_i(T_{\infty}M^B)$ 
\index{N}{SiTinfty@$S_i(T_{\infty}M^B)$}
denote the set of injective maps from $B$ to $(T_{\infty}M \setminus 0)$ up to dilation. Note that $S_i(T_{\infty}M^B)$ is an open submanifold of $S((T_{\infty}M)^B)$. Define \index{N}{FinftyB@$F(\infty;B)$}
$$F(\infty;B)=\breve{C}_{(V \setminus B)}(M) \times S_i(T_{\infty}M^B)$$
where $\breve{C}_{\emptyset}(M)$ has one element.
Any ordered pair $e$ of $V$ defines a canonical map
$p_e$ from $F(\infty;B)$ to $C_2(M)$ in the following way.
\begin{itemize}
\item If $e \subseteq V \setminus B$, then $p_e$ is the composition of the natural
projections
$$F(\infty;B) \longrightarrow \breve{C}_{(V \setminus B)}(M) \longrightarrow  C_e(M)=C_2(M).$$
\item If $e \subseteq B$, then $p_e$ is the composition of the natural
maps
$$F(\infty;B) \longrightarrow S_i(T_{\infty}M^B) \longrightarrow  S_i(T_{\infty}M^e)  \hookrightarrow C_e(M)=C_2(M).$$
\item If $e \cap B =\{b^{\prime}\}$, then $p_e$ is the composition of the natural
maps
$$F(\infty;B) \longrightarrow \breve{C}_{(e \setminus \{b^{\prime}\})}(M) \times  S_i(T_{\infty}M^{\{b^{\prime}\}}) \longrightarrow $$
$$ \longrightarrow  (M \setminus \infty)^{(e \setminus \{b^{\prime}\})} \times  S(T_{\infty}M^{\{b^{\prime}\}}) \hookrightarrow C_e(M)=C_2(M).$$
\end{itemize}

Let $B$ be a subset of $V$ of cardinality $(\geq 2)$.
Let $b \in B$. Let $F(B)$ \index{N}{FB@$F(B)$} 
denote the total space of the fibration over $\left(\breve{C}_{\{b\} \cup(V \setminus B)}(M)\right)$
where the fiber over an element $c$ is $\breve{S}_B(T_{c(b)}M)$.
Again, any ordered pair $e$ of $V$ defines a canonical map
$p_e$ from $F(B)$ to $C_2(M)$.

\begin{itemize}
\item If $e \subseteq (V \setminus B) \cup \{b\}$, then $p_e$ is the composition of the natural
projections
$$F(B) \longrightarrow \breve{C}_{\{b\} \cup(V \setminus B)}(M) \longrightarrow  C_e(M)=C_2(M).$$
\item If $e \subseteq B$, then $p_e$ is the composition of the natural
projections
$$F(B) \longrightarrow \breve{S}_B(T_{c(b)}M) \longrightarrow  \breve{S}_e(T_{c(b)}M)  \longrightarrow C_e(M)=C_2(M).$$
\item If $e \cap B =\{b^{\prime}\}$, let $\tilde{e}$ be obtained from $e$ by replacing  $b^{\prime}$ by $b$, then $p_e=p_{\tilde{e}}$.
\end{itemize}

Set
$$\partial^{\infty}_1(C_{V}(M))= \{F(\infty;B); B \subseteq V; B \neq \emptyset\},$$
$$\partial^d_1(C_{V}(M))= \{F(B); B \subseteq V; \sharp B \geq 2\},$$
and 
$$\partial_1(C_{V}(M))=\partial^{\infty}_1(C_{V}(M)) \cup \partial^d_1(C_{V}(M)).$$
\index{N}{deloneC@$\partial_1(C_{V}(M))$}

The following proposition is proved in Subsection~\ref{subsketchdifcv}.

\begin{proposition}
\label{propconffaceun}
Any $F \in \partial_1(C_{V}(M))$ embeds canonically into $C_{V}(M)$, and its image is a codimension one face of $C_{V}(M)$. Therefore, any such $F$ will be identified
to its image. Then $\partial_1(C_{V}(M))$
\index{N}{deloneC@$\partial_1(C_{V}(M))$} 
is the set of codimension one faces of $C_{V}(M)$.  Furthermore, for any ordered pair $e$ of $V$, for any $F \in \partial_1(C_{V}(M))$, the restriction to $F$ of the canonical map $p_e$ defined from $C_V(M)$ to $C_2(M)$ is the map $p_e$ defined above.
\end{proposition}

Let $B$ be a strict subset of $V$ of cardinality $(\geq 2)$. Let $b \in B$.
Let $X$ be a $3$-dimensional vector space. 
Let \index{N}{fBX@$f(B)(X)$}
$$f(B)(X)=\breve{S}_B(X) \times \breve{S}_{\{b\} \cup (V \setminus B)}(X)$$
be a space of limit configurations where $B$ collapses.

Any ordered pair $e$ of $B$ provides the following canonical projection $p_e$ from $f(B)(X)$ to $S_2(X)$ as follows.

\begin{itemize}
\item If $e \subseteq (V \setminus B) \cup \{b\}$, then $p_e$ is the composition of the natural
projections
$$f(B)(X) \longrightarrow \breve{S}_{\{b\} \cup(V \setminus B)}(X) \longrightarrow  S_e(X)=S_2(X).$$
\item If $e \subseteq B$, then $p_e$ is the composition of the natural
projections
$$f(B)(X) \longrightarrow \breve{S}_B(X) \longrightarrow S_e(X)=S_2(X).$$
\item If $e \cap B =\{b^{\prime}\}$, let $\tilde{e}$ be obtained from $e$ by replacing  $b^{\prime}$ by $b$, then $p_e=p_{\tilde{e}}$.
\end{itemize}
In this article, the sign $\subset$ stands for "$\subseteq$ and $\neq$". Set \index{N}{deloneS@$\partial_1(S_{V}(X))$}
$$\partial_1(S_{V}(X))= \{f(B)(X); B \subset V; \sharp B \geq 2\}.$$
The following proposition is proved in Subsection~\ref{subsketchdifcv}.
\begin{proposition}
\label{propconffacedeux}
Any $F \in \partial_1(S_{V}(X))$ 
\index{N}{deloneS@$\partial_1(S_{V}(X))$}
embeds canonically into $S_{V}(X)$, and its image is a codimension one face of $S_{V}(X)$. Therefore, any such $F$ will be identified
to its image. Then $\partial_1(S_{V}(X))$ is the set of codimension one faces of $S_{V}(X)$.
Furthermore, for any ordered pair $e$ of $V$, for any $F \in \partial_1(S_{V}(X))$, the restriction to $F$ of the canonical map $p_e$ defined from $S_V(X)$ to $S_2(X)$ is the map $p_e$ defined above.
\end{proposition}

\subsection{Sketch of the proof of Theorem~\ref{thkktfra}}
\label{subsketchpkkt}

We shall first see that the wanted invariant $Z$ is the exponential of
a simpler series in $\CA(\emptyset)$, that we are going to present
in another way by means of {\em labelled diagrams\/} that will make the proofs clearer.

A degree $n$ {\em labelled\/} Jacobi diagram 
\index{T}{Jacobi diagram!labelled}
is a Jacobi diagram whose vertices are numbered from $1$ to $2n$, and whose edges are numbered from $1$ to $3n$.

Let $\overline{\Gamma}$ be a labelled Jacobi diagram with underlying Jacobi diagram $\Gamma$. The automorphisms of $\Gamma$ 
\index{T}{Jacobi diagram!automorphism of}
act on the labelling of $\overline{\Gamma}$. In particular, there are exactly $\sharp \mbox{Aut}(\Gamma)$ 
\index{N}{AutGamma@$\sharp \mbox{Aut}(\Gamma)$}
labellings of $\overline{\Gamma}$ that give rise to a labelled Jacobi diagram isomorphic to $\overline{\Gamma}$ as a labelled Jacobi diagram, and the number of labelled Jacobi diagrams with underlying Jacobi diagram $\Gamma$ is $\frac{(2n)!(3n)!}{\sharp \mbox{Aut}(\Gamma)}$.

A Jacobi diagram is {\em edge-oriented\/} 
\index{T}{Jacobi diagram!edge-oriented}
when its edges are oriented. Any labelled Jacobi diagram has $2^{3n}$ such edge-orientations.

A labelled edge-oriented Jacobi diagram inherits a canonical vertex-orientation (up to an even 
number of changes), namely the vertex-orientation that together with the orientation of $V(\Gamma)$ induced by the vertex labels provides a vertex-orientation of $H(\Gamma)$ equivalent to its edge-orientation.
Therefore, an edge-oriented labelled graph $\overline{\Gamma}$ has a well-determined class $[\overline{\Gamma}]$ in $\CA(\emptyset)$.
Furthermore,  an edge-oriented labelled graph $\overline{\Gamma}$ defines a map \index{N}{PGamma@$P(\Gamma)$}
$$P(\overline{\Gamma}):\breve{C}_{2n}(M)  \longrightarrow C_{2}(M)^{3n}$$
whose projection $p_i \circ P(\overline{\Gamma})=P_i(\overline{\Gamma})$
\index{N}{PiGamma@$P_i(\Gamma)$} 
onto the $i^{th}$ factor of $C_{2}(M)^{3n}$ is $p_{e(i)}$
where $e(i)$ denotes the edge labelled by $i$.

Let $\tau_M$ be a trivialisation of $T(M \setminus \infty)$ standard near $\infty$.
For any $i \in \{1, \dots, 3n\}$, let $\omega^{(i)}_M$ be a two-form that is fundamental with respect to $\tau_M$ and to a form $\omega^{(i)}_{S^2}$ such that $\int_{S^2}\omega^{(i)}_{S^2}=1$. Define the $6n$-form on $C_{2}(M)^{3n}$ \index{N}{Omega@$\Omega$}
$$\Omega=\bigwedge_{i=1}^{3n}p_i^{\ast}(\omega^{(i)}_M).$$
Proposition~\ref{propconfunc} allows us to define \index{N}{IGammaMOmega@$I_{\Gamma}(M;\Omega)$} $$I_{\overline{\Gamma}}(M;\Omega)=\int_{{C}_{2n}(M)}P(\overline{\Gamma})^{\ast}(\Omega)=\int_{{C}_{2n}(M)}\bigwedge_{i=1}^{3n}P_i(\overline{\Gamma})^{\ast}(\omega^{(i)}_M).$$
Let ${\cal E}_n$ 
\index{N}{Ecaln@${\cal E}_n$} 
denote the set of all connected edge-oriented labelled Jacobi diagrams with $2n$ vertices. 
We are going to prove the following propositions.

\begin{proposition}
\label{propun}
Under the above assumptions, \index{N}{zntauM@$z_n(\tau_M)$}
$$z_n(\tau_M) = \sum_{\Gamma \in {\cal E}_n} I_{{\Gamma}}(M;\Omega)[{\Gamma}]$$
only depends on $M$ and on $\tau_M$.
\end{proposition}

\begin{proposition}
\label{propdeux}
Let $\omega_T$ be a closed two-form \index{N}{omegaT@$\omega_T$}
that represents the Thom class of $S_2(E_1)$. \index{N}{StwoEone@$S_2(E_1)$}
Then \index{N}{deltan@$\delta_n$}
$$\delta_n= \sum_{\Gamma \in {\cal E}_n} \int_{{S}_{2n}(E_1)}\bigwedge_{i=1}^{3n}P_i({\Gamma})^{\ast}(\omega_T)[{\Gamma}]$$
does not depend on the choice of $\omega_T$. 
\end{proposition}

\begin{proposition}
\label{proptrois}
Under the above assumptions, \index{N}{zntauM@$z_n(\tau_M)$}
$z_n(\tau_M)$
only depends on $M$ and on the homotopy class of $\tau_M$ among the trivialisations that are standard near $\infty$.

For any closed 3-manifold $M$, for any 
trivialisation $\tau_M$ of $T(M \setminus \infty)$ that is standard near $\infty$, and for any $$g: (M \setminus \infty,  M \setminus B_M(1)) \longrightarrow (SO(3),1),$$
define $$\begin{array}{llll} 
\psi(g): &(M \setminus \infty) \times \RR^3 &\longrightarrow  &(M \setminus \infty) \times \RR^3\\
&(x,y) & \mapsto &(x,g(x)(y))\end{array}$$ then \index{N}{deltan@$\delta_n$} 
$$z_n(\psi(g) \circ \tau_M)-z_n(\tau_M)=\frac{1}{2}\mbox{deg}(g)\delta_n.$$
\end{proposition}

Note that Propositions~\ref{propconfunc} and \ref{propconfuns}
ensure that all the mentioned integrals
are well-defined and that all the previous ones are convergent.

Let us now show that Propositions~\ref{propun}, \ref{propdeux}, \ref{proptrois}, \ref{proppont}
prove Theorem~\ref{thkktfra}.

First note that for an antisymmetric $\omega_M$ that is fundamental with respect to $\tau_M$ and to a two-form $\omega_{S^2}$ such that $\int_{S^2}\omega_{S^2}=1$, $\int_{{C}_{2n}(M)}\bigwedge_{i=1}^{3n}P_i(\overline{\Gamma})^{\ast}(\omega_M)[\overline{\Gamma}]$
is independent of the labelling and is equal to $I_{\Gamma}(\omega_M)[\Gamma]$.
Therefore, \index{N}{zntauM@$z_n(\tau_M)$}
$$z_n(\tau_M)=2^{3n}(3n)!(2n)!\sum_{\Gamma \;\mbox{\small connected Jacobi diagrams with 2n vertices}}\frac{I_{\Gamma}(\omega_M)}{\sharp\mbox{Aut}(\Gamma)}
[\Gamma],$$
and $$\tilde{Z}(\tau_M)=\exp\left(\left(\frac{1}{2^{3n}(3n)!(2n)!}z_n(\tau_M)\right)_n\right)$$ only depends on $\tau_M$, according to Proposition~\ref{propun}.
\begin{lemma}
$Z(M;\tau_M)=\tilde{Z}(\tau_M)$. \index{N}{ZMtauM@$Z(M;\tau_M)$}
\end{lemma}
These are two series of combinations of diagrams and it suffices to compare the coefficients
of $[\Gamma]$, for a diagram $\Gamma$ which is a disjoint union of $k_1$ copies of $\Gamma_1$,
 $k_2$ copies of $\Gamma_2$, \dots,  $k_r$ copies of $\Gamma_r$, where $\Gamma_1$, $\Gamma_2$ and $\Gamma_r$ are non-isomorphic connected Jacobi diagrams. The coefficient of $[\Gamma]=\prod_{i=1}^r[\Gamma_i]^{k_i}$ in $Z(M;\tau_M)$ is 
$\frac{I_{\Gamma}(\omega_M)}{\sharp\mbox{Aut}(\Gamma)}$ 
where $I_{\Gamma}(\omega_M) = \prod_{i=1}^r I_{\Gamma_i}(\omega_M)^{k_i}$, and
$\sharp\mbox{Aut}(\Gamma)=\prod_{i=1}^r[(\sharp\mbox{Aut}(\Gamma_i))^{k_i}(k_i)!]$. Therefore, the coefficient in $Z(M;\tau_M)$ is
$$\prod_{i=1}^r\frac{I_{\Gamma_i}(\omega_M)^{k_i}}{\sharp\mbox{Aut}(\Gamma_i)^{k_i}(k_i)!}.$$
Let $k=\sum_{i=1}^rk_i$. The coefficient of $[\Gamma]$ in $\tilde{Z}(\tau_M)$ 
is 
its coefficient in the product
$$\frac{1}{k!}\left(\left(\frac{1}{2^{3n}(3n)!(2n)!}z_n(\tau_M)\right)_n\right)^k$$ where $\prod_{i=1}^r[\Gamma_i]^{k_i}$ occurs $\frac{k!}{ \prod_{i=1}^r (k_i)!}$ times
with the coefficient $\frac{1}{k!}\prod_{i=1}^r\frac{I_{\Gamma_i}(\omega_M)^{k_i}}{\sharp\mbox{Aut}(\Gamma_i)^{k_i}}$.
Thus, the two coefficients coincide.
\eop

Of course, Proposition~\ref{propdeux} implies that \index{N}{ksin@$\xi_n$} \index{N}{deltan@$\delta_n$}
$$\xi_n= \frac{1}{2^{3n}(3n)!(2n)!}\delta_n$$ is independent on the used $\omega_T$.
Now, Propositions~\ref{proptrois}
and \ref{proppont} 
clearly imply that
$(\frac{1}{2^{3n}(3n)!(2n)!}z_n(\tau_M) + \frac{p_1(\tau_M)}{4}\xi_n)$ \index{N}{ksin@$\xi_n$} is 
independent of $\tau_M$, and this in turn implies Theorem~\ref{thkktfra}.

Propositions~\ref{propun}, \ref{propdeux} and \ref{proptrois}
and \ref{proppont} 
will be proved in Subsections~\ref{subproofpropun}, \ref{subsproofdeux}, \ref{submoretriv}
and \ref{subproofpont}, respectively.

\subsection{Proof of Proposition~\ref{propun}, the dependence on the forms.}
\label{subproofpropun}

In this subsection, we prove Proposition~\ref{propun}.

Of course, the only choice in the expression of $z_n(\tau_M)$ is the choice of the $\omega^{(i)}_M$, and it is enough to prove that changing an $\omega^{(i)}_M$ into an $\hat{\omega}^{(i)}_M$ that is fundamental with respect to
 $\tau_M$ and to a form $\hat{\omega}^{(i)}_{S^2}$  such that $\int_{S^2}\hat{\omega}^{(i)}_{S^2}=1$ does not change $z_n=z_n(\tau_M)=z_n(\Omega)$. 

For later use in \cite{sumgen}, we shall rather study how $z_n$ varies when $\omega^{(i)}_M$
varies within a class of forms that is more general than the fundamental forms. 

\begin{definition}
\label{defformad}
\index{T}{form!admissible}
\noindent A two-form $\omega_M$ on $C_2(M)$ or on $\partial C_2(M)$ is {\em admissible\/}  if:\\
$\bullet$ its restriction to $\partial C_2(M) \setminus ST(B_M)$ is $p_M(\tau_M)^{\ast}(\omega_{S^2})$ for some trivialisation $\tau_M$ of $T(M \setminus \infty)$ standard near $\infty$ and for some two-form $\omega_{S^2}$ on $S^2$ with total volume one, and,\\ 
$\bullet$ it is closed.\\
Such a two-form is  {\em antisymmetric\/} if  $\iota^{\ast}(\omega_M)=-\omega_M$. 
\end{definition}

According to Lemma~\ref{lemexisfunap}, an admissible two-form on $\partial C_2(M)$ extends as an admissible two-form on $C_2(M)$, and an admissible antisymmetric two-form on $\partial C_2(M)$ extends as an admissible antisymmetric two-form on $C_2(M)$.
We are going to prove the following proposition.

\begin{proposition}
\label{propzad}
Let \/ $\omega_M$ be an antisymmetric admissible two-form on $C_2(M)$, then
with the notation before Proposition~\ref{propthkktun} \index{N}{ZnomegaM@$Z_n(\omega_M)$}
$$Z_n(\omega_M) =\sum_{\Gamma \;\mbox{\small Jacobi diagram with $2n$ vertices}}\frac{I_{\Gamma}(\omega_M)}{\sharp \mbox{Aut}(\Gamma)}[\Gamma]$$
only depends on $M$ and on the restriction of $\omega_M$ to $ST(B_M)$.
\end{proposition}

In what follows, all the two forms $\omega^{(j)}_M$, for $j \in \{1,2,\dots,2n\}$, and $\hat{\omega}^{(i)}_M$ are admissible with respect to two-forms on $S^2$ denoted by $\omega^{(j)}_{S^2}$, for $j \in \{1,2,\dots,2n\}$ and $\hat{\omega}^{(i)}_{S^2}$, respectively. Note that the restriction of $\omega^{(j)}_M$ on $ST(B_M)$ determines $\omega^{(j)}_{S^2}$, and hence determines $\omega^{(j)}_M$ on $\partial C_2(M)$.
We fix a trivialisation $\tau_M$ of $T(M \setminus \infty)$ standard near $\infty$.

\begin{lemma}
\label{lemnolosseta}
There 
exists a one-form $\eta_{S^2}$ on $S^2$ such that $d \eta_{S^2} = \hat{\omega}^{(i)}_{S^2}-{\omega}^{(i)}_{S^2}$, and
a one-form $\eta$ on $C_2(M)$ such that
\begin{enumerate}
\item $d\eta=\hat{\omega}^{(i)}_M-\omega^{(i)}_M$,
\item the restriction of $\eta$ on $\partial C_2(M) \setminus ST(B_M)$ is $p_M(\tau_M)^{\ast}(\eta_{S^2})$,
\item if $\hat{\omega}^{(i)}_M$ and $\omega^{(i)}_M$ are fundamental with respect to $\tau_M$, then the restriction of $\eta$ on the whole $ \partial C_2(M)$ is $p_M(\tau_M)^{\ast}(\eta_{S^2})$,
\item if $\hat{\omega}^{(i)}_M$ and $\omega^{(i)}_M$ coincide on $ST(B_M)$, then the restriction of $\eta$ on $ \partial C_2(M)$ is zero.
\end{enumerate}
\end{lemma}
\bp
Since $\hat{\omega}^{(i)}_M$ and $\omega^{(i)}_M$ are cohomologous there exists
$\eta$ such that $d\eta=\hat{\omega}^{(i)}_M-\omega^{(i)}_M$ on $C_2(M)$.
Similarly, there exists $\eta_{S^2}$ such that $d \eta_{S^2} = \hat{\omega}^{(i)}_{S^2}-{\omega}^{(i)}_{S^2}$.  If  $\hat{\omega}^{(i)}_M$ and $\omega^{(i)}_M$ coincide on $ST(B_M)$, then $\hat{\omega}^{(i)}_{S^2}={\omega}^{(i)}_{S^2}$, and we choose $ \eta_{S^2}=0$.
Now, $ d(\eta -p_M(\tau_M)^{\ast}(\eta_{S^2}))=0$ on $\partial C_2(M) \setminus ST(B_M)$, and on $ \partial C_2(M)$ if $\hat{\omega}^{(i)}_M$ and $\omega^{(i)}_M$ are fundamental with respect to $\tau_M$, or if $\hat{\omega}^{(i)}_M$ and $\omega^{(i)}_M$ coincide on $ST(B_M)$.

Thanks to the exact sequence
$$ 0=H^1(C_2(M)) \longrightarrow H^1(\partial C_2(M)) \longrightarrow H^2(C_2(M), \partial C_2(M)) \cong H_4(C_2(M))=0, $$
$H^1(\partial C_2(M))=0$. \index{N}{Ctwo@$C_2(M)$} 
It is easy to see that $H^1(\partial C_2(M)\setminus ST(B_M))=0$, too.
Therefore, there exists a function $f$ from $\partial C_2(M)$ to $\RR$ such that $$df =\eta -p_M(\tau_M)^{\ast}(\eta_{S^2})$$
on $\partial C_2(M) \setminus ST(B_M)$, and on $ \partial C_2(M)$ if $\hat{\omega}^{(i)}_M$ and $\omega^{(i)}_M$ are fundamental with respect to $\tau_M$ or if $\hat{\omega}^{(i)}_M$ and $\omega^{(i)}_M$ coincide on $ST(B_M)$.
Extend $f$ to a $C^{\infty}$ map on $C_2(M)$ and change $\eta$ into $(\eta -df)$.
\eop

Set $$z_n=\sum_{\Gamma \in {\cal E}_n} \int_{{C}_{2n}(M)}\bigwedge_{i=1}^{3n}P_i({\Gamma})^{\ast}(\omega^{(i)}_M)[{\Gamma}].$$
Set $\hat{\omega}^{(j)}_M={\omega}^{(j)}_M$ for $j \neq i$, and
let $\hat{z}_n=\sum_{\Gamma \in {\cal E}_n} \int_{{C}_{2n}(M)}\bigwedge_{i=1}^{3n}P_i({\Gamma})^{\ast}(\hat{\omega}^{(i)}_M)[{\Gamma}]$.

Set $$\tilde{\omega}^{(j)}_M= \left\{ 
\begin{array}{ll} {\omega}^{(j)}_M\;&\mbox{if} \; j \neq i\\
\eta \;& \mbox{if} \; j = i, \end{array}
\right.$$
and define the $(6n-1)$-form $\tilde{\Omega}=
\bigwedge_{j=1}^{3n}p_j^{\ast}(\tilde{\omega}^{(j)}_M)$ on $C_{2}(M)^{3n}$.
Then $d\tilde{\Omega}= \bigwedge_{j=1}^{3n}p_j^{\ast}(\hat{\omega}^{(j)}_M)-\bigwedge_{j=1}^{3n}p_j^{\ast}({\omega}^{(j)}_M)$.

For an element $F$ of the set $\partial_1(C_{2n}(M))$ of codimension $1$ faces of $C_{2n}(M)$ described before Proposition~\ref{propconffaceun}, set \index{N}{IGammaF@$I_{\Gamma,F}$}
$$I_{{\Gamma},F}=\int_{F}P({\Gamma})^{\ast}(\tilde{\Omega})=
\int_{F}\bigwedge_{j=1}^{3n}P_j({\Gamma})^{\ast}(\tilde{\omega}^{(j)}_M),$$
where $F$ is oriented as a part of the boundary of the oriented manifold $C_{2n}(M)$,
$$I_{{\Gamma},\partial}=\sum_{F \in \partial_1(C_{2n}(M))}I_{{\Gamma},F}.$$
Then according to the Stokes theorem, $$\hat{z}_n-z_n= \sum_{\Gamma \in {\cal E}_n} I_{{\Gamma},\partial}[\Gamma].$$

We are going to prove that several terms cancel in this sum. More precisely, we shall prove Proposition~\ref{propunb} that obviously implies Proposition~\ref{propun}, and Proposition~\ref{propzad} since
$$Z(\omega_M)=\exp\left(\left(\frac{1}{2^{3n}(3n)!(2n)!}z_n(\omega_M)\right)_n\right)$$
with $$z_n(\omega_M)=\sum_{\Gamma \in {\cal E}_n} \int_{{C}_{2n}(M)}\bigwedge_{i=1}^{3n}P_i({\Gamma})^{\ast}(\omega_M)[{\Gamma}].$$

\begin{proposition}
\label{propunb}
With the notation above,
$$\hat{z}_n-z_n= \sum_{\Gamma \in {\cal E}_n} I_{{\Gamma},F(V)}[\Gamma]$$
where $I_{{\Gamma},F(V)}=0$ for any $\Gamma \in {\cal E}_n$ if $\hat{\omega}^{(i)}_M$ and $\omega^{(i)}_M$ are fundamental with respect to $\tau_M$, or if $\hat{\omega}^{(i)}_M$ and $\omega^{(i)}_M$ coincide on $ST(B_M)$.
\end{proposition}


When $B$ is a subset of $V$, and when a graph $\Gamma$ is given, $E_B$ denotes the set of edges of $\Gamma$ that contain two elements of $B$, 
and $\Gamma_B$ \index{N}{GammaB@$\Gamma_B$}
is the subgraph of $\Gamma$ made of the vertices of $B$ and the edges of $E_B$.

\begin{lemma}
\label{lemfaceinf}
For any non-empty subset $B$ of $V$, for any $\Gamma \in {\cal E}_n$,
$I_{{\Gamma},F(\infty;B)}=0$.
\end{lemma}
\bp
Set $A = V \setminus B$. 
Let $E_C$ be the set of the edges of $\Gamma$ that contain an element of $A$ and an element of $B$.
Let $p_2$ denote the projection of $F(\infty;B)$ onto $S_i(T_{\infty}M^B)$.
For $e \in E^B \cup E^C$, $P_e:(S^2)^{E^B \cup E^C} \longrightarrow S^2$ is the projection onto the factor indexed by $e$.
We show that there exists a smooth map
$$g : S_i(T_{\infty}M^B) \longrightarrow (S^2)^{E^B \cup E^C}$$
such that 
$$\bigwedge_{e \in E^B \cup E^C}p_e^{\ast}(\tilde{\omega}^{(i(e))}_M)
=(g \circ p_2)^{\ast}\left( \bigwedge_{e \in E^B \cup E^C}P_e^{\ast}(\tilde{\omega}^{(i(e))}_{S^2}) \right)$$
where $i(e) \in \{1,2, \dots, 3n\}$ is the label of the edge $e$, and
$$\tilde{\omega}^{(i(e))}_{S^2}= \left\{ 
\begin{array}{ll} {\omega}^{(i(e))}_{S^2}\;&\mbox{if} \; i(e) \neq i\\
\eta_{S^2}\;& \mbox{if} \; i(e) = i. \end{array}
\right.$$
Indeed, if $e \in E^B \cup E^C$, $p_e(F(\infty;B)) \subset \partial C_2(M)$, $$p_e^{\ast}(\tilde{\omega}^{(i(e))}_M)=(p_M(\tau_M) \circ p_e)^{\ast}(\tilde{\omega}^{(i(e))}_{S^2}),$$
and $p_M(\tau_M) \circ p_e$ factors through $S_i(T_{\infty}M^B)$ (and therefore reads 
$((P_e \circ g) \circ p_2)$). Indeed,
if $e \in E^C$, $p_M(\tau_M) \circ p_e$ only depends on the projection on $S(T_{\infty}M)$ of the vertex at $\infty$ (of $B$),
 while, if $e \in E^B$, $p_M(\tau_M) \circ p_e$ factors through $S_i(T_{\infty}M^e)$.

Therefore if the degree of the form $\left( \bigwedge_{e \in E^B \cup E^C}p_e^{\ast}(\tilde{\omega}^{(i(e))}_{S^2}) \right)$ is bigger than the dimension $(3 \sharp B-1)$ of $S_i(T_{\infty}M^B)$, this form vanishes on $F(\infty;B)$.
The degree of the form 
is $( 2\sharp E^B + 2\sharp E^C)$ or $( 2\sharp E^B + 2\sharp E^C-1)$, while
$$(3 \sharp B-1)= 2\sharp E^B + \sharp E^C-1.$$
Therefore, the integral vanishes unless $E^C$ is empty.
In this case, since $\Gamma$ is connected, $B=V$, 
$F(\infty;V)= S_i(T_{\infty}M^V)$, all the $p_M(\tau_M) \circ p_e$
locally factor through the conjugates under the inversion $(x \mapsto x/\norm{x}^2)$ of the translations that make sense, and the form vanishes, too.
\eop

As soon as there exists a smooth map from $F(B)$ to a manifold of strictly smaller
dimension that factorizes $P(\Gamma)$, then $I_{{\Gamma},F(B)}=0$.
We shall use this principle to get rid of some faces.

\begin{lemma}
\label{lemdiscon}
Let $\Gamma \in \CE_n$. For any subset $B$ of $V$ such that $\Gamma_B$ is not connected,
$I_{{\Gamma},F(B)}=0$.
\end{lemma}
\bp Indeed, in the fiber $\breve{S}_B(T_{c(b)}M)$ we may translate
one connected component of $\Gamma_B$ whose set of  vertices is  $C$ independently. This amounts to factorize the $p_e$ through  $\breve{C}_{\{b\} \cup(V \setminus B)}(M)$ if $\sharp B=2$, or through the fibered space over $\breve{C}_{\{b\} \cup(V \setminus B)}(M)$ whose fiber is an open subspace of
$$S\left(\frac{T_{c(b)}M^B}
{\mbox{diag}(T_{c(b)}M^B) \oplus (0^{B \setminus C} \times \mbox{diag}(T_{c(b)}M^C))}
\right).$$ In both cases, all the $p_e$ factor through
a space with smaller dimension.
\eop

\begin{lemma}
\label{lemedge}
Let $\Gamma \in \CE_n$. Let $B$ be a subset of $V$ such that $\sharp B \geq 3$. If some element of $B$ belongs to exactly one edge of $\Gamma_B$, then
$I_{{\Gamma},F(B)}=0$.
\end{lemma}
\bp Let $b$ be the mentioned element, and let $e$ be its edge in  $\Gamma_B$, let $d \in B$ be the other element of $e$. The group $]0,\infty[$ acts on the map $t$ from $B$ to $T_{c(b)}M$ by moving $t(b)$ on the half-line from $t(d)$ through $t(b)$. ($(t(b)-t(d))$ is multiplied by a scalar). When $\sharp B \geq 3$, this action is non trivial on $\breve{S}_B(T_{c(b)}M)$, $P(\Gamma)$
factors through the quotient of $F(B)$ by this action that has one less dimension.
\eop

\begin{lemma}
\label{lemsym}
Let $\Gamma \in \CE_n$. Let $B$ be a subset of $V$ such that at least one element of $B$ belongs to exactly two edges of $\Gamma_B$. Let $\CE(\Gamma)$  denote the set of labelled edge-oriented graph that are isomorphic to $\Gamma$ by an isomorphism that preserves the labels of the vertices, but that may change the labels and the orientations of the edges. Then 
$$\sum_{\tilde{\Gamma}; \tilde{\Gamma} \in \CE(\Gamma)}I_{\tilde{\Gamma},F(B)}[\tilde{\Gamma}]=0.$$
\end{lemma}
\bp Let $v_m$ be the vertex of $B$ with smallest label $m \in \{1,2,\dots, 2n\}$ that belongs to exactly two edges of $\Gamma_B$.
We first describe an orientation-reversing diffeomorphism of the complement of a codimension $3$ submanifold of $F(B)$. Let $v_j$ and $v_k$ denote the (possibly 
equal) two other vertices of the two edges of $\Gamma_B$ that contain $v_m$.
Consider the linear transformation $S$ of the space $S(T_{c(b)}M^B/\mbox{diag}(T_{c(b)}M^B))$ of non-constant maps $f$ from $B$ to $T_{c(b)}M$
up to translations and dilations, that maps $f$ to $S(f)$ where\\
$S(f(v_{\ell}))=f(v_{\ell})$ if $v_\ell \neq v_m$, and,\\
$S(f(v_m))=f(v_j)+f(v_k)-f(v_m)$.\\
This is an orientation-reversing involution of $S(T_{c(b)}M^B/\mbox{diag}(T_{c(b)}M^B))$.
The set of elements of $\breve{S}_B(T_{c(b)}M)$
whose image under $S$ is not in $\breve{S}_B(T_{c(b)}M)$ is a codimension $3$ submanifold of $\breve{S}_B(T_{c(b)}M)$.
The fibered product of $S$ by the identity of the base $\breve{C}_{\{b\} \cup(V \setminus B)}(M)$ is an orientation-reversing smooth involution outside a codimension $3$ submanifold $F_S$ of $F(B)$. It is still denoted by $S$.

Now, let $\sigma(B;\Gamma)(\tilde{\Gamma})$ be obtained from $(\tilde{\Gamma} \in \CE(\Gamma))$ by reversing the orientations of the edges of $\Gamma_B$ that contain $v_m$ and by exchanging 
their labels. Then, as the following picture shows,
$$P(\tilde{\Gamma}) \circ S=P(\sigma(B;\Gamma)(\tilde{\Gamma})).
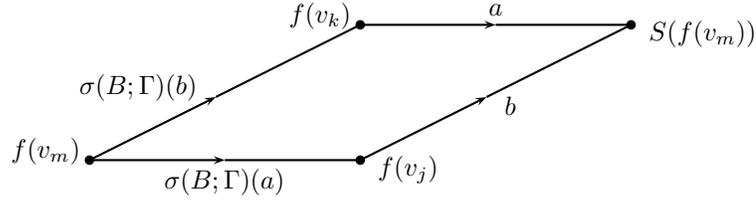
$$

$$\begin{pspicture}[.2](0,0)(8,2.5)
\psset{xunit=1.2cm,yunit=1.2cm}
\rput[b](5.5,2.1){$a$}
\rput[t](2.5,.4){$\sigma(B;\Gamma)(a)$}
\rput[l](5.6,1.1){$b$}
\rput[r](2.2,1.3){$\sigma(B;\Gamma)(b)$}
\rput[r](.9,.6){$f(v_m)$}
\rput[r](3.9,2.1){$f(v_k)$}
\rput[l](4.2,.4){$f(v_j)$}
\rput[l](7.2,1.9){$S(f(v_m))$}
\psline{->}(4,.5)(5.4,1.2)
\psline{-}(5.4,1.2)(7,2)
\psline{->}(1,.5)(2.4,1.2)
\psline{-}(2.4,1.2)(4,2)
\psline{*->}(4,2)(5.5,2)
\psline{-*}(5.5,2)(7,2)
\psline{*->}(1,.5)(2.5,.5)
\psline{-*}(2.5,.5)(4,.5)
\end{pspicture}$$

Therefore, 
$$\begin{array}{ll}
I_{\tilde{\Gamma},F(B)} &=\int_{F(B) \setminus F_S}P(\tilde{\Gamma})^{\ast}(\tilde{\Omega})\\
&=-\int_{F(B) \setminus F_S}S^{\ast}\left(P(\tilde{\Gamma})^{\ast}(\tilde{\Omega})\right)\\
&=-\int_{F(B) \setminus F_S}(P(\tilde{\Gamma}) \circ S)^{\ast}(\tilde{\Omega})\\ &=-\int_{F(B) \setminus F_S}P(\sigma(B;\Gamma)(\tilde{\Gamma}))^{\ast}(\tilde{\Omega})\\
&=-I_{\sigma(B;\Gamma)(\tilde{\Gamma}),F(B)}\end{array}$$
while $[\tilde{\Gamma}]=[\sigma(B;\Gamma)(\tilde{\Gamma})]$.
Now, $\sigma(B;\Gamma)$ defines an involution of $\CE(\Gamma)$, and it is easy to conclude:
$$\sum_{\tilde{\Gamma}; \tilde{\Gamma} \in \CE(\Gamma)}I_{\tilde{\Gamma},F(B)}[\tilde{\Gamma}]=\sum_{\tilde{\Gamma}; \tilde{\Gamma} \in \CE(\Gamma)}I_{\sigma(B;\Gamma)(\tilde{\Gamma}),F(B)}[\sigma(B;\Gamma)(\tilde{\Gamma})]$$
$$=-\sum_{\tilde{\Gamma}; \tilde{\Gamma} \in \CE(\Gamma)}I_{\tilde{\Gamma},F(B)}[\tilde{\Gamma}]=0.$$
\eop

The symmetry used in the above proof was observed by Kontsevich in \cite{ko}.

The three previous lemmas allow us to get rid of the pairs $(B;\Gamma)$ with  $\sharp B \geq 3$ such that at least one element of $B$ does not belong to three edges. Therefore, since the $\Gamma$ are connected, we are left
with the pairs $(B;\Gamma)$ with $B=V$, that are treated by Lemma~\ref{lemfacetot} below, and with the pairs $(B;\Gamma)$
where $B \neq V$,  $\sharp B=2$, and at least one element belongs to
exactly one edge of $\Gamma_B$. The following lemma allows us to get rid of this latter case where $\Gamma_B$ must be an edge.

\begin{lemma}
\label{lemihx}
Let $\Gamma \in \CE_n$. Let $B$ be a subset of $V$
such that $\Gamma_B$ is made of an edge $e(\ell)$  with label $\ell$ oriented from a vertex $v_j$ to a vertex $v_k$.
Let $\Gamma/\Gamma_B$ be the labelled edge-oriented graph obtained from $\Gamma$
by collapsing $\Gamma_B$ down to one point. 
(The labels of the edges of $\;\Gamma/\Gamma_B$ belong to $\{1, 2, \dots, 3n\} \setminus \{\ell\}$, the labels of the vertices of $\;\Gamma/\Gamma_B$ belong to $\{1, 2, \dots, 2n\} \setminus \{k\}$, $\Gamma/\Gamma_B$ has one four-valent vertex $(v_j=v_k)$ and its other vertices are trivalent.) 
Let $\CE(\Gamma;B)$ be the subset of $\CE_n$ that contains the graphs $\tilde{\Gamma}$ whose edge with label $\ell$ goes from $v_j$ to $v_k$ and such that $\Gamma/\Gamma_B$ is equal to $\tilde{\Gamma}/\tilde{\Gamma}_B$. Then 
$$\sum_{\tilde{\Gamma}; \tilde{\Gamma} \in \CE(\Gamma;B)}I_{\tilde{\Gamma},F(B)}[\tilde{\Gamma}]=0.$$
\end{lemma}
\bp
$F(B)$ is fibered over $\breve{C}_{V \setminus \{v_k\}}(M)$ with fiber $ST_{c(v_j)}M$ that contains the direction of the vector from $c(v_j)$ to $c(v_k)$. The oriented face $F(B)$ and the map
$$P(\tilde{\Gamma}): (F(B) \subset C_{2n}(M)) \longrightarrow C_2(M)^{3n}$$
are the same for all the elements $\tilde{\Gamma}$ of $\CE(\Gamma;B)$.
Therefore $I_{\tilde{\Gamma},F(B)}$ is the same for all the elements $\tilde{\Gamma}$ of $\CE(\Gamma;B)$, the sum of the statement is 
$$\sum_{\tilde{\Gamma}; \tilde{\Gamma} \in \CE(\Gamma;B)}I_{\tilde{\Gamma},F(B)}[\tilde{\Gamma}]
=I_{{\Gamma},F(B)}\sum_{\tilde{\Gamma}; \tilde{\Gamma} \in \CE(\Gamma;B)}[\tilde{\Gamma}],$$
and we are left with the study of the set $\CE(\Gamma;B)$.
Let $\tilde{\Gamma} \in \CE(\Gamma;B)$.
Let $a,b,c,d$ be the four half-edges of 
$\Gamma/\Gamma_B$ that contain $(v_j=v_k)$. Let $e_1$ be the first half-edge of $e(\ell)$ that contains
$v_j$, and let $e_2$ be the other half-edge of $e(\ell)$. Then in $\tilde{\Gamma}$, $v_j$ belongs to $e_1$ and to two half-edges of $\{a,b,c,d\}$, and the corresponding unordered pair 
determines $\tilde{\Gamma}$ as an edge-oriented labelled graph.
Thus, there are 6 graphs in $\CE(\Gamma;B)$ labelled by the pairs of elements
of $\{a,b,c,d\}$.
Equip $\Gamma=\Gamma_{ab}$ with a vertex-orientation that reads
$(a,b,e_1)$ at $v_j$ and $(c,d,e_2)$ at $v_k$ and that is consistent
with its given edge-orientation (i.e. such that the edge-orientation of $H(\Gamma)$
is equivalent to its vertex-orientation).
A representative of the orientation of $H(\Gamma)$ reads
$( \dots, a,b,e_1, \dots , c,d,e_2, \dots )$
and is equivalent to the edge-orientation of $H(\Gamma)$ that is the same for all the elements of $\CE(\Gamma;B)$.

Thus, cyclically permuting the letters $b,c,d$ gives rise to two other graphs
in $\CE(\Gamma;B)$ equipped with a suitable vertex-orientation, that respectively reads\\ 
$(a,c,e_1)$ at $v_j$ and $(d,b,e_2)$ at $v_k$, or\\
$(a,d,e_1)$ at $v_j$ and $(b,c,e_2)$ at $v_k$,\\
The three other elements of $\CE(\Gamma;B)$ with their suitable vertex-orientation  are obtained from the three previous ones by exchanging the ordered pair before $e_1$ with the ordered pair before $e_2$. This amounts to exchanging the vertices $v_j$ and $v_k$ in the picture, and does not change the unlabelled vertex-oriented graph. The first three graphs can be represented by three graphs identical outside the pictured disk:

$$
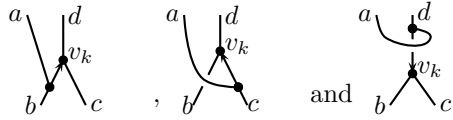
\begin{pspicture}[.2](0,-.2)(1.6,1)
\psset{xunit=1.2cm,yunit=1.2cm}
\rput[r](.05,1){$a$}
\rput[r](.2,0){$b$}
\rput[l](.8,0){$c$}
\rput[l](.55,1){$d$}
\rput[l](.55,.55){$v_k$}
\psline{-*}(.1,1)(.35,.2)
\psline{*-}(.5,.5)(.5,1)
\psline{-}(.75,0)(.5,.5)
\psline{->}(.25,0)(.5,.5)
\end{pspicture}
\;\; \mbox{,} \;\;
\begin{pspicture}[.2](0,-.2)(1.6,1)
\psset{xunit=1.2cm,yunit=1.2cm}
\rput[r](.05,1){$a$}
\rput[r](.2,0){$b$}
\rput[l](.8,0){$c$}
\rput[l](.55,1){$d$}
\rput[l](.55,.65){$v_k$}
\psline{*-}(.5,.6)(.5,1)
\psline{->}(.8,0)(.5,.6)
\psline{-}(.2,0)(.5,.6)
\pscurve[border=2pt]{-*}(.1,1)(.3,.3)(.7,.2)
\end{pspicture}
\;\; \mbox{and} \;\;
\begin{pspicture}[.2](0,-.2)(1.6,1)
\psset{xunit=1.2cm,yunit=1.2cm}
\rput[r](.05,1){$a$}
\rput[r](.2,0){$b$}
\rput[l](.8,0){$c$}
\rput[l](.55,1.05){$d$}
\rput[l](.55,.4){$v_k$}
\psline{<-}(.5,.35)(.5,1)
\psline{-*}(.75,0)(.5,.35)
\psline{-}(.25,0)(.5,.35)
\pscurve[border=2pt]{-*}(.1,1)(.2,.75)(.7,.75)(.5,.85)
\end{pspicture}
$$

Then the sum $\sum_{\tilde{\Gamma}; \tilde{\Gamma} \in \CE(\Gamma;B)}[\tilde{\Gamma}]$ is zero thanks to IHX. \index{N}{IHX}
\eop

\begin{lemma}
\label{lemfacetot}
For any $\Gamma \in \CE_n$,
$I_{{\Gamma},F(V)}=0$ if $\hat{\omega}^{(i)}_M$ and $\omega^{(i)}_M$ are fundamental with respect to $\tau_M$,  or if $\hat{\omega}^{(i)}_M$ and $\omega^{(i)}_M$ coincide on $ST(B_M)$.
\end{lemma}
In the first case, the face $F(V)$ is identified via $\tau_M$ to $\breve{S}_V(\RR^3) \times (M \setminus \infty)$, and the form $\tilde{\Omega}$ to be integrated can be pulled-back through the projection
onto the fiber. In the second case, for any $j \in \{1,2,\dots, 3n\}$, $p_j$ maps $F(V)$ into
$\partial C_2(M)$, and therefore $P(\Gamma)^{\ast}(\tilde{\Omega})=0$ on $F(V)$
thanks to Lemma~\ref{lemnolosseta}.
\eop

This ends the proof of 
Proposition~\ref{propunb}, and hence the proofs of Proposition~\ref{propun}
and \ref{propzad}.

Since for any admissible form $\omega_M$ on $C_2(M)$, $z_n(\omega_M)$ only depends on 
the restriction of $\omega_M$ to $ST(B_M)$, $z_n(\omega_M)$ will also be denoted by $z_n(\omega_{M|{ST(B_M)}})$. 

\begin{proposition}
\label{propunbb}
For any admissible form $\omega$ on $C_2(M)$ and for any one-form $\eta$ on $C_2(M)$ that 
reads $p_M(\tau_M)^{\ast}(\eta_{S^2})$ on $\partial C_2(M) \setminus ST(B_M)$, for some one-form $\eta_{S^2}$ on $S^2$ and for some trivialisation $\tau_M$ that is standard near $\infty$,
$$z_n(\omega +d \eta)-z_n(\omega)$$
$$=\sum_{\Gamma \in {\cal E}_n} \int_{F(V)}\sum_{i=1}^{3n}\left(\bigwedge_{j=1}^{i-1}P_j({\Gamma})^{\ast}(\omega) \wedge P_i({\Gamma})^{\ast}(\eta) \wedge \bigwedge_{j=i+1}^{3n}P_j({\Gamma})^{\ast}(\omega+d\eta)\right)[{\Gamma}].$$
\end{proposition}
\bp
Indeed, according to Proposition~\ref{propunb}, $$z_n(\left(\bigwedge_{j=1}^{i-1}P_j^{\ast}(\omega) \wedge \bigwedge_{j=i}^{3n}P_j^{\ast}(\omega+d\eta)\right))-z_n(\left(\bigwedge_{j=1}^{i}P_j^{\ast}(\omega) \wedge \bigwedge_{j=i+1}^{3n}P_j^{\ast}(\omega+d\eta)\right))=$$
$$= \sum_{\Gamma \in {\cal E}_n}
\int_{F(V)}\left(\bigwedge_{j=1}^{i-1}P_j({\Gamma})^{\ast}(\omega) \wedge P_i({\Gamma})^{\ast}(\eta) \wedge \bigwedge_{j=i+1}^{3n}P_j({\Gamma})^{\ast}(\omega+d\eta)\right)[{\Gamma}],$$
and the above statement is nothing but the sum over the $i$ in $\{1,\dots, 3n\}$ of these equalities.
\eop


\subsection{Forms over $S^2$-bundles}
\label{subsproofdeux}
\noindent{\sc Proof of Proposition~\ref{propdeux}:}
According to Proposition~\ref{propconffacedeux},
the codimension one faces of $S_V(E_1)$ are the fibered spaces over $S^4$
with fibers $f(B)(p^{-1}(x))$, for all the strict subsets $B$ of $V$ with cardinality at least $2$.
Then the independence of $\omega_T$ is proved as in the previous subsection,
using lemmas similar to Lemmas~\ref{lemdiscon}, \ref{lemedge}, \ref{lemsym}, \ref{lemihx}
that treat all the possible faces, and Proposition~\ref{propdeux} is proved.
 \eop

Note that the proofs of these lemmas in fact show that
the image of $S_{2n}(E_1)$ under $\sum_{\Gamma \in \CE_n}P(\Gamma)[\Gamma]$
is a 
cycle whose homology class is in $H_{6n}(S_2(E_1)^{3n};\CA(\emptyset)) $ even if $\CA(\emptyset)$ is defined with integral coefficients. (Its boundary $$\sum_{(\Gamma,F);\Gamma \in \CE_n, F \in \partial_1(S_{V}(E_1))}[P(\Gamma)(F)][\Gamma]$$ vanishes algebraically. ) Then $2^{3n}(2n)!(3n)!\xi_n$ \index{N}{ksin@$\xi_n$} is just
the evaluation of $\bigwedge_{i=1}^{3n}p_i^{\ast}[\omega_T]$ at the class of this cycle.

More generally, we have the following proposition:

\begin{proposition}
\label{propbunind}
Let $E$ be an $\RR^3$-bundle over a base $W$ that is an oriented four-dimensional manifold.
Let $\eta$ denote a one-form on $S_2(E)$ and let $\omega$ denote a closed
two-form on $S_2(E)$.
Let \index{N}{znEomega@$z_n(E;\omega)$}
$$z_n(E;\omega)=\sum_{\Gamma \in {\cal E}_n} \int_{\breve{S}_{2n}(E)}\bigwedge_{i=1}^{3n}P_i({\Gamma})^{\ast}(\omega)[{\Gamma}].$$
Then $z_n(E;\omega+d\eta)-z_n(E;\omega)=\delta_n(E;\omega,\eta)$
with $\delta_n(E;\omega,\eta)=$ 
$$\sum_{\Gamma \in {\cal E}_n} \int_{\breve{S}_{2n}(E_{|\partial W})}\sum_{i=1}^{3n}\left(\bigwedge_{j=1}^{i-1}P_j({\Gamma})^{\ast}(\omega) \wedge P_i({\Gamma})^{\ast}(\eta) \wedge \bigwedge_{j=i+1}^{3n}P_j({\Gamma})^{\ast}(\omega+d\eta)\right)[{\Gamma}].$$
\end{proposition}
\bp
The contributions of the faces coming from the boundary
of $S_{2n}(\RR^3)$ cancel as in the above case and we are left with the contributions
coming from the boundary of $W$. \eop

\begin{lemma}
\label{lemdepboun}
Let $W$ be a connected oriented compact four-dimensional manifold, let $E$ be the trivial $\RR^3$-bundle $E=W \times \RR^3$, and let $\omega$ denote a closed
two-form on $S_2(E)$. If the inclusion induces an injection from $H^2(W)$ to $H^2(\partial W)$ and a surjection from $H^1(W)$ to $H^1(\partial W)$, then
$z_n(E;\omega)$ \index{N}{znEomega@$z_n(E;\omega)$}
only depends on the restriction of $\omega$ to $\partial W \times S^2$.
\end{lemma}
\bp Indeed, a closed form $\omega^{\prime}$ that coincides with $\omega$ on $\partial W \times S^2$ would read $(\omega+d\eta)$ for some one-form $\eta$ whose restriction
to the boundary $\partial W \times S^2$ is closed and may be extended to $W \times S^2$ as a closed form $\eta^{\prime}$. Thus, $\omega^{\prime}=\omega +d(\eta -\eta^{\prime})$ and since $(\eta -\eta^{\prime})$ vanishes on $\partial W \times S^2$, Proposition~\ref{propbunind} guarantees that $z_n(E;\omega)=$
$z_n(E;\omega^{\prime})$.
\eop

Here, a {\em bundle morphism\/} $\psi$ from an $\RR^3$-bundle $E$ to another one $E^{\prime}$ will always restrict to an isomorphism from a fiber of $E$
to a fiber of $E^{\prime}$. Such a bundle morphism induces bundle morphisms
that are still denoted by $\psi$ from $S_n(E)$ to $S_n(E^{\prime})$ for 
every $n$.
Note that such a bundle morphism of $\RR^3$-bundles is determined by 
$\psi: S_2(E) \longrightarrow S_2(E^{\prime})$ 
up to a multiplication by a function from the base of $E$ to $\RR$, that preserves all the maps $\psi: S_n(E) \longrightarrow S_n(E^{\prime})$.

\begin{lemma}
\label{lemtranspform}
Let $E$ be an $\RR^3$-bundle over a base $W$ that is an oriented four-dimensional manifold and let $\omega$ denote a closed
two-form on $S_2(E)$.
Assume that there exist a bundle morphism $\psi$ from $E$ to an $\RR^3$-bundle $E(X)$ over a base $X$, and a closed
two-form $\omega(X)$ on $S_2(E(X))$ such that $\omega=\psi^{\ast}(\omega(X))$.
If $X$ is a manifold of dimension $< 4$, 
then $z_n(E;\omega)=0$, \index{N}{znEomega@$z_n(E;\omega)$}
and if 
$\psi$ is an orientation-preserving diffeomorphism, then $z_n(E;\omega)=z_n(E(X);\omega(X))$.
\end{lemma}
\bp
Indeed, since the maps $\psi$ commute with the $P_i(\Gamma)$, 
$$\int_{\breve{S}_{2n}(E)}\bigwedge_{i=1}^{3n}P_i({\Gamma})^{\ast}(\psi^{\ast}(\omega(X)))
=\int_{\breve{S}_{2n}(E)}\psi^{\ast}\left(\bigwedge_{i=1}^{3n}P_i({\Gamma})^{\ast}(\omega(X))\right).$$
Therefore, if $X$ is of dimension $< 4$, the dimension of $\breve{S}_{2n}(E(X))$
is less than the dimension of $\breve{S}_{2n}(E)$ and the integral vanishes.
If $\psi$ is an orientation-preserving diffeomorphism, then it induces an 
orientation-preserving diffeomorphism from $\breve{S}_{2n}(E)$ to  $\breve{S}_{2n}(E(X))$.
\eop

\subsection{The dependence on the trivializations}
\label{subdeptriv}

The closure of the face $F(V)$ \index{N}{FV@$F(V)$} in $C_V(M)$ is diffeomorphic via $\tau_M$ to $C_1(M) \times S_V(\RR^3)$. When $(S_V=S_V(\RR^3))$ is oriented as in Subsection~\ref{subfra}, the involved  diffeomorphism preserves the orientation.
Since 
any ordered pair $P$ included into $V=\{1,2,\dots, 2n\}$ gives rise to a restriction map
from $S_{2n}=S_V$ to $S_P=S_2=S^2$,
 any edge-oriented labelled graph $\Gamma$ again induces a smooth map 
$$P(\Gamma):S_{2n}(\RR^3) \longrightarrow (S^2)^{3n}$$ whose $i^{th}$ projection $P_i(\Gamma)$ is the map associated to the edge labelled by $i$.

The key-proposition to study how $z_n(\omega_M)$ depends on the restriction of $\omega_M$ to $ST(B_M)$
is the following one.

\begin{proposition}
\label{propkeytriv}
Let $\omega_0$ and $\omega_1$ be two admissible two-forms on $C_2(M)$
 that coincide on $\partial C_2(M) \setminus ST(B_M)$.
Let $\tau_M$ be a trivialisation of 
$(M \setminus \infty)$ that is standard near $\infty$. 
Identify $ST(B_M)$ to $B_M \times S^2$ with respect to $\tau_M$.
Then there exists a closed two-form $\omega$ on $ [0,1] \times B_M \times S^2 $
such that
\begin{itemize}
\item
$\omega$ coincides with 
$\pi_{B_M \times S^2}^{\ast}\omega_1$ on $(\{1\} \times B_M  \cup [0,1] \times \partial B_M) \times S^2$ and, 
\item
$\omega$ coincides with 
$\omega_0$ on  $\{0\} \times B_M \times S^2 $,
\end{itemize}
and, for any such two-form $\omega$,
$$z_n(\omega_1)-z_n(\omega_0)=z_n([0,1] \times B_M \times \RR^3;\omega)$$
where $$z_n([0,1] \times B_M \times \RR^3;\omega)
=\sum_{\Gamma \in {\cal E}_n} \int_{[0,1] \times B_M \times S_V(\RR^3)} \bigwedge_{i=1}^{3n}P_i(\Gamma)^{\ast}(\omega)[\Gamma].$$
\end{proposition}
\bp
First, the two-form $\omega$ exists because 
the restriction induces an isomorphism from $H^2([0,1] \times B_M \times S^2;\RR)$ to $H^2(\partial([0,1] \times B_M \times S^2);\RR)$. See Lemma~\ref{lemdr}.

Next, $z_n([0,1] \times B_M \times \RR^3;\omega)$ is independent of
the chosen closed extension $\omega$ by Lemma~\ref{lemdepboun}.

Now, $\left(z_n(\omega_0)+z_n([0,1] \times B_M \times \RR^3;\omega)\right)$
is independent of $\omega_0$ because $[0,1] \times B_M \times S_V(\RR^3)$
can be glued to $C_V(M)$ along the closure of $F(V)_{|B_M}$ that is identified to $\{0\} \times B_M \times S_V(\RR^3)$ via $\tau_M$. The details of this argument can be written as follows.
Let $\tilde{\omega}_0$ be another admissible form on $\partial C_2(M)$ that coincides with $\omega_1$ outside $ST(B_M)$, and let 
$\tilde{\omega}=\omega +d \eta$ be a closed two-form on 
$[0,1] \times B_M \times S^2$ that
coincides\\
with $\tilde{\omega}_0$ on $\{0\} \times B_M \times S^2 $, and \\
with $\pi_{B_M \times S^2}^{\ast}\omega_1$ on $(\{1\} \times B_M \cup [0,1] \times \partial B_M) \times S^2 $. \\
We assume that $\eta$ vanishes
on $(\{1\} \times B_M  \cup [0,1] \times \partial B_M) \times S^2$  without loss because the $H^1$ of this space is trivial.
In particular, according to Proposition~\ref{propbunind},
$$z_n([0,1] \times B_M \times \RR^3;\tilde{\omega}) - z_n([0,1] \times B_M \times \RR^3;\omega)=$$
$$-\sum_{\Gamma \in {\cal E}_n} \int_{0 \times B_M \times S_V(\RR^3)}\sum_{i=1}^{3n}\left(\bigwedge_{j=1}^{i-1}P_j({\Gamma})^{\ast}(\omega_0) \wedge P_i({\Gamma})^{\ast}(\eta) \wedge \bigwedge_{j=i+1}^{3n}P_j({\Gamma})^{\ast}(\tilde{\omega}_0)\right)[{\Gamma}].$$
Similarly, $\tilde{\omega}_0$ and $\omega_0$  extend to $C_2(M)$ as $\tilde{\omega}_M$ and ${\omega}_M$, respectively, and there exists a one-form
$\eta^{\prime}$ on $C_2(M)$ such that $\tilde{\omega}_M=\omega_M +d \eta^{\prime}$
where $\eta^{\prime}$ vanishes on $\partial C_2(M) \setminus ST(B_M)$, and coincides with $\eta$ on $ST(B_M)$. Then 
according to Proposition~\ref{propunbb},
$$z_n(\tilde{\omega}_M) - z_n(\omega_M)= - (z_n([0,1] \times B_M \times \RR^3;\tilde{\omega}) - z_n([0,1] \times B_M \times \RR^3;\omega)).$$

In particular, when $\tilde{\omega}_0=\omega_1$, we can choose the extension $\tilde{\omega}=\pi_{B_M \times S^2}^{\ast}(\omega_1)$ where 
$$\pi_{B_M \times S^2}: [0,1] \times B_M \times \RR^3 \longrightarrow \{1\} \times B_M \times \RR^3$$
is the natural bundle morphism
and we
have $$z_n(\omega_0)+z_n([0,1] \times B_M \times \RR^3;\omega)=$$
$$=z_n(\omega_1)+z_n([0,1] \times B_M \times \RR^3;\pi_{B_M \times S^2}^{\ast}(\omega_1))$$
where, according to Lemma~\ref{lemtranspform}, $$z_n([0,1] \times B_M \times \RR^3;\pi_{B_M \times S^2}^{\ast}(\omega_1))=0.$$

\eop

\begin{lemma}
If $\tilde{\tau}_M$ is a trivialization homotopic to $\tau_M$, then
$z_n(\tau_M)=z_n(\tilde{\tau}_M)$.
\end{lemma}
\bp When $\tilde{\tau}_M$ is homotopic to $\tau_M$,
there exists $g:  [0,1] \times B_M  \longrightarrow GL^+(\RR^3)$ such that
$g$ maps a neighborhood of $([0,1] \times(B_M \setminus B_M(1)) \cup \{1\} \times B_M )$ to $1$, and, if $\tau_M(v \in T_m M)=(m,u \in \RR^3)$, then
$\tilde{\tau}_M(v \in T_m M)=(m,g(0,m)(u))$. The map $g$ induces the bundle-morphism
$$\begin{array}{llll}\phi(g):&[0,1] \times B_M \times \RR^3 &\longrightarrow &\RR^3\\
&(t,m,v)& \mapsto & g(t,m)(v).
\end{array}$$
such that $\phi(g)^{\ast}(\omega_{S^2})$ satisfies the hypotheses of Proposition~\ref{propkeytriv}, and 
$$z_n(\tau_M)-z_n(\tilde{\tau}_M)=z_n([0,1] \times B_M \times S^2;\phi(g)^{\ast}(\omega_{S^2})).$$
Then thanks to Lemma~\ref{lemtranspform}, the right-hand side vanishes.
\eop

This lemma concludes the proof of the first part of Proposition~\ref{proptrois}.

\begin{lemma}
\label{lemtrivind}
Let $G: M \setminus \infty \longrightarrow GL^+(\RR^3)$ map $(M \setminus \infty) \setminus B_M(1)$ to $1$. Then 
$z_n(\psi(G) \circ \tau_M)-z_n(\tau_M)$ is independent of $\tau_M$.
\end{lemma}
\bp
Let $\tilde{\tau}_M$ be another trivialisation. Then, there exists 
$g:(M \setminus \infty) \longrightarrow GL^+(\RR^3)$
such that $\tilde{\tau}_M=\psi(g) \circ \tau_M $.
The map $g$ induces automorphisms $\psi(g)$ on all $((M \setminus \infty) \times S_V)$.
Furthermore on $B_M \times S^2$,
$p_M(\tau_M) =p_{S^2}$,
$p_M(\tilde{\tau}_M)= p_{S^2} \circ \psi(g) $,
$p_M( \psi(G) \circ \tau_M)= p_{S^2} \circ \psi(G)$, and
$p_M(\psi(G)  \circ \tilde{\tau}_M)=p_{S^2} \circ \psi(G) \circ \psi(g)$.
Thus, when $\omega$ is suitable to compute $\left(z_n(\tau_M)-z_n(\psi(G)  \circ \tau_M)\right)$, according to Proposition~\ref{propkeytriv}, $\psi(g)^{\ast}(\omega)$ is suitable to compute $\left(z_n(\tilde{\tau}_M)-z_n( \psi(G)  \circ \tilde{\tau}_M)\right)$, and since
this amounts to pull-back $\bigwedge_{i=1}^{3n}P_i(\Gamma)^{\ast}(\omega)$ by the orientation-preserving diffeomorphism $\psi(g)$ acting on $ I \times B_M \times S_V$, it does 
not change the integrals. Therefore, 
$$z_n(\tau_M)-z_n(\psi(G)  \circ \tau_M)=z_n(\tilde{\tau}_M)-z_n( \psi(G)  \circ \tilde{\tau}_M)
$$ 
and we are done.
\eop

The above lemma allows us to define 
$$z_n^{\prime}(G)=z_n(\psi(G) \circ \tau_M)-z_n(\tau_M)$$
for any $G: M \setminus \infty \longrightarrow GL^+(\RR^3)$ that maps $(M \setminus \infty) \setminus B_M(1)$ to $1$.

\begin{lemma}
\label{lemvalrho}
If $G$ maps the complement of a ball $B^3$ to the identity, and if $G$ is homotopic to $\rho$ on the quotient of this 3-ball by its boundary, then
$$z_n^{\prime}(G)=\delta_n.$$
\end{lemma}
\bp
Indeed, in this case, there exists a two-form $\omega$ on $ [0,1] \times  B_M\times S^2 $
that coincides with $p_{S^2}^{\ast}(\omega_{S^2})$ near $\{1\} \times B_M \times S^2 $ and $[0,1] \times (B_M \setminus B^3) \times S^2$ and that coincides with 
$(p_{S^2} \circ \psi(\tilde{\rho})^{-1})^{\ast}(\omega_{S^2})$ near $\{0\} \times B^3 \times S^2 $ where $\tilde{\rho}$ denotes the restriction of $G$ to $B^3$ that is homotopic to $\rho$. Then
$$z_n(\tau_M)-z_n(\psi(G)^{-1} \circ \tau_M)=\sum_{\Gamma \in {\cal E}_n} \int_{[0,1] \times B^3 \times S_V(\RR^3)}\bigwedge_{i=1}^{3n}P_i(\Gamma)^{\ast}(\omega)[\Gamma]$$
since the forms vanish on $[0,1] \times (B_M \setminus B^3) \times S_V$.
Now, view the bundle $E_1$ of Subsection~\ref{subfra} as \index{N}{Eone@$E_1$}
$$E_1 \cong B^4 \times \RR^3 \cup_{(\partial B^4=S^3=B^3 \cup_{S^2} -B^3) \times S^2} -B^4 \times \RR^3$$
where
$(x,y)$ of the first copy $\partial B^4 \times \RR^3$ is identified with
$(x,y)$ of the second copy  if $x$ is in $(-B^3)$,
 and with $\psi(\tilde{\rho})(x,y)$ otherwise. Then $\omega$ can be extended by 
$p_{S^2}^{\ast}(\omega_{S^2})$ outside $ [0,1] \times (B^3) \times S^2 \subset -B^4 \times S^2$ on $S_2(E_1)$. The integrals over 
$S_V\left(E_1 \setminus p^{-1}\left(([0,1] \times B^3) \subset -B^4\right)\right)$ are zero.
Therefore, $$z_n(\tau_M)-z_n(\psi(G)^{-1} \circ \tau_M)=\sum_{\Gamma \in {\cal E}_n} \int_{S_V(E_1)} \bigwedge_{i=1}^{3n}P_i(\Gamma)^{\ast}(\omega)[\Gamma]=\delta_n.$$
Now,  $\omega$  represents 
the Thom class of $E_1$, and we conclude with the help of Lemma~\ref{lemtrivind}.
\eop

\subsection{More on trivialisations of $3$-manifolds}
\label{submoretriv}

Let us now recall some more standard facts about homotopy classes of orientation-respecting trivialisations of $3$-manifolds.

Fix a trivialisation $\tau_M$ of $T(M \setminus \infty)$ that is standard near $\infty$. 
Any other such will read $\psi(G) \circ \tau_M$ for a unique
$G: ((M \setminus \infty), M \setminus (\infty \cup B_M(1)))  \longrightarrow (GL^+(\RR^3),1)$, with the notation of Proposition~\ref{proppont}.
Then $(G \mapsto \psi(G) \circ \tau_M)$ induces a (non-canonical) bijection between the homotopy classes of trivialisations of $T(M \setminus \infty)$ that are standard near $\infty$, and the homotopy classes of maps from $(M,M \setminus B_M(1))$ to $(GL^+(\RR^3),1)$.
This latter set is denoted by $[(M,M \setminus B_M(1)),(GL^+(\RR^3),1)]$.
It canonically coincides with $[(M,M \setminus B_M(1)),(SO(3),1)]$.

Let $\Gamma$ be a topological group, and let $X$ be a topological space.
Define the product of two maps $f$ and $g$ from $X$ to $\Gamma$ as
$$\begin{array}{llll}fg: &X & \longrightarrow &\Gamma\\
&x &\mapsto & f(x)g(x).\end{array}$$ 
This product induces a group structure on the set $[X,\Gamma]$ of homotopy classes
of maps from $X$ to $\Gamma$. When $X=M$, this product induces a group structure on $[(M,M \setminus B_M(1)),(GL^+(\RR^3),1)]$. Recall the easy lemma.

\begin{lemma}
\label{prodpi}
The usual product of $\pi_n(\Gamma)$ coincides with the product induced by the multiplication in $\Gamma$ (defined above with $X=S^n$).
\end{lemma}
\eop

Let $G_M(\rho):M \longrightarrow SO(3)$ \index{N}{GM@$G_M(\rho)$} 
\index{N}{rho@$\rho$} be a map that sends the complement of a ball $B^3 \subset B_M(1)$ to the identity, and that is homotopic to $\rho$ on the quotient of this 3-ball by its boundary. Note that all such maps induce the same element
$[G_M(\rho)]$ in $[(M,M \setminus B_M(1)),(GL^+(\RR^3),1)]$.

The elements of $[(M,M \setminus B_M(1)),(SO(3),1)]$ have a well-defined degree
that is the degree of one of their representative from $M$ to $SO(3)$.

\begin{lemma}
\label{lempreptrivun}
Let $M$ be a closed oriented $3$-manifold.
\begin{enumerate}
\item Any map $G$ from $(M,M \setminus B_M(1))$ to $(SO(3),1)$, such that 
$$\pi_1(G): \pi_1(M) \longrightarrow \pi_1(SO(3))\cong \ZZ/2\ZZ$$
is trivial, belongs to the subgroup $<[G_M(\rho)]>$ of 
$[(M,M \setminus B_M(1)),(SO(3),1)]$
generated by $[G_M(\rho)]$. \index{N}{GM@$G_M(\rho)$}
\item For any $[G] \in [(M,M \setminus B_M(1)),(SO(3),1)]$, 
$$[G]^2 \in <[G_M(\rho)]>.$$
\item The group $[(M,M \setminus B_M(1)),(SO(3),1)]$ is abelian.
\item The degree is a group homomorphism from $[(M,M \setminus B_M(1)),(SO(3),1)]$ to $\ZZ$.
\item The morphism $$\begin{array}{llll}\frac{\mbox{deg}}{2}:&[(M,M \setminus B_M(1)),(SO(3),1)]\otimes_{\ZZ} \QQ &\longrightarrow
&\QQ[G_M(\rho)]\\
&[g] \otimes 1 &\mapsto &\frac{\mbox{deg}(g)}{2}[G_M(\rho)]\end{array}$$
is an isomorphism.
\end{enumerate}
\end{lemma}
\bp Assume that $\pi_1(G)$ is trivial. Choose a cell decomposition of $B_M$ with respect to its boundary
with no zero-cell, only one three-cell, one-cells and two-cells. Then after a homotopy, we may assume that $G$ maps the one-skeleton of $B_M$ to $1$. Next, since $\pi_2(SO(3)) = 0$, we may assume that $G$ maps the two-skeleton of $B_M$ to $1$, and therefore that $G$ maps the exterior of some $3$-ball to $1$. Now $G$ becomes a map from $B^3/\partial B^3=S^3$ to $SO(3)$, and its homotopy class is $k[\rho]$ in $\pi_3(SO(3))=\ZZ[\rho]$, where $(2k)$ is the degree of the map $G$ from $S^3$ to $SO(3)$. Therefore $G$ is homotopic to $G_M(\rho)^k$, and this proves the first assertion.

Since $\pi_1(G^2)=2\pi_1(G)$ is trivial, the second assertion follows.

For the third assertion, first note that $[G_M(\rho)]$ belongs to the center of
$[(M,M \setminus B_M(1)),(SO(3),1)]$ because it can be supported in a small ball disjoint 
from the support (preimage of $SO(3) \setminus \{1\}$) of a representative of any other element. Therefore, according to the second assertion any square will be in the center. Furthermore, since any commutator induces the trivial map on $\pi_1(M)$, any commutator is in $<[G_M(\rho)]>$.
In particular, if $f$ and $g$ are elements of $[(M,M \setminus B_M(1)),(SO(3),1)]$, 
$$(gf)^2=(fg)^2=(f^{-1}f^2g^2f)(f^{-1}g^{-1}fg)$$ 
where the first factor equals $f^2g^2=g^2f^2$. Exchanging $f$ and $g$ yields
$f^{-1}g^{-1}fg=g^{-1}f^{-1}gf$. Then the commutator that is a power of $[G_M(\rho)]$ has a vanishing square, and thus a vanishing degree. Then it must be trivial.

For the fourth assertion, it is easy to see that $\mbox{deg}(fg)=\mbox{deg}(f)+\mbox{deg}(g)$ when $f$ or $g$
is a power of $[G_M(\rho)]$, and that $\mbox{deg}(f^k)=k\mbox{deg}(f)$ for any $f$. In general,
$\mbox{deg}(fg)=\frac{1}2\mbox{deg}((fg)^2)=\frac{1}2\mbox{deg}(f^2g^2)=\frac{1}2\left( \mbox{deg}(f^2)+\mbox{deg}(g^2)\right)$, and the fourth assertion is proved.

In particular, 
$$\begin{array}{llll}\frac{\mbox{deg}}{2}:&[(M,M \setminus B_M(1)),(SO(3),1)]\otimes_{\ZZ} \QQ &\longrightarrow
&\QQ[G_M(\rho)]\\
&[g] \otimes 1 &\mapsto &\frac{\mbox{deg}(g)}{2}[G_M(\rho)]\end{array}$$
is an isomorphism, and the last assertion follows, too.
\eop

\begin{lemma}
The map $z_n^{\prime}$ from $[(M,M \setminus B_M(1)),(SO(3),1)]$ to $\CA_n(\emptyset)$ is a group homomorphism.
\end{lemma}
\bp According to Lemma~\ref{lemtrivind},
$z_n^{\prime}(fg)=z_n(\psi(f) \psi(g) \tau_M) -z_n( \psi(g) \tau_M)+z_n( \psi(g) \tau_M) -z_n(\tau_M)=z_n^{\prime}(f)+z_n^{\prime}(g)$. \eop

This lemma, together with Lemma~\ref{lemvalrho} that asserts that $z_n^{\prime}(G_M(\rho))=\delta_n$ and Lemma~\ref{lempreptrivun},
concludes the proof of Proposition~\ref{proptrois}. \eop

\subsection{Proof of Proposition~\ref{proppont}}
\label{subproofpont}


Recall that the first \indexT{Pontryagin class} $p_1(W)$ of a closed oriented 4-manifold $W$ is the obstruction to trivialise the complexification of its tangent bundle. It is defined like in Subsection~\ref{subpont}. See also \cite{milnorsta}.
According to \cite[Example 15.6]{milnorsta}, $p_1(\CC P^2)=3$. \index{N}{pone@$p_1$}
We shall use the following Rohlin theorem that compares the two cobordism invariants of closed 4-manifolds.

\begin{theorem}[Rohlin]
When $W$ is a closed oriented 4-manifold,
$$p_1(W)=3\mbox{signature(W)}.$$
\end{theorem}

\begin{lemma}
\label{lempunrel}
Let $M$ be a $\QQ$-sphere. Let $\tau_M$ be a trivialisation of $T(M \setminus \infty)$ that is trivial near $\infty$. Let $W$ and $W^{\prime}$ be two cobordisms between $B^3(3)$ and $B_M$. Then \index{N}{pone@$p_1$}
$$p_1(W; \tau(\tau_M))-p_1(W^{\prime}; \tau(\tau_M))= 3\left(\mbox{signature}(W) -\mbox{signature}(W^{\prime})\right).$$
\end{lemma}
\bp 
Let $N(\partial W)$ be a regular neighborhood of $\partial W$ in $W$, or in $W^{\prime}$.
Let $\tau$ be a trivialisation of $TW \otimes \CC$ defined in $N(\partial W)$.
Set $\tilde{W}=W \setminus \mbox{Int}(N(\partial W))$, and $\tilde{W}^{\prime}=W^{\prime} \setminus \mbox{Int}(N(\partial W))$. 
Then 
$$p_1(W;\tau)-p_1(W^{\prime};\tau)=p_1(\tilde{W};\tau)-p_1(\tilde{W}^{\prime};\tau)$$
does not depend on the trivialisation $\tau$ and equals 
$p_1(\tilde{W}\cup_{\partial \tilde{W}} - \tilde{W}^{\prime})$.
According to Rohlins's theorem, this is $3\; \mbox{signature}(\tilde{W}\cup_{\partial \tilde{W}} - \tilde{W}^{\prime})$,
where $\tilde{W}\cup_{\partial \tilde{W}} - \tilde{W}^{\prime}$ is homeomorphic to $W \cup_{\partial W}(- W^{\prime})$ and $\partial W$ is homeomorphic to $M$.

Since $M$ is a $\QQ$-sphere, the Mayer-Vietoris sequence makes clear that $$H_2(W \cup_M(- W^{\prime});\RR)=H_2(W;\RR) \oplus H_2(W^{\prime};\RR),$$
and it is easy to see that
$$\mbox{signature}(W \cup_M (-W^{\prime}))= \mbox{signature}(W) -\mbox{signature}(W^{\prime}),$$
and to conclude.
\eop

In particular, the definition of $p_1$ does not depend
on the chosen $4$-cobordism $W$ with signature $0$. It is clear that $p_1(\tau_M)$ only depends on the homotopy class of $\tau_M$. Proposition~\ref{proppont} is now the direct consequence of Lemmas~\ref{lempunind},
\ref{lemvarpun} and \ref{lembijpun} below.

Let $\KK= \RR$ or $\CC$. Let $n \in \NN$. The stabilisation maps induced
by the inclusions
$$\begin{array}{llll}i: & GL(\KK^n) & \longrightarrow & GL(\KK \oplus \KK^n)\\
& g & \mapsto & (i(g): (x,y) \mapsto (x,g(y))\end{array}$$
will be denoted by $i$. The $\KK$ (euclidean or hermitian) oriented vector space with the direct orthonormal basis $(v_1, \dots, v_n)$ will be denoted by $\KK<v_1, \dots, v_n>$.
The inclusions $SO(n) \subset SU(n)$ will be denoted by $c$.
The projection from $SO(\RR^4=\RR<1,i,j,k>)$ to $S^3$ that maps 
$g$ to $g(1)$ is denoted by $p$.
In particular, the long exact sequence associated to the fibration
$SO(3) \hookrightarrow SO(4) \hfl{p} S^3$ gives rise 
to the exact sequence \index{N}{rho@$\rho$}
$$\pi_3(SO(3))=\ZZ[\rho] \hfl{i_{\ast}} \pi_3(SO(4)) \hfl{p_{\ast}} \pi_3(S^3)=\ZZ[\mbox{Id}] \longrightarrow \{0\}$$
Let $m_r$ \index{N}{mr@$m_r$}
denote the map from $S^3=S(\HH)$ to $SO(\RR^4=\HH)$ be induced by the
right-multiplication. When $v \in S^3$ and $x \in \HH$, $m_r(v)(x)=x.v$.
Define a section $\sigma$ of $p_{\ast}$, by setting
$$\sigma([\mbox{Id}])=[m_r].$$
In particular, $\pi_3(SO(4))$ is generated by $i_{\ast}([\rho])$ and $[m_r]$.

Let $m^{\CC}_r$ \index{N}{mrC@$m^{\CC}_r$}
denote the homeomorphism from $S^3=S(\HH)$ to $SU(\CC^2=\CC<1,j>= \HH)$ be induced by the
right-multiplication. When $v \in S^3$ and $x \in \HH$, $m^{\CC}_r(v)(x)=x.v$.
$$m^{\CC}_r(z_1 + z_2 j)=\left[ \begin{array}{cc} z_1 &-\overline{z}_2\\ z_2 &\overline{z}_1\end{array}\right].$$
$$\pi_3(SU(2))=\ZZ[m^{\CC}_r]$$
Finally recall that $i^n_{\ast}: \pi_3(SU(2)) \longrightarrow \pi_3(SU(n+2))$
is an isomorphism for any natural number $n$, 
and in particular, that \index{N}{itwo@$i^2(m^{\CC}_r)$}
$$\pi_3(SU(4))=\ZZ[i^2(m^{\CC}_r)].$$

The following lemma determines the map $$c_{\ast}: \pi_3(SO(4)) \longrightarrow \pi_3(SU(4)).$$
\begin{lemma}
\label{lempitroissoquatre}
$$c_{\ast}([m_r])=2[i^2(m^{\CC}_r)].$$
$$c_{\ast}(i_{\ast}([\rho]))=-4[i^2(m^{\CC}_r)].$$ \index{N}{rho@$\rho$}
$$\pi_3(SO(4))=\ZZ[m_r] \oplus \ZZ i_{\ast}([\rho]).$$
\end{lemma}
\bp
Let $m_{\ell}$ denote the map from $S^3=S(\HH)$ to $SO(\RR^4=\HH)$ induced by the
left-multiplication. When $v \in S^3$ and $x \in \HH$, $m_{\ell}(v)(x)=v.x$.
Let $\overline{m}_r=m_r^{-1}$. 
When $v \in S^3$ and $x \in \HH$, $\overline{m}_r(v)(x)=x.\overline{v}$.
Then in $\pi_3(SO(4))$, $$i_{\ast}([\rho])=[m_{\ell}] + [\overline{m}_r]=
[m_{\ell}] - [{m}_r],$$
thanks to Lemma~\ref{prodpi}.
Now, using the conjugacy of quaternions, $m_{\ell}(v)(x)=v.x=\overline{\overline{x}.\overline{v}}=\overline{\overline{m}_{r}(v)(\overline{x})}$.
Therefore $m_{\ell}$ is conjugated to $\overline{m}_{r}$ via 
the conjugacy of quaternions that acts on $\RR^4$ as a hyperplan symmetry.

Now, observe that since $U(4)$ is connected, the conjugacy by an element of $U(4)$ induces the identity on $\pi_3(SU(4))$. 
Thus, $$c_{\ast}([m_{\ell}])=c_{\ast}([\overline{m}_{r}])=-c_{\ast}([{m}_{r}]),$$
and 
$$c_{\ast}(i_{\ast}([\rho]))=-2 c_{\ast}([{m}_{r}]).$$
Therefore, we are left with the proof of the following sublemma
that implies that $i_{\ast}: \pi_3(SO(3)) \longrightarrow \pi_3(SO(4))$
is injective and thus, that 
$$\pi_3(SO(4))=\ZZ[m_r] \oplus \ZZ i_{\ast}([\rho]).$$

\begin{sublemma} \index{N}{itwo@$i^2(m^{\CC}_r)$}
$$c_{\ast}([m_r])=2[i^2(m^{\CC}_r)].$$
\end{sublemma}
\bp
Let $\HH + I \HH$ denote the complexification of $\RR^4= \HH= \RR<1,i,j,k>$.
Here, $\CC=\RR \oplus I\RR$.
When $x \in \HH$ and $v \in S^3$, $c(m_r)(v)(Ix)=Ix.v$, and $I^2=-1$.
Let $\varepsilon=\pm 1$, define
$$ \CC^2(\varepsilon)=\CC<\frac{\sqrt{2}}{2}(1+ \varepsilon Ii),\frac{\sqrt{2}}{2}(j+ \varepsilon Ik)>.$$
Consider the quotient $\CC^4/\CC^2(\varepsilon)$.
In this quotient, $Ii=-\varepsilon 1$, $Ik =-\varepsilon j$, and since $I^2=-1$,
$I1=\varepsilon i$ and $Ij=\varepsilon k$. Therefore this quotient is isomorphic to $\HH$ as a real vector space with its complex structure $I= \varepsilon i$.
Then it is easy to see that $c(m_r)$ maps $\CC^2(\varepsilon)$ to $0$ in this quotient. Thus $c(m_r)(\CC^2(\varepsilon)) = \CC^2(\varepsilon)$.
Now, observe that $\HH + I \HH$ is the orthogonal sum of $\CC^2(1)$ and $\CC^2(-1)$.
In particular, $\CC^2(\varepsilon)$ is isomorphic to the quotient $\CC^4/\CC^2(-\varepsilon)$ that is isomorphic to $(\HH;I= -\varepsilon i)$ and 
$c(m_r)$ acts on it by the right multiplication. Therefore,
with respect to the orthonormal basis $\frac{\sqrt{2}}{2}(1-Ii, j-Ik, 1+Ii, j+Ik )$, $c(m_r)$ reads
$$c(m_r)(z_1+z_2j) =\left[\begin{array}{cccc} 
z_1 & -\overline{z}_2 & 0 & 0\\ 
z_2 & \overline{z}_1 & 0 & 0\\
0 & 0 & \overline{z}_1=x_1-Iy_1 & -z_2\\
0 & 0 & \overline{z}_2 & z_1=x_1+Iy_1\\
\end{array} \right]$$
Therefore, the homotopy class of $c(m_r)$ (invariant under conjugacy by an element of $U(4)$) is the sum of the homotopy classes of
$$(z_1+z_2j) \mapsto \left[\begin{array}{cc} 
m^{\CC}_r & 0 \\ 
0 & 1
\end{array} \right] \;\; \mbox{and}\;\; (z_1+z_2j) \mapsto \left[\begin{array}{cc} 
1 & 0 \\ 
0 & m^{\CC}_r \circ \iota
\end{array} \right] $$
where $\iota (z_1 +z_2 j)= \overline{z}_1 +\overline{z}_2 j$.
Since the first map is conjugate by a fixed element of $SU(4)$ to $i^2_{\ast}(m^{\CC}_r)$, it is homotopic to $i^2_{\ast}(m^{\CC}_r)$, and since $\iota$ induces the 
identity on $\pi_3(S^3)$, the second map is homotopic to $i^2_{\ast}(m^{\CC}_r)$, too.
\eop

\begin{lemma}
\label{lempunind}
Consider $g: (B_M,  ]1,3] \times S^2) \longrightarrow (SO(3),1)$ and
$$\begin{array}{llll} 
\psi(g): &B_M \times \RR^3 &\longrightarrow  &B_M \times \RR^3\\
&(x,y) & \mapsto &(x,g(x)(y))\end{array}$$ then $\left(p_1(\psi(g) \circ \tau_M)-p_1(\tau_M)\right)$ is independent of $\tau_M$. 
\end{lemma}
\bp Indeed, $\left(p_1(\psi(g) \circ \tau_M)-p_1(\tau_M)\right)$ can be defined as the obstruction to extend the following trivialisation of the tangent bundle of
$[0,1] \times B_M$ restricted to the boundary. This trivialisation is
$T[0,1] \oplus \tau_M$ on $(\{0\} \times B_M) \cup ([0,1] \times \partial B_M)$ and $T[0,1] \oplus \psi(g) \circ \tau_M$ on $\{1\} \times B_M$. But this obstruction is the obstruction
to extend the map $\tilde{g}$ from $\partial([0,1] \times B_M)$ to $SO(4)$ that maps  $(\{0\} \times B_M) \cup ([0,1] \times \partial B_M)$ to $1$ and that coincides with $i(g)$ on $\{1\} \times B_M$, viewed as a map from $\partial([0,1] \times B_M)$ to $SU(4)$, on $([0,1] \times B_M)$. This obstruction that lies in 
$\pi_3(SU(4))$ since $\pi_i(SU(4))=0$, for $i<3$, is independent of $\tau_M$.
\eop

Define $p^{\prime}_1:[(M, M \setminus B_M(1)),(SO(3),1)] \longrightarrow \ZZ$
by $$p^{\prime}_1(g)=p_1(\psi(g) \circ \tau_M)-p_1(\tau_M).$$

\begin{lemma}
\label{lemvarpun}
$$p^{\prime}_1(g)=p_1(\psi(g) \circ \tau_M)-p_1(\tau_M)=-2\mbox{deg}(g).$$
\end{lemma}
\bp Lemma~\ref{lempunind} guarantees that $p^{\prime}_1$ is a group homomorphism. According to Lemma~\ref{lempreptrivun}, $p^{\prime}_1$ must read
$p^{\prime}_1(G_M(\rho)) \frac{\mbox{deg}}{2}$.
Thus, we are left with the proof that \index{N}{GM@$G_M(\rho)$} $$p^{\prime}_1(G_M(\rho))=-4.$$
Let $g=G_M(\rho)$, we can extend $\tilde{g}$ (defined in the proof of Lemma~\ref{lempunind}) by the constant map with value 1 outside
$[\varepsilon, 1] \times B^3 \cong B^4$ and, in $\pi_3(SU(4))$
$$[c(\tilde{g}^{-1}_{|\partial B^4})]=-(p_1(\psi(g) \circ \tau_M)-p_1(\tau_M))[i^2(m^{\CC}_r)].$$
Since $\tilde{g}^{-1}_{|\partial B^4}$ is homotopic to $i(\rho)^{-1}$, Lemma~\ref{lempitroissoquatre} allows us to conclude.
\eop

\begin{lemma}
\label{lembijpun}
\begin{itemize}
\item If $M$ is a given $\ZZ$-sphere, then $p_1$ \index{N}{pone@$p_1$} defines a bijection from the set of homotopy classes of trivialisations of $M$ that are standard near $\infty$ to $4 \ZZ$.
\item For any $\ZZ$-sphere $M$, for any trivialisation $\tau_M$ of $M$ that is standard near $\infty$, 
$$\left(p_1(\tau_M)-\mbox{dimension}(H_1(M;\ZZ/2\ZZ))\right) \in 2 \ZZ.$$
\end{itemize}
\end{lemma}
\bp Any closed oriented $3$-manifold $M$ bounds a $4$-dimensional manifold $W$
obtained from $B^4=[0, \varepsilon] \times B^3$ by attaching $b_2(W)$ two-handles
with even self-intersection \cite{kaplan}. We are going to prove the following sublemma.

\begin{sublemma}
There exists a trivialisation $\tau_M$ of $T(M \setminus \infty)$ that is standard near $\infty$ such that
$$p_1(W;\tau(\tau_M)) \equiv 2 b_2(W)\; \mbox{\rm mod}\; 4.$$
\end{sublemma} 
\bp For our $W$, there exists a Morse function that coincides with the projection
onto $[0,1]$ near the boundary where $W$ looks like $[0,1] \times B^3$ or
$[0,1] \times B_M$ and whose only critical points are index two critical points
that correspond to the $b_2(W)$ two-handles. Let $X$ be the gradient field of this
function that is defined outside the critical points.
Let $B^4$ be a $4$-ball of $W$ that intersects $\partial W$ along a $3$-ball $B^3 \subset M$ and that contains all the critical points.
$W \setminus B^4$ is homotopy equivalent
to $W$ and is obtained from a regular neighborhood of $(\{0\} \times B^3) \cup
(-[0,1] \times S^2)$ by attaching two-handles. 
The obstruction to extend the trivialisation $X \oplus \tau_S^3$ of $TW$ defined near $(\{0\} \times B^3) \cup
(-[0,1] \times S^2)$ 
to these handles is in $\pi_1(SO(4))=i_{\ast}(\pi_1(SO(3)))=\ZZ/2\ZZ$, it
is the self-intersection of the handles mod 2, and it vanishes.
Therefore, the trivialisation $X \oplus \tau_S^3$ of $TW$ defined near $(\{0\} \times B^3) \cup
(-[0,1] \times S^2)$ 
 extends to $(W \setminus B^4)$ as a trivialisation of the form $X \oplus \tau$.
In particular, $\tau$ provides a trivialisation $\tau_M$ on $B_M \setminus B^3$  that is standard near $\infty$, and that can be extended to $B^3$ since $\pi_2(SO(3))=\{0\}$. Now, $X \oplus \tau$ is a frame on $\partial B^4$ that is viewed as a map from $\partial B^4$ to $SO(4)$, and, in $\pi_3(SU(4))$
$$[c_{\ast}(X \oplus \tau)]=-p_1(W;\tau(\tau_M))[i^2(m^{\CC}_r)].$$
Note that $p(X \oplus \tau)=X$ and that $[X]=b_2(W)[\mbox{Id}]$ in $\pi_3(S^3)=H_3(S^3)$.
Indeed, $X$ defines a map from the complement $C$ in $B^4$ of small balls centered at
the critical points to $S^3$. In $C$,  $\partial B^4$ is homologous to the sum of the boundaries of these small balls. Therefore, when $X_{\ast}$ denotes the map
from $H_3(C)$ to $H_3(S^3)$ induced by $X$, $[X]=X_{\ast}[\partial B^4]$
is the sum of the degrees of $X$ on the boundaries of the small balls.
Since $X$ is obtained from the outward normal field by a multiplication by
a matrix with two negative eigenvalues on the boundaries of these small balls, the degree is one for all these critical points, and we have proved that $[X]=b_2(W)[\mbox{Id}]$.
Therefore,
$$(X \oplus \tau) \in \left(b_2(W)m_r \oplus i_{\ast}(\pi_3(SO(3)))\right) \subset \pi_3(SO(4)),$$
and, according to Lemma~\ref{lempitroissoquatre},
$$[c_{\ast}(X \oplus \tau)] \in 2 b_2(W)\ZZ [i^2(m^{\CC}_r)] + 4 \ZZ [i^2(m^{\CC}_r)],$$
This concludes the proof of the sublemma.
\eop

Now, it follows from Lemma~\ref{lempunrel} that
$$\begin{array}{ll}p_1(\tau_M) &=p_1(W; \tau(\tau_M))-3\,\mbox{signature}(W)\\
&\equiv 2 b_2(W) -3\,\mbox{signature}(W)\; \mbox{\rm mod} \;4.\end{array}$$
Since $M$ is a $\QQ$-sphere, $(\mbox{signature}(W)-b_2(W)) \in 2 \ZZ$, and therefore $$p_1(\tau_M)\equiv \mbox{signature}(W)\; \mbox{\rm mod} \;4.$$
When $M$ is a $\ZZ$-sphere the intersection form of $W$ is unimodular, therefore since the form is even the signature of $W$ is divisible by $8$ (see~\cite[Chap. V]{serre}), 
and $p_1(\tau_M) \in 4 \ZZ$. Thus, by Lemmas~\ref{lempreptrivun} and \ref{lemvarpun}, $p_1$ maps the homotopy classes of trivialisations of $M$ that are standard near $\infty$ onto
$4\ZZ$. These lemmas also show that $p_1$ is bijective from the set of homotopy classes of trivialisations of $M$ that are standard near $\infty$ to
$4\ZZ$.

Lemma~\ref{lemvarpun} implies that for any pair $(\tau_M, \tau_M^{\prime})$ of trivialisations of $M$ that are standard near $\infty$,
$(p_1(\tau_M)-p_1(\tau_M^{\prime}))$ is even.
Now, since the intersection matrix of $W$ mod $2$ is a presentation matrix for
$H_1(M;\ZZ/2 \ZZ)$ and since it can be written as the orthogonal sum of matrices
$\left[\begin{array}{cc}0&1\\1&0\end{array}\right]$ and a null matrix of dimension $\mbox{rank}(H_1(M;\ZZ/2 \ZZ))$, $$\mbox{signature}(W) \equiv \mbox{rank}(H_1(M;\ZZ/2 \ZZ)) \;\mbox{\rm mod}\; 2$$
and we are done. This concludes the proof of Lemma~\ref{lembijpun} and the proof of Proposition~\ref{proppont}.
\eop

\subsection{Computation of $\xi_1$} \index{N}{ksione@$\xi_1$}

We use notation introduced in Subsection~\ref{subfra}.

\begin{proposition}
\label{propcptrois}
The projective space $\CC P^3$ is homeomorphic to $-S_2(E_1)$. \index{N}{Eone@$E_1$} \index{N}{StwoEone@$S_2(E_1)$}
\end{proposition}
\begin{lemma}
The projective space $\CC P^3$ is an $S^2$-bundle over $S^4$.
\end{lemma}
\bp
Let $\HH= \CC \oplus \CC j$ be the quaternionic field, and let $\HH P^1$ be the quotient of $\HH^2 \setminus 0$ by the left multiplication by $(\HH^{\ast}=\HH \setminus \{0\})$. $$\HH P^1=S^4=\{(h_1:1); h_1 \in \HH\} \cup_{\HH^{\ast}} \{(1:h_2); h_2 \in \HH\}.$$
where $(h_1:1)=(1:h_1^{-1})$ when $h_1 \neq 0$. The complex projective space
$\CC P^3$ is the quotient of $(\CC^4 \setminus \{0\}=\HH^2 \setminus 0)$ by the  left multiplication by $\CC^{\ast} \subset \HH^{\ast}$. The projection from  $\HH^2 \setminus \{0\}$ to $S^4$ factors through $\CC P^3$
that becomes a bundle over $S^4$ with fiber $_{\CC^{\ast}}\setminus (\HH \setminus \{0\})=\CC P^1=S^2$.
\eop

\begin{lemma}
\label{lemcompglu}
Let $P_{13}=\left[\begin{array}{ccc}0&0&1\\0&1&0\\-1&0&0\end{array}\right] \in SO(3)$. Let 
$$\begin{array}{llll}g_3: &S^3 & \longrightarrow &SO(3)\\
& h_1 & \mapsto & P_{13} \rho(h_1)^{-1} P_{13}^{-1}\end{array}$$
$$\CC P^3=B^4 \times S^2 \cup_{\partial B^4 \times S^2 \stackrel{\psi(g_3)}{\rightarrow} \partial (-B^4) \times S^2} (-B^4 \times S^2)$$
\end{lemma}
\bp
Let $h_1 \in \HH^{\ast}$. The fiber of $\CC^4 \setminus \{0\}$ over $(h_1:1)$ is $\{(kh_1,k); k \in \HH^{\ast}\}$.
The fiber of $\CC^4 \setminus \{0\}$ over $(1:h_1^{-1})$ is $\{(\ell,\ell h_1^{-1}); \ell \in \HH^{\ast}\}$ with $\ell=kh_1$.
Therefore,
$$\CC P^3=B^4 \times \CC P^1 \cup_{\psi(\gamma_3)} (-B^4 \times \CC P^1)$$
where $\psi(\gamma_3)((h_1;[k]) \in \partial B^4 \times
\CC P^1)=(\overline{h}_1;\gamma_3(h_1)([k]))$
and $\gamma_3(h_1)([k])=[k.h_1]$ in $_{\CC^{\ast}} \setminus \HH^{\ast}=\CC P^1$, with $[k=z_1+z_2 j]=(z_1:z_2)$.
To express the action $g_3(h_1)$ of $\gamma_3(h_1)$ on $$S^2=\{(z \in \CC;h \in \RR);|z|^2 +h^2=1\},$$ we will use the inverse diffeomorphisms
$$\begin{array}{llll} \xi:& \CC P^1 & \longrightarrow & S^2\\
&(z_1:z_2)& \mapsto & (\frac{2 z_1 \overline{z}_2}{|z_1|^2 +|z_2|^2},
h=\frac{|z_2|^2-|z_1|^2}{|z_1|^2 +|z_2|^2})\\
 \xi^{-1}:& S^2& \longrightarrow &\CC P^1\\
&(z;h) & \mapsto & \begin{array}{ll}(z:1+h)\;&\mbox{if}\; h \neq -1\\
(1-h:\overline{z})\;&\mbox{if}\; h \neq  1 \end{array}\;\;\end{array}$$
and write
$$g_3(h_1)=\xi \circ \gamma_3(h_1) \circ \xi^{-1}.$$
Let $(z;h) \in S^2$, $h \neq -1$. Let $h_1=z_3+z_4j \in S^3 \subset \HH$. 
$$ (z+(1+h)j)(z_3+z_4j)=z^{\prime}_1 +z^{\prime}_2 j$$
with $z^{\prime}_1=zz_3-(1+h)\overline{z}_4$, $z^{\prime}_2= zz_4 +(1+h)\overline{z}_3$, and
$$|z^{\prime}_1|^2 + |z^{\prime}_2|^2=|z|^2 +(1+h)^2=2+2h.$$
Then $$g_3(z_3+z_4j)(z;h)= \xi(\gamma_3(z_3+z_4j)(z:1+h))=\xi((z^{\prime}_1:z^{\prime}_2))=
(z^{\prime};h^{\prime}).$$

$$|z^{\prime}_2|^2=|z|^2|z_4|^2 +(1+h)^2|z_3|^2 +(1+h)(z z_3 z_4 + \overline{z z_3 z_4}).$$
$$\frac{|z^{\prime}_2|^2}{1+h}=1 +h(|z_3|^2-|z_4|^2) +(z z_3 z_4 + \overline{z z_3 z_4}).$$
$$h^{\prime}=\frac{2|z^{\prime}_2|^2-(|z^{\prime}_1|^2 + |z^{\prime}_2|^2)}{|z^{\prime}_1|^2 + |z^{\prime}_2|^2}
=\frac{|z^{\prime}_2|^2}{1+h}-1= h(|z_3|^2-|z_4|^2) +(z z_3 z_4 + \overline{z z_3 z_4}).$$
$$z^{\prime}_1\overline{z}^{\prime}_2=|z|^2 z_3 \overline{z}_4 -(1+h)^2 z_3 \overline{z}_4 + (1+h)z_3^2 z -(1+h) \overline{z}_4^2 \overline{z} $$
$$z^{\prime}=\frac{2 z^{\prime}_1\overline{z}^{\prime}_2}{|z^{\prime}_1|^2 + |z^{\prime}_2|^2}=\frac{z^{\prime}_1\overline{z}^{\prime}_2}{1+h}
=- 2h  z_3 \overline{z}_4 + z_3^2 z -  \overline{z}_4^2\overline{z}.$$
In particular, the map $g_3(z_3+z_4j)$ from $S^2$ to $S^2$ extends as an element of $GL(\RR^3)$ still denoted by $g_3(z_3+z_4j)$ with the matrix

$$g_3(z_3+z_4j)=\left[\begin{array}{ccc}
\mbox{Re}(z_3^2-z_4^2)       & \mbox{Im}(z_4^2-z_3^2) &
                                           -2\mbox{Re}({z}_3\overline{z}_4) \\
\;\;\mbox{Im}(z_3^2+z_4^2)  & \mbox{Re}(z_3^2+z_4^2)          &  
                                              -2\mbox{Im}({z}_3\overline{z}_4)\\
2\mbox{Re}({z}_3{z}_4) & -2\mbox{Im}({z}_3{z}_4)  &   |z_3|^2-|z_4|^2
\end{array}\right].$$
Let us now compute the matrix of the conjugacy 
$$\rho(z_3 + z_4 j) : v \mapsto (z_3 + z_4 j)v(\overline{z}_3- z_4 j).$$
$$(z_3 + z_4 j)i(\overline{z}_3- z_4 j)=i(|z_3|^2 - |z_4|^2) - 2 z_3z_4 k$$
$$(z_3 + z_4 j)j(\overline{z}_3- z_4 j)= (z_3^2+z_4^2) j + {z}_3\overline{z}_4 - z_4 \overline{z}_3$$
$$(z_3 + z_4 j)k(\overline{z}_3- z_4 j)=i (z_4 \overline{z}_3 + z_3\overline{z}_4) -z_4^2 k +z_3^2 k $$

$$\rho(z_3 + z_4 j)=\left[\begin{array}{ccc}
|z_3|^2 - |z_4|^2       & 2\mbox{Im}({z}_3\overline{z}_4) &
                                           2\mbox{Re}({z}_3\overline{z}_4) \\
\;\;2\mbox{Im}({z}_3{z}_4)  & \mbox{Re}(z_3^2+z_4^2)          &  
                                                      \mbox{Im}(z_4^2-z_3^2)\\
-2\mbox{Re}({z}_3{z}_4) & \mbox{Im}(z_3^2+z_4^2)          & 
                                                      \mbox{Re}(z_3^2-z_4^2)
\end{array}\right].$$
Therefore, $g_3(z_3 + z_4 j)=P_{13}\rho(z_3 + z_4 j)^{-1}P_{13}$, and we are done.
\eop

It is now easy to conclude the proof of Proposition~\ref{propcptrois}. Since $SO(3)$ is connected, the gluing map of Lemma~\ref{lemcompglu} is homotopic to $(v \mapsto \rho^{-1}(v))$. Now,
to conclude define the orientation-reversing diffeomorphism
$S$ from
$$ \CC P^3\cong B^4 \times S^2 \cup_{\partial B^4 \times S^2 \stackrel{\psi(\rho^{-1})}{\rightarrow} \partial (-B^4) \times S^2} (-B^4 \times S^2)$$
to $$S_2(E_1)=B^4 \times S^2 \cup_{\partial B^4 \times S^2 \stackrel{\psi(\rho)}{\rightarrow} \partial (-B^4) \times S^2} (-B^4 \times S^2)$$
by
$$S((x,v) \in B^4 \times S^2 \subset \CC P^3)=(x,v) \in -B^4 \times S^2 \subset S_2(E_1)$$
and 
$$S((x,v) \in -B^4 \times S^2 \subset \CC P^3)=(x,v) \in B^4 \times S^2 \subset S_2(E_1).$$

\eop

\begin{proposition}
\label{propxiun} \index{N}{ksione@$\xi_1$}
$$\xi_1= -\frac{1}{12}[\theta].$$
\end{proposition}
\bp The only degree one Jacobi diagram is 
$$\theta=\begin{pspicture}[.4](-1,-.5)(2,1.5)
\pscircle(0.5,0.5){.5}
\psline{*-*}(0,.5)(1,.5)
\rput[r](-.1,.5){$1$}
\rput[l](1.1,.5){$2$}
\rput[b](.25,.6){\small b}
\rput[b](.75,.6){\small B}
\rput[r](.05,.8){\small c}
\rput[r](0.05,.2){\small a}
\rput[l](.95,.8){\small C}
\rput[l](.95,.2){\small A}
\end{pspicture}.$$ 
Orient its edges from $1$ to $2$, and orient $V(\theta)=\{1,2\}$ with its natural order.
Then the edge-orientation of $\theta$ is given by the order $(a,A,b,B,c,C)$ that is equivalent to the order
$(a,b,c,B,A,C)$ of the vertex-orientation where the vertices of $\theta$ are oriented by the picture.
Therefore,
$$\xi_1=\frac{1}{12}\int_{S_2(E_1)}\omega_T^3[\theta].$$
Recall that $H^2(S_2(E_1))\cong H^2(\CC P^3)= \ZZ[\omega_{\CC P^3}]$
where $\omega_{\CC P^3}$ is Poincar\'e dual to $\CC P^2$ and $\int_{\CC P^1}\omega_{\CC P^3}=1$.
Since the orientation-reversing diffeomorphism $S$ from $\CC P^3$ to $S_2(E_1)$ restricts to an orientation-preserving diffeomorphism from a fiber
$\CC P^1$ of $\CC P^3$ to a fiber $S^2$ of $S_2(E_1)$,
$$\int_{S(\CC P^1)}(S^{-1})^{\ast}(\omega_{\CC P^3})=1=\int_{S(\CC P^1)}\omega_T.$$
Since $H^2(S_2(E_1))\cong \ZZ$, this shows that  $\omega_T=(S^{-1})^{\ast}(\omega_{\CC P^3})$. \index{N}{omegaT@$\omega_T$} Then
$$12\xi_1=\int_{S_2(E_1)}(S^{-1})^{\ast}(\omega_{\CC P^3})^3[\theta]=-\int_{\CC P^3}\omega_{\CC P^3}^3[\theta]=-[\theta].$$
\eop

\newpage
\section{Compactifications of configuration spaces}
\label{seccomp}
\setcounter{equation}{0}

In this section, we give a detailed description of the compactifications of the configuration spaces mentioned in Subsection~\ref{substaconf} and we prove all the statements of this subsection that is the introduction to this section.
These compactifications are similar to the Fulton and MacPherson compactifications \cite{fmcp} first used by Bott and Taubes in \cite{bt}. Here, we use the Poirier approach \cite{Po} to present them.
The used definitions and the used properties of blow-ups will be given in Subsection~\ref{subsdifblowup}.

\subsection{Topological definition of the compactifications}

For any subset $A$ of
$V$, recall the restriction map \index{N}{pA@$p_A$}
$$p_A: \breve{C}_V(M) \longrightarrow \breve{C}_A(M).$$

Let $M^A(\infty^A)$ \index{N}{MAinftyA@$M^A(\infty^A)$}
be the manifold obtained from $M^A$ by blowing-up $\infty^A=(\infty, \infty, \dots, \infty)$. When $\sharp A=1$, set $C(A;M)=M^A(\infty)$. When $\sharp A > 1$, define $C(A;M)$ \index{N}{CAM@$C(A;M)$} 
from $M^A(\infty^A)$ by blowing-up the closure of
the {\em strict diagonal\/} of $(M \setminus \infty)^A$ made of the  constant maps
from $A$ to $(M \setminus \infty)$. Proposition~\ref{propblodifdeux} asserts that $C(A;M)$ inherits 
a canonical differentiable structure from the differentiable structure of $M^A$.
Let $\Pi_A: C(A;M) \longrightarrow M^A$ \index{N}{PiA@$\Pi_A$} be the canonical projection.

Consider the embedding
$$\iota =\prod_{A \subseteq V, A \neq \emptyset}p_A: \breve{C}_V(M) \longrightarrow
\prod_{A \subseteq V, A \neq \emptyset}C(A;M)$$
and identify $\breve{C}_V(M)$ with its image under $\iota$.
Define $C_V(M)$ \index{N}{CV@$C_V(M)$}
as a topological space as the closure of $\iota(\breve{C}_V(M))$
in the compact space $\prod_{A \subseteq V, A \neq \emptyset}C(A;M)$.
Note that when $\sharp V=1$, $C_V(M)$ is 
homeomorphic to $C_1(M)$. \index{N}{Cone@$C_1(M)$}
We have the following lemma.

\begin{lemma}
\label{lemcond1}
Any $c=(c_A)_{A \subseteq V, A \neq \emptyset} \in C_V(M)$, satisfies the following property $(C1)$:
The restriction of $\Pi_V(c_V)$ to $A$ is equal to $\Pi_A(c_A)$.  
\end{lemma}
\bp Indeed, the set made of the configurations that satisfy $(C1)$ for a given $A$ is closed since it is the preimage
of the diagonal of $(M^A)^2$ under a continuous map. Furthermore, this set contains $\breve{C}_V(M)$. Therefore, it contains $C_V(M)$.\eop

Since we shall use the differentiable structure of the $C(A;M)$ to define the structure of $C_V(M)$. We first study the former one in detail.

\subsection{Differentiable structure on a blow-up}
\label{subsdifblowup}

\begin{definition}
A {\em \indexT{dilation}} is a homothety with ratio in $]0,\infty[$.
\end{definition}

In general, when $V$ is a vector space $SV=S(V)=\frac{V \setminus \{0\}}{]0,\infty[}$
denotes the quotient of $(V \setminus \{0\})$ by the action of $]0,\infty[$ that always operates by  scalar multiplication.
Recall that the {\em unit normal bundle\/} $SN_X(Z)$ of a submanifold $Z$ 
in a smooth manifold $X$ is a bundle over $Z$ whose fiber over $(z \in Z)$ is $S(\frac{T_zX}{T_zZ})$.

\begin{definition}
\label{defblodif}
 
As a set, the {\em \indexT{blow-up}\/} of $X$ along $Z$ is \index{N}{XZ@$X(Z)$}
$$X(Z)=(X \setminus Z) \cup SN_X(Z).$$ 
It is equipped with a canonical projection from $X(Z)$ to $X$ that is the identity
outside $SN_X(Z)$ and that is the bundle projection from $SN_X(Z)$ to $Z$ on $SN_X(Z)$. The following proposition defines the canonical smooth structure
of a blow-up.
\end{definition}

\begin{proposition}
\label{propblodifun}
Let $Z$ be a $C^{\infty}$ submanifold of a $C^{\infty}$ manifold $X$ that is transverse to the possible boundary $\partial X$ of $X$. 
The \indexT{blow-up} $X(Z)$ \index{N}{XZ@$X(Z)$} 
has a unique smooth structure of a manifold with corners such that 
\begin{enumerate}
\item the canonical projection from $X(Z)$ to $X$ is smooth and restricts to a diffeomorphism from $X \setminus Z$ to its
image in $X$,
\item any 
smooth diffeomorphism $\phi: [0,\infty[^c \times \RR^n \longrightarrow X$ from $[0,\infty[^c \times \RR^n$ to an open subset $\phi([0,\infty[^c \times \RR^n)$ in $X$ whose image intersects $Z$ exactly along
$\phi([0,\infty[^c \times \RR^{d-c} \times 0)$, for natural integers $c,d,k$ with $c \leq d$, provides a smooth embedding  
$$\begin{array}{lll}([0,\infty[^c \times \RR^{d-c}) \times [0, \infty[ \times S^{n+c-d-1} &\hfl{\tilde{\phi}} &X(Z) \\ (x,\lambda \in ]0, \infty[,v) &\mapsto & \phi(x,\lambda v) \\ (x,0,v) &\mapsto & D\phi(x,0)(v) \in SN_X(Z) \end{array}$$ 
with open image in $X(Z)$.
\end{enumerate}
\end{proposition}

\bp
We use local diffeomorphisms of the form $\tilde{\phi}$ and charts
on $X \setminus Z$ to build an atlas for $X(Z)$.
These charts are obviously compatible over $X \setminus Z$, and we need to 
check compatibility for charts 
 $\tilde{\phi}$ and $\tilde{\psi}$ induced by embeddings $\phi$ and $\psi$ as in the statement.
For these, transition maps 
read:
$$(x,\lambda,u) \mapsto (\tilde{x},\tilde{\lambda}, \tilde{u})$$
where
$$\tilde{x}=p_1 \circ \psi^{-1} \circ \phi(x,\lambda u)$$ 
$$ \tilde{\lambda} = \parallel p_2 \circ \psi^{-1} \circ \phi(x,\lambda u )\parallel$$
$$ \tilde{u} = \left\{ \begin{array}{ll}\frac{p_2 \circ \psi^{-1} \circ \phi(x,\lambda u )}{\tilde{\lambda}} &\mbox{if} \; \lambda \neq 0 \\
\frac{D\left(p_2 \circ \psi^{-1} \circ \phi(x,0 )\right)(u)dt}{\parallel D\left(p_2 \circ \psi^{-1} \circ \phi(x,0 )\right)(u)dt\parallel} &\mbox{if} \; \lambda = 0\end{array} \right. $$
In order to check that this is smooth, write
$$p_2 \circ \psi^{-1} \circ \phi(x,\lambda u ) = \lambda \int_0^1D\left(p_2 \circ \psi^{-1} \circ \phi(x,t \lambda u )\right)(u)dt $$
where the integral does not vanish when $\lambda$ is small enough.

More precisely, assuming $c=0$ for simplicity in the notation, since the restriction to $S^{n-d-1}$ of $D\left(p_2 \circ \psi^{-1} \circ \phi(x,0 )\right)$ is an injection, for any $u_0 \in S^{n-d-1}$, there exists a neighborhood of $(0,u_0)$ in $[0,\infty[ \times S^{n-d-1}$ such that for any $(\lambda,u)$ in this neighborhood, we have the following condition about the scalar product $$\langle D\left(p_2 \circ \psi^{-1} \circ \phi(x,\lambda u )\right)(u),D\left(p_2 \circ \psi^{-1} \circ \phi(x,0)\right)(u) \rangle \; > 0.$$
Therefore, there exists $\varepsilon  > 0$ such that for any $\lambda \in [0,\varepsilon[$, and for any $u \in S^{n-d-1}$,
$$\langle D\left(p_2 \circ \psi^{-1} \circ \phi(x,\lambda u )\right)(u),D\left(p_2 \circ \psi^{-1} \circ \phi(x,0)\right)(u) \rangle  \;  > 0.$$
Then 
$$ \tilde{\lambda} = \lambda \norm{\int_0^1D\left(p_2 \circ \psi^{-1} \circ \phi(x,t \lambda u )\right)(u)dt}$$
is a smooth function (defined even when $\lambda \leq 0$)
and 
$$ \tilde{u} =
\frac{\int_0^1D\left(p_2 \circ \psi^{-1} \circ \phi(x,t \lambda u )\right)(u)dt}{\parallel\int_0^1D\left(p_2 \circ \psi^{-1} \circ \phi(x,t \lambda u )\right)(u)dt\parallel}$$ is smooth, too.
\eop

\begin{proposition}
\label{propblodifdeux}
Let $Y$ be a $C^{\infty}$ submanifold of a $C^{\infty}$ manifold $X$ without boundary, and let $Z$ be a $C^{\infty}$ submanifold of $Y$. 
\begin{enumerate}
\item The boundary $\partial X(Z)$ of $X(Z)$ is canonically diffeomorphic to $SN_X(Z)$. 
\item The closure $\overline{Y}$ of $(Y \setminus Z)$ in $X(Z)$ is a submanifold of $X(Z)$ that intersects $\partial X(Z)$ as the unit normal bundle $SN_Y(Z)$ of $Z$ in $Y$. 
\item The blow-up $X(Z)(\overline{Y})$ of $X(Z)$ along $\overline{Y}$ has a canonical differential structure of a manifold with corners, and the preimage of \/ $\overline{Y} \subset X(Z)$ in $X(Z)(\overline{Y})$ under the canonical projection
$$X(Z)(\overline{Y}) \longrightarrow X(Z) $$
is a fibered space over $\overline{Y}$ with fiber the spherical normal bundle
of $Y$ in $X$ pulled back by $(\overline{Y} \longrightarrow Y)$.
\end{enumerate}
\end{proposition}
\bp

\begin{enumerate}
\item The first assertion is easy to observe from the charts in Proposition~\ref{propblodifun}.
\item Now, it is always possible to choose a chart $\phi$ as above such that 
furthermore the image of $\phi$ intersects $Y$ exactly along
$\phi(\RR^k \times 0)$, $k > d$.
Then, let us look at the induced chart $\tilde{\phi}$ of $X(Z)$ near a point of $\partial X(Z)$. \\
The intersection of $(Y \setminus Z)$ with the image of $\tilde{\phi}$ is $\tilde{\phi}\left(\RR^d \times ]0,\infty[  \times (S^{k-d-1}\subset S^{n-d-1})\right)$. Thus, the closure of $(Y \setminus Z)$ intersects the image of $\tilde{\phi}$ as $$\tilde{\phi}\left(\RR^d \times [0,\infty[  \times (S^{k-d-1}\subset S^{n-d-1})\right).$$
\item
Together with the above mentioned charts of $\overline{Y}$, the smooth  injective map $$S^{k-d-1} \times \RR^{n-k} \longrightarrow S^{n-d-1}$$ $$(u,y) \mapsto \frac{(u,y)}{\parallel (u,y) \parallel}$$ 
identifies $\RR^{n-k}$
with the fibers of the normal bundle of $\overline{Y}$ in $X(Z)$. The blow-up process
will therefore replace $\overline{Y}$ by the quotient of the corresponding $(\RR^{n-k} \setminus \{0\})$-bundle
by $]0,\infty[$ which is of course the pull-back under the natural projection $(\overline{Y} \longrightarrow Y)$ of the spherical normal bundle of $Y$ in $X$.
\end{enumerate}
\eop

\noindent{\sc Proof of Lemma~\ref{lemextproj}:}
According to Proposition~\ref{propblodifun}, near the diagonal of $\RR^3$, we have a chart of $C_2(S^3)$
$$\psi: \RR^3 \times [0,\infty[  \times S^2 \longrightarrow C_2(S^3)$$
that maps $( x \in \RR^3, \lambda \in ]0,\infty[ ,y \in S^2)$ to $(x, x + \lambda y) \in (\RR^3)^2$. Here, $p_{S^3}$ extends as the projection onto the $S^2$ factor.\\
Consider the embedding \index{N}{phiinfty@$\phi_{\infty}$}
$$\begin{array}{llll}\phi_{\infty}: &\RR^3 &\longrightarrow &S^3\\
& \mu (x \in S^2) & \mapsto & \left\{\begin{array}{ll} \infty \;&\;\mbox{if}\; \mu=0\\
\frac{1}{\mu}x \;&\;\mbox{otherwise.} \end{array}\right.\end{array}$$
This chart identifies $S(T_{\infty}S^3)$ to $S(\RR^3)$. When $\mu \neq 0$,
$$p_{S^3}(\phi_{\infty}(\mu x), y \in \RR^3)= \frac{\mu y-x}{\norm{\mu y- x}}.$$
Then $p_{S^3}$ can be smoothly extended on $S(T_{\infty}S^3) \times \RR^3$ by
$$p_{S^3}(D\phi_{\infty}(x) \in S(T_{\infty}S^3), y \in \RR^3) = -x.$$
Similarly, set
$$p_{S^3}( x \in \RR^3, D\phi_{\infty}(y \in S(\RR^3)) \in S(T_{\infty}S^3)) = y.$$
Now, when $$(x,y) \in \left(S((\RR^3)^2) \setminus S(\mbox{diag}((\RR^3)^2)) \stackrel{(D\phi_{\infty})^2}{\cong} S((T_{\infty}S^3)^2) \setminus S(\mbox{diag}((T_{\infty}S^3)^2))\right),$$
and when $x$ and $y$ are not equal to zero,
$$p_{S^3}(\phi_{\infty}(\lambda x),\phi_{\infty}(\lambda y))=\frac{\frac{y}{\norm{y}^2}-\frac{x}{\norm{x}^2}}
{\norm{\frac{y}{\norm{y}^2}-\frac{x}{\norm{x}^2}}}
=\frac{\norm{x}^2y-\norm{y}^2x}
{\norm{\norm{x}^2y-\norm{y}^2x}}.$$ 
Therefore, $p_{S^3}$ smoothly extends on $M^2(\infty,\infty)$ outside the boundaries of $\infty \times C_1(M)$, 
$C_1(M) \times \infty$ and $\mbox{diag}(C_1(M))$ as
$$p_{S^3}( (D\phi_{\infty})^2((x,y) \in S^5)) =\frac{\norm{x}^2y-\norm{y}^2x}
{\norm{\norm{x}^2y-\norm{y}^2x}}.$$
Let us check that $p_{S^3}$ smoothly extends over the boundary of the diagonal of $C_1(M)$. There is a chart of $C_2(M)$ near the preimage of this boundary
in  $C_2(M)$
$$\psi_2: [0,\infty[ \times [0,\infty[ \times S^2 \times S^2 \longrightarrow C_2(S^3)$$
that maps $(\lambda \in ]0,\infty[ , \mu \in  ]0,\infty[,  x \in S^2, y \in S^2)$
to $(\phi_{\infty}(\lambda x), \phi_{\infty}(\lambda(x + \mu y)))$ where $p_{S^3}$ reads
$$(\lambda,\mu,x,y) \mapsto \frac{y- 2\langle x,y \rangle x -\mu x}
{\norm{y- 2\langle x,y \rangle x -\mu x}},$$ and therefore smoothly extends when $\mu=0$. We similarly check that $p_{S^3}$ smoothly extends over the boundaries of $(\infty \times C_1(M))$ and 
$(C_1(M) \times \infty)$.
\eop

\subsection{The differentiable structure of C(A;M)} \index{N}{CAM@$C(A;M)$}

Recall that $M^A(\infty^A)$ is the manifold obtained from $M^A$ by blowing-up $\infty^A=(\infty, \infty, \dots, \infty)$. As a set, $M^A(\infty^A)$ is the union of $M^A \setminus \infty^A$ with the spherical tangent bundle $S\left((T_{\infty}M)^A \right)$ of $M^A$ at $\infty^A$.
Let $\overline{\mbox{diag}((M\setminus \infty)^A)}$ denote the closure in $M^A(\infty^A)$ of
the {\em strict diagonal \/} of $(M \setminus \infty)^A$ made of the  constant maps. The boundary of $\overline{\mbox{diag}((M\setminus \infty)^A)}$
is the strict diagonal of $(T_{\infty}M \setminus 0 )^A $ up to dilation.
This allows us to see all the elements of $\overline{\mbox{diag}((M\setminus \infty)^A)}$ as constant maps from $A$ to $C_1(M)$, and provides a canonical diffeomorphism $p_1:\overline{\mbox{diag}((M\setminus \infty)^A)} \longrightarrow C_1(M)$.

Now, $C(A;M)$ is obtained from $M^A(\infty^A)$ by blowing-up $\overline{\mbox{diag}((M\setminus \infty)^A)}$. 
Thus, as a set, $C(A;M)$ is the union of 
\begin{itemize}
\item the set of non constant maps from $A$ to $M$,
\item the space $\frac{(T_{\infty}M)^A \setminus \mbox{diag}((T_{\infty}M)^A)}{]0,\infty[}$, 
and,
\item the bundle over $\overline{\mbox{diag}((M\setminus \infty)^A)}=C_1(M)$ whose fiber at a constant map with value $x \in C_1(M)$ is 
$$S\left(\frac{T_{\Pi_1(x)}M^A}{\mbox{diag}((T_{\Pi_1(x)}M)^A)} \right)$$ 
\end{itemize}
Note that $(\frac{T_{\Pi_1(x)}M^A}{\mbox{diag}((T_{\Pi_1(x)}M)^A)} \setminus 0)$ may be identified with
$((T_{\Pi_1(x)}M)^{A \setminus b} \setminus 0)$ for any $b \in A$ through $$[(v_a)_{a \in A}] \mapsto
(v_a - v_b)_{a \in (A \setminus b)}.$$
Recall that for $A \subset V$, $\Pi_A: C(A;M) \longrightarrow M^A$
 denotes the canonical projection.
When $V$ is a euclidean vector space, $S(V)$ is simply the unit sphere of $V$.

\begin{example}
\label{exa1}{\bf Charts near $\Pi_A^{-1}(\mbox{diag}((M \setminus \infty)^A))$.}\\
Let $$\phi: \RR^3 \longrightarrow M \setminus \infty$$ be a smooth embedding
that is a chart of $M$ near $\phi(0)=x$. 
Let $A$ be a finite set of cardinality $\sharp A \geq 2$. Let $b \in A$.
Let us construct an explicit chart $\psi(A;\phi;b)$ of $C(A;M)$ near a point of 
$\Pi_A^{-1}(x^A)$ where $x^A$ denotes the constant map of $M^A$ with value $x$.

We have the chart
$$\begin{array}{llll}\tilde{\psi}(A;\phi;b):& \RR^3 \times (\RR^3)^{(A \setminus b)} &\longrightarrow &M^A\\
&(y,(y_c)_{c \in A \setminus b}) & \mapsto & \left(c \mapsto \left\{ \begin{array}{ll}\phi(y)& \mbox{if}\; c=b\\\phi(y+y_c)& \mbox{if}\; (c \in A \setminus b)\end{array}\right.\right)
\end{array}$$
of submanifold for the strict diagonal, 
and this induces the chart \index{N}{psiAphiB@$\psi(A;\phi;b)$}
$$\begin{array}{lll} \RR^3 \times [0,\infty[ \times S((\RR^3)^{(A \setminus b)}) &\hfl{\psi(A;\phi;b)} &C(A;M)\\
(y, \lambda \in ]0,\infty[, (y_c)_{c \in (A \setminus b)}) & \mapsto & 
\left(c \mapsto \left\{ \begin{array}{ll}\phi(y)& \mbox{if}\; c=b
\\\phi(y+ \lambda y_c)& \mbox{if}\; (c \in A \setminus b)\end{array}\right.\right)\\
(y, 0, (y_c)_{c \in (A \setminus b)}) & \mapsto & \left(D\phi(y)(y_c)\right)_{c \in (A \setminus b)} \in S\left(\frac{(T_{\phi(y)}M)^A}{\mbox{\small diag}((T_{\phi(y)}M)^A)}\right)
\end{array}$$
for $C(A;M)$ in $\Pi_A^{-1}(\phi(\RR^3)^A)$. \end{example}
Let $S\left(\frac{T(M\setminus \infty)^A}{\mbox{\small diag}(T(M\setminus \infty)^A)}\right)$ denote the total space of the fibration over 
$(M\setminus \infty)$ whose fiber over $x \in (M\setminus \infty)$ is $S\left(\frac{T_xM^A}{\mbox{\small diag}(T_xM^A)}\right)$.
Let
$$\Pi_d: \Pi_A^{-1}(\mbox{diag}(M \setminus \infty)^A) \longrightarrow S\left(\frac{T(M\setminus \infty)^A}{\mbox{\small diag}(T(M\setminus \infty)^A)} \right),$$
denote the canonical projection.
An element in the target of $\Pi_d(\Pi_A^{-1}(x^A))$ will be seen as a non-constant map from $A$ to $T_xM$ up to translation and up to dilation.

\begin{lemma}
\label{lemcond2}
Any $c=(c_A)_{A \subseteq V, A \neq \emptyset} \in C_V(M)$, satisfies the following property $(C2)$:
For any two subsets $A$ and $B$ of $V$ such that the cardinality of $A$ is greater than $1$ and $A \subset B$, 
if $c_B \in \Pi_B^{-1}(\mbox{diag}(M \setminus \infty)^B)$, then the restriction
to $A$ of 
$\Pi_d(c_B)$ is a (possibly null) positive multiple of $\Pi_d(c_A)$. 
\end{lemma}
\bp 
Choose a basepoint $b \in A$ for $A$ and $B$. 
Consider the projection $\Pi_{AB}$ of $\prod_{C \subseteq V, C \neq \emptyset}C(C;M) $ onto $C(A;M) \times C(B;M)$
in a neighborhood of some $c$ such that $x^B=\Pi_B(c_B)$ and $x^A=\Pi_A(c_A)$,
with $x \in M\setminus \infty$.
Then $$(\psi(A;\phi;b)^{-1} \times \psi(B;\phi;b)^{-1}) \circ \Pi_{AB}$$ map the elements of $\breve{C}_V(M)$ to elements of the form
$(y, \lambda_A, u_A, y, \lambda_B, u_B)$ where $y \in \RR^3$, $\lambda_A,\lambda_B \in ]0,+\infty[$, $u_A \in (\RR^3)^{A \setminus b}$, 
$u_B \in (\RR^3)^{B \setminus b}$, $\parallel u_A \parallel=\parallel u_B \parallel=1$ and,
$$\lambda_B p(u_B) =\lambda_A u_A$$
where $p$ is the natural projection (or restriction) from $(\RR^3)^{B \setminus b}$ to $(\RR^3)^{A \setminus b}$.
In particular, $p(u_B)$ and $u_A$ are colinear in $(\RR^3)^{A \setminus b}$,
and their scalar product is $\geq 0$. These two conditions define a closed 
subset of $\left((\RR^3)^{A \setminus b}\right)^2$. Therefore, they must be satisfied in the image of the closure $C_V(M)$. Since they read as stated when $c_B \in \Pi_B^{-1}(\mbox{diag}(M \setminus \infty)^B)$ , that is when $\lambda_B=0$, (and hence $\lambda_A=0$, too) we are done. \eop

Let
$$\Pi_{\infty}: \Pi_A^{-1}(\infty^A) \longrightarrow S\left((T_{\infty}M)^A \right) \subset M^A(\infty^A)$$ denote the canonical projection.

An element in the target of $\Pi_{\infty}$ will be seen as a non-zero map from $A$ to $T_{\infty}M$  up to dilation. 

\begin{example}
\label{exa2}{\bf Charts of $M^A(\infty^A)$ near $\Pi_A^{-1}(\infty^A)$.}\\
Let $$\phi_{\infty}: \RR^3 \longrightarrow M $$ be a smooth embedding
such that $\phi_{\infty}(0)=\infty$. Then the composition
$$]0,\infty[ \times S((\RR^3)^A) \hfl{\mbox{multiplication}} (\RR^3)^A
\hfl{(\phi_{\infty})^A} M^A$$
induces the chart \index{N}{psiAphiinfty@$\psi(A;\phi_{\infty})$}
$$\begin{array}{llll}\psi(A;\phi_{\infty}):& [0,\infty[ \times S((\RR^3)^A) &\longrightarrow & M^A(\infty^A)\\
&(\lambda,u) & \mapsto &\phi_{\infty} \circ \lambda u \;\;\;\;\;\;\; \mbox{when} \; \lambda \neq 0.
\end{array}$$
Here, $u$ is seen as a map from $A$ to $\RR^3$.\\
Note that $\Pi_{\infty}(\psi(A;\phi_{\infty})(0,u))=D_0 \phi_{\infty} \circ u$.
\end{example}

\begin{lemma}
\label{lemcond3}
Any $c=(c_A)_{A \subseteq V, A \neq \emptyset} \in C_V(M)$, satisfies the following property $(C3)$:
For any two non-empty subsets $A$ and $B$ of $V$ such that $A \subset B$, 
if $c_B \in \Pi_B^{-1}(\infty ^B)$, then the restriction
to $A$ of 
$\Pi_{\infty}(c_B)$ is a (possibly null) positive multiple of $\Pi_{\infty}(c_A)$. 
\end{lemma}
\bp This can be proved along the same lines as Lemma~\ref{lemcond2} using the chart of Example~\ref{exa2}, and this is left to the reader. 
\eop

\begin{example}
\label{exa3}{ \bf Charts of $C(A;M)$ near the intersection of $\Pi_A^{-1}(\infty^A)$ and the closure of the strict diagonal of $(M \setminus \infty)^A$.}\\
Use the notation of the previous example~\ref{exa2}.
Let $b \in A$. Assume $\sharp A  > 1$.
From $\tilde{\psi}(A;\phi_{\infty};b)$
$$\begin{array}{lll} ]0,\infty[ \times S\left(\RR^3 \times (\RR^3)^{(A \setminus b)}\right) &\longrightarrow &M^A\\
(\lambda;y,(y_c)_{c \in (A \setminus b)}) & \mapsto & \left(c \mapsto \left\{ \begin{array}{ll}\phi_{\infty}( \frac{\lambda}{\sqrt{\sharp A}}y)& \mbox{if}\; c=b\\\phi_{\infty}(\lambda(\frac{1}{\sqrt{\sharp A}}y+y_c))& \mbox{if}\; c\neq b\end{array}\right.\right)
\end{array}$$
we get a chart \index{N}{psiAphiinftyb@$\psi(A;\phi_{\infty};b)$}
$${\psi}(A;\phi_{\infty};b):[0,\infty[ \times S^2 \times [0,\infty[
\times S\left((\RR^3)^{(A \setminus b)}\right) \longrightarrow C(A;M)$$
with the property that
$$\Pi_A({\psi}(A;\phi_{\infty};b)(\lambda,u,\mu,v))=
\phi_{\infty} \circ \lambda\left((\frac{1}{\sqrt{\sharp A}}u)^A + \mu v \right)$$
as a map from $A$ to $M$, where $v(b)=0$.
In particular $$\Pi_A^{-1}(\infty^A) \cap \mbox{Im}\left(\psi(A;\phi_{\infty};b)\right)=\psi(A;\phi_{\infty};b)\left(\{0\} \times S^2 \times [0,\infty[
\times S\left((\RR^3)^{(A \setminus b)}\right)\right),$$
and
$$\begin{array}{llll} \Pi_{\infty}:& \Pi_A^{-1}(\infty^A)  \cap \mbox{Im}\left(\psi(A;\phi_{\infty};b)\right) &\longrightarrow &S\left((T_{\infty}M)^A\right) \subset M^A(\infty^A)
\\ &
\psi(A;\phi_{\infty};b)(0,u,\mu,v) &\mapsto &  D_0\phi_{\infty} \circ (\frac{1}{\sqrt{\sharp A}}u^A + \mu v)\end{array}$$
where $u^A$ stands for the constant map with value $u$.\\
The boundary $\Pi_{\infty}^{-1}(\mbox{diag}(T_{\infty}M)^A)$ of 
$\overline{\mbox{diag}((M\setminus \infty)^A)}$
is $\psi(A;\phi_{\infty};b)(\{0\} \times S^2 \times \{0\}
\times S\left((\RR^3)^{(A \setminus b)}\right))$. The projection $p_1$ of $\psi(A;\phi_{\infty};b)(0,u,0,v)$ onto the boundary of $C_1(M)$ is $D_0\phi_{\infty}(u) \in S(T_{\infty}(M))$.
\end{example}

Let
$$\Pi_{\infty,d}:\Pi_{\infty}^{-1}(\mbox{diag}(T_{\infty}M)^A)) \longrightarrow S\left(\frac{T_{\infty}M^A}{\mbox{diag}(T_{\infty}M)^A}\right)$$
denote the canonical map.
Note that it reads
$$\psi(A;\phi_{\infty};b)(0,u,0,v) \mapsto D_0\phi_{\infty}\circ v$$
in the above charts.

\begin{lemma}
\label{lemcond4}
Any $c=(c_A)_{A \subseteq V, A \neq \emptyset} \in C_V(M)$, satisfies the following property $(C4)$:
For any two subsets $A$ and $B$ of $V$ such that the cardinality of $A$ is greater than $1$ and $A \subset B$, 
if $c_B \in \Pi_B^{-1}(\infty ^B)$, and if $\Pi_{\infty}(c_B)$ is a constant map (or is diagonal) then the restriction
to $A$ of 
$\Pi_{\infty,d}(c_B)$ is a (possibly null) positive multiple of $\Pi_{\infty,d}(c_A)$. 
\end{lemma}
\bp Again, this can be seen on the charts given in the previous example. 
Consider the projection $\Pi_{AB}$ of $\prod_{C \subseteq V, C \neq \emptyset}C(C;M) $ onto $C(A;M) \times C(B;M)$
in a neighborhood of some $c$ such that $\infty^B=\Pi_B(c_B)$, $\Pi_{\infty}(c_B)$ is constant, $\infty^A=\Pi_B(c_A)$ and $\Pi_{\infty}(c_A)$ is constant.
Then $$(\psi(A;\phi_{\infty};b)^{-1} \times \psi(B;\phi_{\infty};b)^{-1}) \circ \Pi_{AB}$$ 
maps the elements of $\breve{C}_V(M)$ to elements of the form
$$(\lambda_A, u_A, \mu_A, v_A, \lambda_B, u_B,\mu_B, v_B)$$ 
where  $\lambda_A,\lambda_B, \mu_A, \mu_B \in ]0,+\infty[$, $u_A,u_B \in S^2$, 
$v_A \in S((\RR^3)^{A \setminus b})$, $v_B \in S((\RR^3)^{B \setminus b})$,  and,
$$ \frac{\lambda_B}{\sqrt{\sharp B}}u_B =\frac{\lambda_A}{\sqrt{\sharp A}}u_A$$
$$ {\lambda_B\mu_B}p(v_B) ={\lambda_A\mu_A}v_A$$
where $p$ is the natural projection (or restriction) from $(\RR^3)^{B \setminus b}$ to $(\RR^3)^{A \setminus b}$.
Now, it is easy to conclude as before.\eop

\subsection{Sketch of construction of the differentiable structure of $C_V(M)$}
\label{subsketchdifcv} \index{N}{CV@$C_V(M)$}

In this subsection, we sketch the construction of the differentiable structure of $C_V(M)$ and we reduce the proofs of Propositions~\ref{propconfunc}, \ref{propconfuns}, \ref{propconffaceun}, \ref{propconffacedeux} to the proofs of Lemmas~\ref{lempropcons},~\ref{lempropconsadface} and Proposition~\ref{propcdeuxcoinc} stated
below.

We shall use the notation $A \subseteq B$ (resp. $A \subset B$) to say that $A$ is a subset (resp. strict subset) of $B$.
Define $\tilde{C}_V(M)$ \index{N}{CtildeV@$\tilde{C}_V(M)$} 
to be the set of the elements $c=(c_A)_{A \subseteq V, A \neq \emptyset}$ of $$\prod_{A \subseteq V, A \neq \emptyset}C(A;M)$$
that satisfy the properties $(C1)$, $(C2)$, $(C3)$ $(C4)$, of Lemmas~\ref{lemcond1},
\ref{lemcond2}, \ref{lemcond3} and \ref{lemcond4}. These lemmas ensure that $C_V(M)$ is a subset of $\tilde{C}_V(M)$.

An element of $\tilde{C}_V(M)$ is a map $(\Pi_V(c_V)\in M^V)$ from $V$ to $M$ with additional
data that allow us to see 
\begin{itemize}
\item the restricted configurations corresponding to a multiple point $x \neq \infty$ at any scale $A \subseteq \Pi_V(c_V)^{-1}(x)$ as a non-constant map $\Pi_d(c_A)$ from  $A$ to
$T_xM$ up to dilation and translation,
\item  the restricted configuration at a scale $A \subseteq \Pi_V(c_V)^{-1}(\infty)$ first as a non-zero map $\Pi_{\infty}(c_A)$ from $A$
to $T_{\infty}M$ up to dilation, and, if this latter map is constant,
\item with an additional zoom, the restricted configuration at a smaller scale as another independent non-constant map $\Pi_{\infty,d}(c_A)$ from $A$
to $T_{\infty}M$ up to dilation and translation.
\end{itemize}
with respective compatibity conditions $(C2)$, $(C3)$, $(C4)$.
Therefore, elements of $\tilde{C}_V(M)$ will be called {\em limit configurations.\/}

We are going to prove that $C_V(M)$ \index{N}{CV@$C_V(M)$}
is equal to $\tilde{C}_V(M)$ \index{N}{CtildeV@$\tilde{C}_V(M)$} 
and to construct a differentiable structure for $C_V(M)$ by proving the following proposition.

\begin{proposition}
\label{propdifconf}
For any $c^0 \in \tilde{C}_V(M)$, \index{N}{CtildeV@$\tilde{C}_V(M)$} 
there exist
\begin{enumerate}
\item $k \in \NN$, and an open neighborhood $O$ of $0$ in   $ [0,\infty[^k$,
(set $[0,\infty[^0=]0,\infty[^0=\{0\}$ if $k=0$)
\item an open neighborhood $W$ of a point $w^0$ in a smooth manifold $\tilde{W}$ without boundary,
\item an open neighborhood $U$ of $c^0$ in $\prod_{A \subseteq V, A \neq \emptyset}C(A;M)$,
\item  a smooth map $\xi: O \times W \longrightarrow U$ such that  $\xi(0;w^0)=c^0$,
$\xi(O \times W) \subset \tilde{C}_V(M)$, and $\xi((O \cap ]0,\infty[^k)\times W) = \breve{C}_V(M) \cap \xi(O \times W)$,
\item a smooth map $r: U \longrightarrow \RR^k \times \tilde{W}$ such that 
\begin{itemize}
\item $r \circ \xi$ is the identity of $O\times W$, 
\item $r(U \cap \tilde{C}_V(M)) \subseteq O \times W$, and 
\item the restriction of $\xi \circ r$ to $U \cap \tilde{C}_V(M)$ is the identity of $U \cap \tilde{C}_V(M)$.
\end{itemize}
(This implies that $\xi(O \times W)=U \cap \tilde{C}_V(M)$.)
\end{enumerate}
\end{proposition}

Proposition~\ref{propdifconf} easily implies that our $\xi$ form an atlas for
$\tilde{C}_V(M)$ that becomes a smooth manifold with corners and that 
$C_V(M)=\tilde{C}_V(M)$. Furthermore, with such an atlas, the inclusion $\iota$ from $C_V(M)$ to $\prod_{A \subseteq V, A \neq \emptyset}C(A;M)$ will be smooth, and a map $f$ from a smooth manifold $X$ to $C_V(M)$ will be smooth if and only if $\iota \circ f$ is smooth.

When $\sharp V=1$, we observe at once that $\tilde{C}_V(M)$ is diffeomorphic to $C(V;M)$. Therefore, our two definitions of $C_1(M)$ coincide.
We shall prove the following proposition in Subsection~\ref{subpropcdeuxcoinc}.

\begin{proposition}
\label{propcdeuxcoinc}
Let $V=\{1,2\}$. Let $C_2(M)$ \index{N}{Ctwo@$C_2(M)$} 
denote the manifold obtained from $M^2(\infty,\infty)$ by blowing up $\infty \times C_1(M)$, 
$C_1(M) \times \infty$ and $\mbox{diag}(C_1(M))$
as in Subsection~\ref{subconfap}. Let $C_V(M)$ be the compactification of 
$\breve{C}_2(M)$ defined in this subsection.
Then $C_2(M)$ is canonically diffeomorphic to $C_V(M)$.
\end{proposition}

\begin{lemma}
\label{lemdifconfimpconfunc}
Propositions~\ref{propdifconf} and \ref{propcdeuxcoinc} imply Proposition~\ref{propconfunc}.
\end{lemma}
\bp
It is obvious from Proposition~\ref{propdifconf} that $(\breve{C}_V(M) = \iota(\breve{C}_V(M)))$ is the interior of $C_V(M)$. On $\breve{C}_{V}(M)$, for $e=(a,b)$, $p_e$ is given by the projection on 
$C(e;M)$ that determines the projections on $C(\{a\};M)$ and $C(\{b\};M)$.
These projections naturally extend from $\tilde{C}_V(M)$ to the closure of the image of $\breve{C}_{V}(M)$ in $C_e(M)$, and they will be smooth because they come from the smooth projections and because of the forms of our charts.
\eop

Proposition~\ref{propconfuns} is easier to prove than Proposition~\ref{propconfunc} and could be proved before. 
Nevertheless, we shall focus on the proof of Proposition~\ref{propconfunc} and see
Proposition~\ref{propconfuns} as a particular case of Proposition~\ref{propconfunc} with the help of the following proposition
\ref{propdifconfad}.

Let $0^V$ denote the constant map with value $0$ in $(\RR^3)^V$, where $\RR^3=S^3 \setminus \infty$.
The preimage of $0^V$ under the canonical projection $\Pi_V:C(S^3;V) \longrightarrow (S^3)^V$ is the set of non-constant maps from $V$ to $T_0(\RR^3)$ up to dilation and translation. This allows us to see
$\breve{S}_V(\RR^3)$ as an open submanifold of $\Pi_V^{-1}(0^V)$. 
Furthermore, for a given element $s_V$ of $\breve{S}_V(\RR^3)$, there is a unique element of $\tilde{C}_V(S^3)$ whose projection on $C(V;M)$ is $s_V$ (by (C2) that determines its other projections). This allows us to see $\breve{S}_V(\RR^3)$ as a subset of $\tilde{C}_V(S^3)$. 
Set \index{N}{SVR@$S_V(\RR^3)$}
$$S_V(\RR^3)=(\Pi_V \circ p_V)^{-1}(0^V) \cap \tilde{C}_V(S^3).$$
$S_V(\RR^3)$ is a compact set that contains $\breve{S}_V(\RR^3)$.
Proposition~\ref{propconfuns} now becomes the consequence of the following proposition (together with Proposition~\ref{propdifconf}) by a proof similar to the proof of Lemma~\ref{lemdifconfimpconfunc}
above.

\begin{proposition}
\label{propdifconfad}
For any $c^0 \in \tilde{C}_V(S^3)$ such that $\Pi_V \circ p_V(c^0)=0^V$, in Proposition~\ref{propdifconf}, we have 
\begin{enumerate}
\item $k \geq 1$, $\tilde{W}=\RR^3 \times \tilde{\tilde{W}}$,
\item $S_V(\RR^3) \cap U=\xi\left(O \times W \cap \left(\{0\} \times [0,\infty[^{k-1} \times \{0\} \times \tilde{\tilde{W}}\right)\right)$,
\item $\xi\left(O \times W \cap \left(\{0\} \times ]0,\infty[^{k-1} \times \{0\} \times \tilde{\tilde{W}}\right)\right) = \breve{S}_V(\RR^3) \cap U$.
\end{enumerate}
\end{proposition}

Proposition~\ref{propdifconf} and Proposition~\ref{propdifconfad} are a consequence of the two following lemmas.

\begin{lemma}
\label{lempropcons}
Proposition~\ref{propdifconf} and Proposition~\ref{propdifconfad} are true when $\Pi_V(c^0_V)$ is a constant map of $M^V$.
\end{lemma}

\begin{lemma}
\label{lemimp}
(1) Lemma~\ref{lempropcons} implies Proposition~\ref{propdifconf}.\\
(2) Assume Lemma~\ref{lempropcons} is true.
Let $(A_i)_{i=1,2,\dots, s}$ be a partition of $V$ into nonempty subsets
$$V = \coprod_{i=1}^s A_i.$$
Let $\phi_i:\RR^3 \longrightarrow M$, for $i=1, \dots, s$, be embeddings with disjoint images in $M$.
Let $U_A$ be the following open subset of $C(A;M)$.
$$U_A=\{c_A \in C(A;M); \Pi_A(c_A)( A \cap A_i ) \subseteq \phi_i(\RR^3)\}$$
and define $$Q_V=\tilde{C}_V(M) \cap \prod_{\emptyset \subset A \subseteq V}U_A\;\; \mbox{and}
\;\; Q_{A_i}=\tilde{C}_{A_i}(M) \cap \prod_{\emptyset \subset A \subseteq {A_i}}U_A$$
Then the map $(Q_V \longrightarrow \prod_{i=1}^sQ_{A_i})$ induced by the restrictions is a diffeomorphism.
\end{lemma}
\noindent{\sc Proof of Lemma~\ref{lemimp}:}
(1) Let $c^0=(c^0_A)_{A \subseteq V, A \neq \emptyset} \in \tilde{C}_V(M)$.
Consider the map $\Pi_V(c_V)$ from $V$ to $M$ and set
$$\Pi_V(c_V)(V)=\{m_1, m_2, \dots, m_s\}$$
and $$A_i=\Pi_V(c_V)^{-1}(m_i)$$
Choose embeddings $\phi_i:\RR^3 \longrightarrow M$, for $i=1, \dots, s,$ with disjoint images in $M$ such that $\phi_i(0)=m_i$.
Let $c^i=c^0_{|A_i}=(c^0_A)_{A \subseteq A_i, A \neq \emptyset} \in \tilde{C}_{A_i}(M)$ denote the restriction of $c^0$ to $A_i$. 
According to Lemma~\ref{lempropcons}, we may find $k_i,U_i,O_i,W_i,w^0_i,\xi_i,r_i$ satisfying the conclusions
of Proposition~\ref{propdifconf} with $(c^i,A_i)$ instead of $(c^0,V)$, and after
a possible reduction of $U_i,O_i,W_i$, we may assume that $U_i \subseteq \prod_{A \subseteq A_i}U_A$.
Then, set $k=\sum_{i=1}^s k_i$, $O=\prod_{i=1}^s O_i$, $W=\prod_{i=1}^s W_i$,
$w^0=(w^0_i)_{i \in \{1, \dots, s\}}$.
$$U=\prod_{i=1}^s U_i \times \prod_{A \subseteq V; \forall i, A \cap (A \setminus A_i) \neq \emptyset}C(A;M)$$
Define $\xi((v,w)=(v_1, \dots, v_s,w_1, \dots, w_s))=(\xi(v,w)_A)_{A \subseteq V; A \neq \emptyset}$
by $\xi(v,w)_A=\xi_i(v_i,w_i)_A$ if $A \subseteq A_i$
and $\Pi_A(\xi(v,w)_A)(a\in A_i)=\Pi_{A_i}(\xi(v,w)_{A_i})(a)$.
When $A$ intersects all the $(A \setminus A_i)$, $\Pi_A(\xi(v,w)_A)$ is not constant. Since the restriction of $\Pi_A$ to the preimage of the set of non-constant maps is a diffeomorphism onto its image, $\xi(v,w)_A$ is smoothly well-determined for these $A$. Therefore $\xi$ is well-determined and smooth. Furthermore, $\xi(v,w)$ satisfies
$(C1)$ by construction and $\xi(v,w)$ satisfies the other conditions $(C2)$, $(C3)$ and $(C4)$ that are (thanks to $(C1)$ and to the choice of the $U_i$) conditions on some $\xi(v,w)_A$ and $\xi(v,w)_B$
for $A \subset B \subseteq A_i$. It is easy to see that $\xi(0,w^0)=c^0$,
and $\xi((O \cap ]0,\infty[^k)\times W) \subset \breve{C}_V(M)$ since the elements of $\breve{C}_V(M)$ are the elements $c$ of $\tilde{C}_V(M)$ such that $\Pi_V(c_V) \in M^V$ is an injective map from $V$ to $(M \setminus \infty)$.
We also easily see that $$\breve{C}_V(M) \cap \xi( O \times W) \subseteq \xi((O \cap ]0,\infty[^k)\times W).$$
When $r_i(u_i \in U_i)=(r^1_i(u_i) \in \RR^{k_i}; r_i^2(u_i) \in W_i)$, define $$r((u_i)_{i \in \{1, \dots, s\}};(c_A)_{A \subseteq V; \forall i, A \cap (A \setminus A_i) \neq \emptyset})$$ $$=((r^1_i(u_i))_{i \in \{1, \dots, s\}};(r^2_i(u_i))_{i \in \{1, \dots, s\}}).$$
Now, it is easy to see that Lemma~\ref{lempropcons} implies Proposition~\ref{propdifconf}.
The second part (2) of the lemma follows from the above proof.
\eop

Assume that Proposition~\ref{propdifconf} is true and come back to the faces defined in Subsection~\ref{substaconf}.
First recall that 
$F(\infty;V)=S_i(T_{\infty}M^V) \subseteq S((T_{\infty}M)^V)$
embeds in $C(V;M)$. This embedding is smooth and canonical.
Furthermore, by $(C1)$ and $(C3)$, there is a unique map of $F(\infty;V)$
into
$C_V(M)$ whose composition with the projection on $C(V;M)$ is the above embedding. Since the restrictions are smooth from $C(V;M) \cap S_i(T_{\infty}M^V)$ to the $C(A;M)$ for $A \subset V$, the charts of Proposition~\ref{propdifconf} for $C_V(M)=\tilde{C}_V(M)$ make clear that $F(\infty;V)$ smoothly injects into $C_V(M)$. Lemma~\ref{lemimp} allows us 
to conclude that for any non-empty subset $B$ of $V$, $F(\infty;B)$ injects
into $C_V(M)$, smoothly and canonically. It is easy to see that the projections
$p_e$ associated to pairs of elements of $V$ restrict to the image of $F(\infty;B)$
as described in Subsection~\ref{substaconf}. The reader can similarly check
that, for any subset $B$ of $V$ with $(\sharp B \geq 2)$, $F(B)$ smoothly and canonically injects into $C_V(M)$ and that the restrictions of the $p_e$ to the images of the $F(B)$ are described in  Subsection~\ref{substaconf}.
Obviously, the images of the $F \in \partial_1(C_V(M))$ are disjoint.

Let us inject $(f(B)=f(B)(\RR^3))$ into $S_V(\RR^3)$ where $B$ is a strict subset of $V$, $\sharp B \geq 2$.
Identify $\breve{S}_{\{b\} \cup(V \setminus B)}(\RR^3)$ with a subspace of $S\left(\frac{(\RR^3)^V}{\mbox{\small diag}((\RR^3)^V)}\right)$ made of maps that are constant on $B$, by setting $c(B)=c(b)$. In particular, $\breve{S}_{\{b\} \cup(V \setminus B)}(\RR^3)$ smoothly embeds into $C(V;S^3) \cap \Pi_V^{-1}(0^V)$. When $A$ is a non-empty subset of $V$ that is not a subset of $B$, $\breve{S}_{\{b\} \cup(V \setminus B)}(\RR^3)$ smoothly projects to
$C(A;S^3) \cap \Pi_A^{-1}(0^A)$ by the restrictions imposed by
$(C1)$ and $(C2)$ (that do not determine anything for the subsets of $B$ where $c$
is constant).
Now, $\breve{S}_{B}(\RR^3)$ smoothly embeds into $C(B;S^3) \cap \Pi_B^{-1}(0^B)$  and smoothly projects to
$C(A;S^3) \cap \Pi_A^{-1}(0^A)$, when $A$ is a non-empty subset of $B$, by the restrictions imposed by
$(C1)$ and $(C2)$. This allows us to define a canonical smooth injection of $f(B)(\RR^3)$ into
$S_V(\RR^3)$, and the $p_e$ have the desired form on the image. When $B$ and $B^{\prime}$ are two disjoint subsets of $V$, $f(B)$ and $f(B^{\prime})$
are disjoint.

The $F(\infty;B)$, $F(B)$ and $f(B)(\RR^3)$ will be identified with their
images.

\begin{lemma}
\label{lempropconsadface}
Assume that Proposition~\ref{propdifconf} is true. In Proposition~\ref{propdifconf},
\begin{itemize}
\item when $\Pi_V(c^0_V)$ is a constant map with value $m \in (M \setminus \infty)$,\\ $k=1$ if and only if $c^0 \in F(V)$,
\item when $\Pi_V(c^0_V)$ is the constant map with value $\infty$,\\ $k=1$ if and only if $c^0 \in F(\infty;V)$, 
\item when $\Pi_V(c^0_V)$ is the constant map $0^V$ of $(S^3)^V$,\\ $k=2$ if and only if  $c^0 \in f(B)(\RR^3)$ for some strict subset $B$ of $V$ with $\sharp B \geq 2$.
\end{itemize}
\end{lemma}

\noindent{\sc Proof of Proposition~\ref{propconffaceun} assuming Proposition~\ref{propdifconf} and Lemma~\ref{lempropconsadface}:}

Let $c$ belong to a codimension one face of $\tilde{C}_V(M)$.
As in the proof of Lemma~\ref{lemimp}, set
$\Pi_V(c_V)(V)=\{m_1, m_2, \dots, m_s\}$,
and $A_i=\Pi_V(c_V)^{-1}(m_i)$.
Choose embeddings $\phi_i:\RR^3 \longrightarrow M$, for $i=1, \dots, s$ with disjoint images in $M$ such that $\phi_i(0)=m_i$. Then by Lemma~\ref{lemimp},
if $c_{|A_i}$ belongs to a codimension $d(i)$ face, then $c$ belongs to 
a codimension $\left(\sum_{i=1}^sd(i) \right)$-face. Therefore, there exists
a unique $j$ such that $c_{|A_j}$ belongs to a codimension one face. Set $B=A_j$.
When $i \neq j$, $c_{|A_i}$ belongs to the interior $\breve{C}_{A_i}(M)$ of $C_{A_i}(M)$, and since $c_{|A_i}$ is constant, $A_i$ contains a unique element 
and $c_{|A_i}$ does not map it to $\infty$. 
Two cases occur. Either $c_B(B)=\{\infty\}$ and $c \in F(\infty;B)$, or $c_B(B)=\{c_B(b)\} \subset (M \setminus \infty)$ and $c \in F(B)$.
Therefore the union of the codimension one faces is a subset 
of $\coprod_{F \in \partial_1(C_V(M))}F$. 
Conversely, Lemma~\ref{lempropconsadface} and the local product structure of
Lemma~\ref{lemimp} make clear that $\coprod_{F \in \partial_1(C_V(M))}F$ is a subset of the union of codimension one faces.
Now, it is clear that every $F \in \partial_1(C_V(M))$ is connected.
Furthermore, the closure of any such $F$ does not meet any other 
$F^{\prime}  \in \partial_1(C_V(M))$.\\
Let us prove this for $F=F(\infty;B)$.
In the closure of $F(\infty;B)$ all the configurations map $B$
to $\infty$ therefore $\overline{F(\infty;B)}$ may only meet the
$F(\infty;A)$ such that $B \subset A$. Consider a configuration $c$ in $(\overline{F(\infty;B)} \cap F(\infty;A))$. With the notation of Example~\ref{exa2},
since $c \in \overline{F(\infty;B)}$, $p_A(c)=c_A=\psi(A;\phi_{\infty})(\lambda, u \in S((\RR^3)^A))$, where $u$ maps $B$ to $0$;
then $\Pi_{\infty}(c_A)$
maps $B$ to $0$, but in this case $c \notin F(\infty;A)$.
A similar proof left to the reader leads to the same conclusion for $F=F(B)$.
Therefore, the $F$ are closed in the finite disjoint union $\coprod_{F \in \partial_1(C_V(M))}F$. Thus, they are the codimension one faces of $C_V(M)$,
and consequently, they smoothly embed in $C_V(M)$.
\eop

\noindent{\sc Proof of Proposition~\ref{propconffacedeux} assuming Lemma~\ref{lempropcons} and Lemma~\ref{lempropconsadface}:}
It is immediate from Lemma~\ref{lempropconsadface} and Proposition~\ref{propdifconfad}
that the disjoint union of the elements $f(B)(\RR^3)$ of $\partial_1(S_V(\RR^3))$
coincides with the union of the codimension one faces. A proof similar to the above one shows that the $f(B)$ are the connected components of this union.
Therefore, they are the codimension one faces of $S_V(\RR^3)$ and they smoothly embed there.
\eop

Proposition~\ref{propcdeuxcoinc} will be proved in Subsection~\ref{subpropcdeuxcoinc}.
Apart from Proposition~\ref{propcdeuxcoinc}, we are left with the proofs of Lemmas~\ref{lempropcons} and \ref{lempropconsadface} about the structure of $\tilde{C}_V(M)$ near a configuration $c^0$ such that
$\Pi_V(c^0_V)$ is the constant map $m^V$
with value $m$. The case where $(m \neq \infty)$ will be treated in the next subsection. The case $(m = \infty)$ is similar though more complicated, it will be treated in Subsection~\ref{subsecinf}, but some arguments will not be repeated.

\subsection{Proof of Proposition~\ref{propdifconf} when $\Pi_V(c^0_V)=m^V$, $m \in M \setminus \infty$.}
\label{subsecnoninf}

Let $\phi: \RR^3 \longrightarrow (M  \setminus \infty)$ be a smooth embedding, $\phi(0)=m$.
If $\sharp V=1$, set $k=0$, $W=\RR^3$, $w^0=0$, $U=\phi(\RR^3) \subset C_1(M)$, $\xi=\phi$ and $r=\phi^{-1}$, and we are done. Assume $\sharp V \geq 2$.
Let $c^0=(c^0_A)_{A \subseteq V; A \neq \emptyset} \in \tilde{C}_V(M)$ be such that $\Pi_V(c^0_V)=m^V$.

\medskip

\noindent{\bf The tree $\tau(c^0)$ associated to the limit configuration $c^0$.\/} \index{N}{tauc@$\tau(c^0)$}

\medskip
We shall define a set $\tau(c^0)$ of subsets of $V$ with cardinality $\geq 2$ as follows.

The set is organized as a tree with $(V \in \tau(c^0))$ as a root.
The other elements of $\tau(c^0)$ are constructed inductively as follows.
Every element $A$ of $\tau(c^0)$ is the {\em daughter\/} of its unique {\em mother\/} $\hat{A}$ in $\tau(c^0)$, except for $V$ that
has no mother, and some elements have {\em daughters\/} (i.e. are the mother of these).
A daughter is strictly included into its mother, and any two daughters are disjoint.
Therefore, it is enough to construct the daughters of an element $A$.
By assumption, 
$c^0_{A} \in \Pi_{A}^{-1}(\mbox{diag}(M \setminus \infty)^{A}) \subset
C(A;M)$.
Thus, $$\Pi_d(c^0_{A}) \in S\left(\frac{T_mM^{A}}{\mbox{diag}
(T_mM^{A})} \right)$$
defines a map from $A$ to $T_mM$ up to translation
and dilation.
The daughters of $A$ will be the preimages of multiple points.
The preimages of non-multiple points will be the {\em sons\/} of $A$.

$\tau(c^0)$ has the property that whenever $\{A,B\} \subseteq \tau(c^0)$,
either $A \subset B$, or $B \subset A$, or $A \cap B = \emptyset$.

Fix $c^0$, and $\tau=\tau(c^0)$.
For any $A \in \tau$ choose a basepoint $b(A)=b(A;\tau)$, \index{N}{bA@$b(A)$} 
such that
if $A \subset B$, if $B \in \tau$, and if $b(B) \in A$, then $b(A)=b(B)$.
When $A \in \tau$, $D(A)$ denotes the set of daughters of $A$.

\medskip

\noindent{\bf Configuration spaces associated to $\tau$.\/}

\medskip
For any $A \in \tau$, 
consider the following subsets of the unit sphere $S((\RR^3)^V)$ of $(\RR^3)^V$ equipped with its usual scalar product.
Define the set $C(A;b(A);\tau)$ \index{N}{CAbAtau@$C(A;b(A);\tau)$} 
of maps $w:V \longrightarrow \RR^3$ such that
\begin{itemize}
\item $\parallel w \parallel=1$
\item $w(b(A))=0$, $w(V \setminus A)=\{0\}$, and
\item $w$ is constant on any daughter of $A$.
\end{itemize}
It is easy to see that $C(A;b(A);\tau)$ \index{N}{CAbAtau@$C(A;b(A);\tau)$} 
has a canonical differentiable structure
(and is diffeomorphic to a sphere of dimension $\left(3 (\sharp A -\sum_{i=1}^n\sharp A_i +n-1)-1\right)$ where $A_1, \dots, A_n$ are the daughters of $A$.

Note that $c^0_A=\psi(A;\phi;b(A))(0;0;w^0_A)$ with the notation of Example~\ref{exa1}
where the natural extension $w^0_A$ (by some zeros) of 
$w^0_A \in (\RR^3)^{A \setminus b(A)} \subset (\RR^3)^V$ is in $C(A;b(A);\tau)$.

Define the set $O(A;b(A);\tau)$ \index{N}{OAbAtau@$O(A;b(A);\tau)$} 
of maps $w:V \longrightarrow \RR^3$ such that
\begin{itemize}
\item $\parallel w \parallel=1$
\item $w(b(A))=0$, $w(V \setminus A)=\{0\}$, and
\item  Two elements of $A$ that belong to different children (daughters and sons)
of $A$ are mapped to different points of $\RR^3$.
\end{itemize}

It is clear that $O(A;b(A);\tau)$ is an open subset of $S((\RR^3)^{A \setminus b(A)})$ that contains $w^0_A$.
Set 
$$W_A=O(A;b(A);\tau) \cap C(A;b(A);\tau)$$
$W_A$ is an open subset of the sphere $C(A;b(A);\tau)$.

\medskip

\noindent{\bf The data $U$, $W$, $w^0$ and $k$.\/}

\medskip
\begin{itemize}
\item $k=\sharp \tau$.
\item $\tilde{W}=\RR^3 \times \prod_{A \in \tau} W_A$
\item $W$ will be an open neighborhood of $w^0=(0;(w_A^0)_{A \in \tau})$ 
in $\tilde{W}$.
\item $\tilde{U}=\prod_{A \in \tau}\psi(A;\phi;b(A))\left(\RR^3 \times [0,\infty[ \times O(A;b(A);\tau)\right) \times \prod_{A \notin \tau}C(A;M)$
\item $U$ will be an open neighborhood of $c^0$ in $\tilde{U}$.
\end{itemize}

\medskip

\noindent{\bf Construction of $\xi$.\/}

\medskip
Let $$P=((\mu_A)_{A \in \tau};u; (w_A)_{A \in \tau}) \in \RR^{\tau} \times W$$
and $$P^0=((0)_{A \in \tau};0; (w^0_A)_{A \in \tau}) =(0;w^0).$$
When $A \in \tau$, define
$$v_A = v_A(P)=\sum_{C \in \tau; C \subseteq A}\left(\prod_{D \in \tau; C \subseteq D \subset A} \mu_D \right) w_C \in S((\RR^3)^{A \setminus b(A)}) \subset S((\RR^3)^V)$$

Note that $v_A$ is a smooth function defined on $\RR^k \times 
\tilde{W}$, and that $v_A(P^0)=w_A^0$.
In particular, $\parallel v_A(P^0) \parallel=1$ and $\frac{v_A(P^0)}{\parallel v_A(P^0) \parallel}$ is in $O(A;b(A);\tau)$.
Therefore, we can choose neighborhoods $O$ of $0$ in $[0,\infty[^k$
and $W$ of $w^0$ in $\tilde{W}$, so that for any $P$ in $O \times W$,   $\parallel v_A(P) \parallel \neq 0$ and 
$\frac{v_A(P)}{\parallel v_A(P) \parallel}$ is in $O(A;b(A);\tau)$.

We choose $O$ and $W$ so that these properties are satisfied for any $A \in \tau$.

In order to define $\xi$, we define its projections $\xi_A(P)$ onto the factors 
$C(A;M)$. First set
$$\xi_V(P)=\psi(V;\phi;b(V))(u;\mu_V;\frac{v_V}{\parallel v_V \parallel})$$
Then $$\Pi_V(\xi_V(P))(a)= \phi(u + \frac{\mu_V}{\parallel v_V \parallel}v_V(a))$$
When $A \in \tau$, set
$$\xi_A(P)=\psi(A;\phi;b(A))\left(u + \frac{\mu_V}{\parallel v_V \parallel}v_V(b(A));
\frac{{\parallel v_A \parallel}\prod_{D \in \tau; A \subseteq D \subseteq V} \mu_D}{\parallel v_V \parallel};\frac{v_A}{\parallel v_A \parallel}\right)$$

The latter definition makes sense because 
$v_{\hat{A}}$ is not constant on $A$ since $\frac{v_{\hat{A}}(P)}{\parallel v_{\hat{A}}(P) \parallel}$ belongs to $O(A;b(A);\tau)$.
Indeed, either $\Pi_{\hat{A}}(\xi_{\hat{A}}(P))$ is non constant and then its restriction to $A$ is non constant, and we take the usual smooth
restriction, or $\Pi_{\hat{A}}(\xi_{\hat{A}}(P))$ is constant with value $\phi(v)$, and we take the restriction of the map $D_v \phi \circ {v}_{\hat{A}}$ from $\hat{A}$ to $T_{\phi(v)}M$ up to translation and dilation.
It is easy to check that this restriction is smooth from this open subset of $C({\hat{A}};M)$ to $C(A;M)$, by using appropriate charts of $C(A;M)$ and $C(\hat{A};M)$ as in Example~\ref{exa1} with the same basepoint 
for $A$ and $\hat{A}$.
Thus, we defined a smooth map $\xi$ from $O \times W$ to $\tilde{U}$ such that $\xi(0;w^0)=c^0$.

\medskip

\noindent{\bf Checking that $\xi$ satisfies $(C1)$.  \/}

\medskip
It is enough to check that $\Pi_A(\xi(P)_A)= \Pi_V(\xi_V(P))_{|A}$ for any $A \in \tau$.
Let $A \in \tau$, $a \in A$.

$$\Pi_V(\xi_V(P))(a)=
\phi\left(u + \frac{\mu_V}{\parallel v_V \parallel}\sum_{C \in \tau; a \in C}\left(\prod_{D \in \tau; C \subseteq D \subset V} \mu_D \right) w_C(a)\right)$$
where the elements $C$ of $\tau$ that contain $a$, are 
\begin{enumerate}
\item the $C$ of $\tau$ such that $A \subset C$ that satisfy $w_C(a)=w_C(b(A))$, and
\item the $C$ of $\tau$ such that $a \in C \subseteq A$ that satisfy $w_C(b(A))=0$.
\end{enumerate}
In particular,$$u + \frac{\mu_V}{\parallel v_V \parallel}v_V(b(A))=
u + \frac{\mu_V}{\parallel v_V \parallel}\sum_{C \in \tau; A \subset C}\left(\prod_{D \in \tau; C \subseteq D \subset V} \mu_D \right) w_C(a).$$
Therefore,
$$\phi^{-1}(\Pi_V(\xi_V(P))(a))-\left(u + \frac{\mu_V}{\parallel v_V \parallel}v_V(b(A))\right) $$
$$=\frac{1}{\parallel v_V \parallel}\sum_{C \in \tau; a\in C; C \subseteq A}\left(\prod_{D \in \tau; C \subseteq D \subseteq V} \mu_D \right) w_C(a)$$
$$=\frac{\prod_{D \in \tau; A \subseteq D \subseteq V} \mu_D}{\parallel v_V \parallel}v_A(a).$$

\medskip

\noindent{\bf Checking that $\xi(O \times W) \subset \tilde{C}_V(M)$.  \/}

\medskip
It is enough to check that $\xi(P)$ satisfies $(C2)$ since $\Pi_V(\xi_V(P))(V) \subset (M \setminus \infty)$.
Let $A \subset B \subseteq \Pi_V(\xi_V(P))^{-1}(x)$.
If $B$ is not in $\tau$, then $\hat{B} \subseteq \Pi_V(\xi_V(P))^{-1}(x)$ 
(see the construction of $\xi_B(P)$), and $\xi_B(P)$ is the non-trivial restriction of $\xi_{\hat{B}}(P)$. Therefore, for this proof, we may assume that $B \in \tau$. Similarly,
we may assume that $A \in \tau$. Then it is enough to check that the restriction of $v_B$ to $A$ up to translation is a $(\geq 0)$ multiple of $v_A$, and this is easy to observe in the defining formula for $v_A$.

\medskip

\noindent{\bf Checking that $\xi((O \cap ]0,\infty[^k)\times W) = \breve{C}_V(M) \cap \xi(O \times W)$.  \/}

\medskip
Let us first prove that $\xi((O \cap ]0,\infty[^k)\times W) \subset \breve{C}_V(M)$. Since $(C1)$ is fulfilled in the image of $\xi$, it is enough to prove that $\Pi_V(\xi_V(P))$
is injective when the $\mu_A$ are non zero.
Let $a$ and $b$ be in $V$, and let $A$ be the smallest element of $\tau$ that contains both of them. Then $v_A(a) \neq v_A(b)$
since $\frac{v_A}{\parallel v_A \parallel}$ is in $O(A;b(A);\tau)$, thus $\Pi_A(\xi_A(P))$ separates $a$ and $b$, and we are done thanks to $(C1)$. Conversely, since as soon as a $\mu_A$ vanishes, the corresponding $\Pi_A(\xi_A(P))$ is constant, $$\breve{C}_V(M) \cap \xi(O \times W) \subseteq \xi((O \cap ]0,\infty[^k)\times W).$$

\medskip

\noindent{\bf Construction of $r$.\/}

\medskip
For any $A \in \tau$, choose 
$b^{\prime}(A) \neq b(A) \in A$ \index{N}{bprimeA@$b^{\prime}(A)$} 
to be either the element of a son of $A$ or a basepoint of a daughter of $A$ that does not contain $b(A)$. Note that $w^0_A(b^{\prime}(A)) \neq 0$.

The map $r$ will factor through the projection onto $$\prod_{A \in \tau}\psi(A;\phi;b(A))\left(\RR^3 \times [0,\infty[ \times O(A;b(A);\tau)\right).$$
Let $$Q=\left(\psi(A;\phi;b(A))(u_A; \lambda_A; y_A) \right)_{A \in \tau}$$ be a point of this space and let
$$r(Q)=((\mu_A)_{A \in \tau}; u_V; (w_A)_{A \in \tau})$$ denote its image in $\RR^k \times \tilde{W}$.
The map $u_V$ is already defined and smooth, and we need to define the $\mu_A$ and the $w_A$ as smooth functions of $(u_A; \lambda_A; y_A)_{A \in \tau}$.
Define
$w^1_A \in (\RR^3)^A$ by
$$w^1_A(a)=\left\{\begin{array}{ll} y_A(a) & \mbox{if}\; a \in \left(A \setminus
(\cup_{B \in D(A)}B) \right)\\
y_A(b(B)) & \mbox{if}\; a \in B\;\mbox{and if} \;B \in D(A) \end{array}\right..$$
 Then set
$$w_A=\frac{w^1_A}{\parallel w^1_A \parallel}$$
Since $y_A \in O(A;b(A);\tau)$, $\parallel w^1_A \parallel \neq 0$, and $w_A$ is smooth.
Then define $\mu_V=\lambda_V$, and for $A \in \tau$, $A \neq V$,

$$\mu_A=\frac{\parallel w_{\hat{A}}(b^{\prime}({\hat{A}})) \parallel  \langle y_{\hat{A}}(b^{\prime}(A))-y_{\hat{A}}(b(A)),w_A(b^{\prime}(A)) \rangle }{\parallel y_{\hat{A}}(b^{\prime}({\hat{A}})) \parallel \; \parallel w_A(b^{\prime}(A)) \parallel^2}$$

Then it is clear that $r$ is smooth from $\tilde{U}$ to $\RR^k \times \tilde{W}$.

\medskip

\noindent{\bf Checking that $r \circ \xi$ is the identity of $O \times W$.\/}

\medskip
We compute $$r \circ \xi\left(P=((\mu_A)_{A \in \tau};u; (w_A)_{A \in \tau})\right)=r((\xi_A(P))_{A \in \tau})$$
$$=((\hat{\mu}_A)_{A \in \tau}; u_V; (\hat{w}_A)_{A \in \tau}).$$
where 
$$\xi_V(P)=\psi(V;\phi;b(V))(u;\mu_V;\frac{v_V}{\parallel v_V \parallel}),$$  
$$\xi_A(P)=\psi(A;\phi;b(A))\left(u_A; \lambda_A;y_A=\frac{1}{\parallel v_A \parallel}v_A\right),$$
and $v_A$ has been defined in the construction of $\xi$.\\
We easily find $u=u_V$ and $\hat{\mu}_V=\mu_V$, and $\hat{w}_A=w_A$.
Since $$\left(v_{\hat{A}}(b^{\prime}(A)) - v_{\hat{A}}(b(A))\right)=\mu_A w_A(b^{\prime}(A))\;\;\;\mbox{and}\;\;\; v_{\hat{A}}=\parallel v_{\hat{A}} \parallel y_{\hat{A}}$$
$$\parallel v_{\hat{A}} \parallel\left(y_{\hat{A}}(b^{\prime}(A)) - y_{\hat{A}}(b(A))\right)=\mu_A w_A(b^{\prime}(A)).$$
Furthermore, 
$$y_{\hat{A}}(b^{\prime}({\hat{A}}))  =\frac{1}{\parallel v_{\hat{A}} \parallel}v_{\hat{A}}(b^{\prime}({\hat{A}})) =\frac{1}{\parallel v_{\hat{A}} \parallel}w_{\hat{A}}(b^{\prime}({\hat{A}})).$$
Therefore
$\hat{\mu}_A=\mu_A$.
Thus, $r \circ \xi$ is the identity of $O \times W$.

\medskip

\noindent{\bf Checking that $r(\tilde{U} \cap \tilde{C}_V(M)) \subseteq [0,\infty[^k \times \tilde{W}$.\/}

\medskip
When $Q=\left(\psi(A;\phi;b(A))(u_A; \lambda_A; y_A) \right)_{A \in \tau}$
comes from an element of $\tilde{C}_V(M)$, if $A$ and $B$ are two elements of $\tau$ such that $A \subset B$,
then for any $a \in A$,
$$u_B + \lambda_B y_B(a)=u_A + \lambda_A y_A(a),$$ and the map from $A$ to $\RR^3$ that maps $a$ to 
$(y_B(a)-y_B(b(A)))$ is a $(\geq 0)$ multiple of $y_A$.
$$r(Q)=((\mu_A)_{A \in \tau};u_V; (w_A)_{A \in \tau})$$ 

where 
$$\mu_A=\frac{\parallel w_{\hat{A}}(b^{\prime}({\hat{A}})) \parallel \; \parallel y_{\hat{A}}(b^{\prime}(A))-y_{\hat{A}}(b(A)) \parallel}{\parallel y_{\hat{A}}(b^{\prime}({\hat{A}})) \parallel \; \parallel w_A(b^{\prime}(A)) \parallel}$$
when $A \in \tau$, $A \neq V$.
Indeed, $\left((y_{\hat{A}})_{|A} -(y_{\hat{A}}(b(A)))^A\right)$ is a $(\geq 0)$ multiple of $y_A$, and $y_A(b^{\prime}(A))$  is a $(\geq 0)$ multiple of $w_A(b^{\prime}(A))$. In particular, $\mu_A \geq 0$, $\mu_V=\lambda_V$ is also positive.

Also, note that $r(c^0)=P^0$.

Now,  choose
$(\varepsilon  > 0)$ such that $[0,\varepsilon[^k \subset O$,
reduce $O$ into
$[0,\varepsilon[^k$ and set 
$$U=r^{-1}(]-\varepsilon,\varepsilon[^k \times W).$$

Then $r(U \cap \tilde{C}_V(M)) \subseteq O \times W$.

\medskip

\noindent{\bf Checking that $\xi \circ r_{|U \cap \tilde{C}_V(M)}$ is the identity of $U \cap \tilde{C}_V(M)$.\/}

Keep the above notation for $Q$ and $r(Q)$. Assume $Q \in U \cap \tilde{C}_V(M)$.
$$(\xi \circ r(Q))_A=\psi(A;\phi;b(A))(\tilde{u}_A;\tilde{\lambda}_A;\frac{v_A}{\norm{v_A}})$$
where $v_A$ is the vector associated to $r(Q)$ in the construction of $\xi$.

\medskip

\noindent{\em Proof that $y_A=\frac{v_A}{\norm{v_A}}$.}\\

\medskip
Let $b_0$ be an element of $A$. Inductively define 
$$B^1=\hat{\{b_0\}} \subset B^2=\hat{B^{1}} \subset \dots \subset B^{i+1}=\hat{B^{i}}  \subset \dots \subset B^k=A
.$$
Set $y_i= y_{B^i}$, $b_i=b(B^i)$, $b^{\prime}_i=b^{\prime}(B^{i})$ and $w_i= w_{B^i}$.
Then

$$v_A(b_0)=\sum_{i=1}^k \left( \prod_{j=i}^{k-1} 
\frac{\parallel w_{j+1}(b^{\prime}_{j+1}) \parallel \; \parallel y_{j+1}(b^{\prime}_j)-y_{j+1}(b_j) \parallel}{\parallel y_{j+1}(b^{\prime}_{j+1}) \parallel \; \parallel w_j(b^{\prime}_j) \parallel}
\right)w_i(b_0)$$

where
$$w_i(b_0)=
\frac{\parallel w_i(b^{\prime}_i) \parallel}{\parallel y_{i}(b^{\prime}_i) \parallel }
y_i(b_{i-1}).$$
Therefore

$$v_A(b_0)=\sum_{i=1}^k \parallel w_{A}(b^{\prime}(A)) \parallel \left( \frac{ \prod_{j=i}^{k-1} \parallel y_{j+1}(b^{\prime}_j)-y_{j+1}(b_j) \parallel}{\prod_{j=i}^{k}\parallel y_{j}(b^{\prime}_j) \parallel }
y_i(b_{i-1})\right)$$

while

$$y_A(b_0)= \sum_{i=1}^k(y_A(b_{i-1})-y_A(b_i)),$$
and since, for $i \leq j$, 
$$y_{j+1}(b_{i-1})-y_{j+1}(b_i)=\frac{\parallel y_{j+1}(b^{\prime}_j)-y_{j+1}(b_j)\parallel}{\parallel y_{j}(b^{\prime}_j) \parallel}\left(y_{j}(b_{i-1})-y_{j}(b_i)\right)
$$ 
$$y_A(b_{i-1})-y_A(b_i)=\prod_{j=i}^{k-1}\frac{\parallel y_{j+1}(b^{\prime}_j)-y_{j+1}(b_j)\parallel}{\parallel y_{j}(b^{\prime}_j) \parallel}
y_i(b_{i-1})$$

Thus, 
$$v_A(b_0)= \frac{\parallel w_{A}(b^{\prime}(A)) \parallel}{\parallel y_A(b^{\prime}(A)) \parallel}y_A(b_0)$$
and $y_A=\frac{v_A}{\norm{v_A}}$.
\eop

Note that $\tilde{\lambda}_V=\mu_V=\lambda_V$, $\tilde{u}_V=u_V$, and therefore 
$(\xi \circ r(Q))_V$ is the restriction of $Q$ to $V$ whose value in $M$ at $b(A)$ determines $\tilde{u}_A$ by $(C1)$ that is fulfilled in the image of $\xi$. Therefore $\tilde{u}_A=u_A$.

\medskip

\noindent{\em Proof that $\tilde{\lambda}_A ={\lambda}_A$ for $A \neq V$.}\\

\medskip
Now, let us compute $\tilde{\lambda}_A$ for $A \in \tau$, $A \neq V$.
Define
$$B^1=A \subset B^2=\hat{B^{1}} \subset \dots \subset B^{i+1}=\hat{B^{i}}  \subset \dots \subset B^k=V
.$$
Again, set $y_i= y_{B^i}$, $b_i=b(B^i)$, $b^{\prime}_i=b^{\prime}(B^{i})$ and $w_i= w_{B^i}$.
$$\tilde{\lambda}_A=\frac{\lambda_V\parallel v_A \parallel}{\parallel v_V \parallel}
\prod_{i=1}^{k-1} 
\left( 
\frac{\parallel w_{i+1}(b^{\prime}_{i+1}) \parallel \; \parallel y_{i+1}(b^{\prime}_i)-y_{i+1}(b_i) \parallel}{\parallel w_i(b^{\prime}_i) \parallel \; \parallel y_{i+1}(b^{\prime}_{i+1}) \parallel }
\right)$$
where
$$\parallel v_A \parallel= \frac{\parallel w_A(b^{\prime}(A)) \parallel}{\parallel y_A(b^{\prime}(A)) \parallel}.$$
Therefore
$$\tilde{\lambda}_A=\frac{\lambda_V\parallel y_V(b^{\prime}(V)) \parallel}{\parallel y_A(b^{\prime}(A)) \parallel}
\prod_{i=1}^{k-1} 
\left( 
\frac{ \parallel y_{i+1}(b^{\prime}_i)-y_{i+1}(b_i) \parallel}{\parallel y_{i+1}(b^{\prime}_{i+1}) \parallel }
\right)$$
where
$$
\lambda_{B^i}= 
\frac{\parallel y_{i+1}(b^{\prime}_i)-y_{i+1}(b_i) \parallel}
{ \parallel y_{i}(b^{\prime}_{i}) \parallel}
\lambda_{B^{i+1}}$$
$$\lambda_A= \lambda_V \prod_{i=1}^{k-1} 
\left( \frac{\parallel y_{i+1}(b^{\prime}_i)-y_{i+1}(b_i) \parallel}
{ \parallel y_{i}(b^{\prime}_{i}) \parallel} \right)$$
Thus $\tilde{\lambda}_A=\lambda_A$.

When $A \notin \tau$, $(\xi \circ r(Q))_A$ is the restriction of $Q_{\hat{A}}$ to $A$, and we can conclude that the restriction of $\xi \circ r$ to ${U} \cap \tilde{C}_V(M)$ is the identity.
\eop

This concludes the proof of Proposition~\ref{propdifconf} in this case.
\eop

Proposition~\ref{propdifconfad} follows from a careful reading of the previous proof. Since $V \in \tau$, $k=\sharp V \geq 1$.
Choose the natural embedding $$\phi: \RR^3 \longrightarrow S^3 = \RR^3 \cup \{\infty\}.$$ 
The elements of ${S}_V(\RR^3) \cap U$ (resp. $\breve{S}_V(\RR^3) \cap U$) are the elements whose projection onto $C(V;M)$ is of the form
$\psi(V;\phi;b(V))(0 \in \RR^3;\mu_V=0;y_V)$ for some $y_V$ (resp. for some injective $y_V$). In particular, the second item is true where $\mu_V$ is the distinguished real parameter that vanishes if and only if $\Pi_V(c_V)$ is constant.
Now, since
$v_V$ is injective if and only if all the $\mu_{D}$ are non-zero for $D \in \tau \setminus V$, the third item is true. 
Proposition~\ref{propdifconfad} is proved. \eop

In this case, Lemma~\ref{lempropconsadface} also follows from a careful reading of the previous proof. Indeed, $k=1$ if and only if $\tau=\{V\}$, that is if and only
if $c_V^0 \in \breve{S}_V(T_{c(b)}M)$, that is if and only if $c^0 \in F(V)$.
Now, $k=2$ if and only if $\tau=\{V,B\}$ for some strict subset $B$ of $V$
with $\sharp B \geq 2$, that is if and only if :\\
$\Pi_d(c^0_V)$ is constant on $B$
and injective on $\{b\} \cup (V \setminus B)$, and $\Pi_d(c^0_B)$ is injective.\\
This is equivalent to say that under the assumptions of Lemma~\ref{lempropconsadface}, $c^0 \in f(B)(\RR^3)$. 
Lemma~\ref{lempropconsadface} is proved in this case.\eop

\subsection{Proof of Proposition~\ref{propdifconf} when $\Pi_V(c^0_V)=\infty^V$.}
\label{subsecinf}

Let $\phi_{\infty}: \RR^3 \longrightarrow M$ be a smooth embedding, $\phi_{\infty}(0)=\infty$.
Let $c^0=(c^0_A)_{A \subseteq V; A \neq \emptyset} \in \tilde{C}_V(M)$ be such that $\Pi_V(c^0_V)$ maps every point of $V$ to $\infty$.

\medskip

\noindent{\bf The tree $\tau(c^0)$ associated to $c^0$.\/} \index{N}{tauc@$\tau(c^0)$}

\medskip
We shall define a set $\tau=\tau(c^0)$ of non-empty subsets of $V$ as follows.

The set is organized as a tree with $V$ as a root.
The other elements of $\tau$ are constructed inductively as follows.
Again, every element $A$ of $\tau$ is the {\em daughter\/} of its unique {\em mother\/} $\hat{A}$ in $\tau$, except for $V$ that
has no mother, and some elements have {\em daughters\/} (i.e. are the mother of these).
In order to define the daughters of $A \in \tau$, consider the map defined up to dilation \index{N}{Piinfty@$\Pi_{\infty}(c^0_{A})$}
$$\Pi_{\infty}(c^0_{A}):A \longrightarrow T_{\infty}(M).$$
\begin{itemize}
\item If this map $\Pi_{\infty}(c^0_{A})$
 is non-constant, or if $A$ has only one element, then $A$ is {\em non-degenerate\/}.
In this case, let $A_0$ denote the preimage of $\{0\}$ under $\Pi_{\infty}(c^0_{A})$. If $A_0$ is non-empty,
$A_0$ is a daughter of $A$ and this daughter is said to be {\em special;\/}
the other daughters of $A$ will be the preimages of multiple points different from $0$ under $\Pi_{\infty}(c^0_{A})$.
The preimages of non-multiple points different from zero will be the {\em sons\/} of $A$.
\item If the map $\Pi_{\infty}(c^0_{A})$ is constant, and if $\sharp A \geq 2$, then $A$ is {\em degenerate\/}, and we consider the non-constant map defined up to translation and dilation \index{N}{Piinftyd@$\Pi_{\infty,d}(c^0_{A})$}
$$\Pi_{\infty,d}(c^0_{A}):A \longrightarrow T_{\infty}(M).$$
The daughters of $A$ are the preimages of multiple points under this map, and its sons are the preimages of the other points.
\end{itemize}
By definition $V$ is special, and an element $A \neq V$ of $\tau$ is special
if and only if $\Pi_{\infty}(c^0_{\hat{A}})(A)=\{0\}$.

Let $\tau_d$ \index{N}{taud@$\tau_d$} 
be the set of the degenerate elements of $\tau$, 
and let $\tau_s$ \index{N}{taus@$\tau_s$} 
be the set of the special elements of $\tau$. When $A \in \tau$, $D(A)$ denotes the set of daughters of $A$. Note that 
\begin{itemize}
\item $V \in \tau_s$,
\item $D(A \in \tau_d) \subset \tau_d \cap (\tau \setminus \tau_s)$
\item $D(A \notin \tau_d) \subseteq (\tau_d \cup \{A_0\})$,
\item $\tau =\tau_s \cup \tau_d$
\item If $A \neq V$, $A \in \tau_s$, then $\hat{A} \notin \tau_d$, therefore $\hat{A} \in \tau_s$.
\end{itemize}
In particular, 
$$\tau_s=\{V=V(1),V(2), \dots, V(\sigma)\}$$
where $V(i)_0=V(i+1) \neq \emptyset$ if $i <  \sigma$, and $V(\sigma)_0= \emptyset$.
Also note that $\tau_s \cap \tau_d \subseteq \{V(\sigma)\}$.
\index{N}{taud@$\tau_d$} \index{N}{taus@$\tau_s$}

Fix $c^0$, and $\tau=\tau(c^0)$.
For any $A \in \tau$ choose a basepoint 
$b(A)=b(A;\tau)$, \index{N}{bA@$b(A)$} 
such that
\begin{itemize} 
\item $b_i=b(V(i))=b(V(\sigma))$ for any $i=1, \dots, \sigma$, and,
\item if $A \subset B$, if $B \in \tau$, and if $b(B) \in A$, then $b(A)=b(B)$. 
\end{itemize}

\medskip

\noindent{\bf Configuration spaces associated to $\tau$.\/}

\medskip
Let $i \in \{1, \dots, \sigma\}$. 
Define the smooth manifold $C(V(i);\tau)$ \index{N}{CVitau@$C(V(i);\tau)$}
as the following submanifold of the unit sphere $S((\RR^3)^V)$ of $(\RR^3)^V$ equipped with its usual scalar product.
The set $C(V(i);\tau)$ is the set of maps $w:V \longrightarrow \RR^3$ such that
\begin{itemize}
\item $\parallel w \parallel=1$
\item $w(V(i)_0)=\{0\}$, $w(V \setminus V(i))=\{0\}$, and
\item \begin{itemize}
\item if $V(i) \notin \tau_d$, $w$ is constant on any daughter of $V(i)$,  and, 
\item if $V(i) \in \tau_d$, $w$ is constant.
\end{itemize}
\end{itemize}

Define the open subset $O(V(i);\tau)$ \index{N}{OVitau@$O(V(i);\tau)$} 
of $S((\RR^3)^{V(i)})$ as
the set of maps $w:V \longrightarrow \RR^3$ such that
\begin{itemize}
\item $\parallel w \parallel=1$
\item $w(V \setminus V(i))=\{0\}$,
\item If $V(i) \notin \tau_d$, two elements of $V(i)$ that belong to different children (daughters and sons)
of $V(i)$ are mapped to different points of $\RR^3$, and
\item $0 \notin w(V(i) \setminus V(i)_0)$.
\end{itemize}
Set $W^s_i=O(V(i);\tau) \cap C(V(i);\tau)$. $W^s_i$ is an open submanifold of the sphere $C(V(i);\tau)$.

Then after a proper scalar multiplication, the natural extension $s^0_i$ (by some zeros) of $(D_0\phi_{\infty})^{-1} \circ \Pi_{\infty}(c^0_{V(i)})$ is in $W^s_i$.

For any $A \in \tau_d$, consider the smooth manifold $C(A;b(A);\tau)$ and the open subset $O(A;b(A);\tau)$ of $S((\RR^3)^{A \setminus b(A)})$  defined as in Subsection~\ref{subsecnoninf}. Set 
$$W_A=O(A;b(A);\tau) \cap C(A;b(A);\tau).$$

Then after a proper normalization, the natural extension $w^0_A$ (by some zeros) of $(D_0\phi_{\infty})^{-1} \circ \Pi_{\infty,d}(c^0_A)$ is in $W_A$.


\medskip

\noindent{\bf The data $U$, $W$, $w^0$ and $k$.\/}

\medskip
\begin{itemize}
\item $k= \sharp \tau_s + \sharp \tau_d= \sigma +  \sharp \tau_d$.
\item $\tilde{W}=\prod_{i=1}^{\sigma} W^s_i \times \prod_{A \in \tau_d} W_A$
\item $W$ will be an open neighborhood of $w^0=((s_i^0)_{i \in \{1, \dots, \sigma\}};(w_A^0)_{A \in \tau_d})$ 
in $\tilde{W}$.
\item $\tilde{U}=\prod_{A \in \tau_d}\psi(A;\phi_{\infty};b(A))\left([0,\infty[
\times S^2 \times [0,\infty[ \times O(A;b(A);\tau)\right) $\\
$\times 
\prod_{A \in (\tau \setminus \tau_d)}\psi(A;\phi_{\infty})\left([0,\infty[
\times O(A;\tau)\right) \times \prod_{A \notin \tau}C(A;M)$\\
(with the charts of Examples~\ref{exa2} and \ref{exa3}).
\item $U$ will be an open neighborhood of $c^0$ in $\tilde{U}$.
\end{itemize}

\noindent{\bf Forgetting $\prod_{A \notin \tau}C(A;M)$.\/}

When $A \notin \tau$, let $\hat{A}$ be the smallest element in $\tau$ that contains
$A$. Let $C=\prod_{A \subseteq V, A \neq \emptyset}C(A;M)$, $C^{\tau}=\prod_{A \in \tau}C(A;M)$ and let $C^{\tau}_V(M)$ be the subspace of $C^{\tau}$ made of the elements that satisfy the restriction conditions $(C1)$, $(C2)$, $(C3)$, $(C4)$ of Lemmas~\ref{lemcond1},
\ref{lemcond2}, \ref{lemcond3}, \ref{lemcond4} that involve elements of $\tau$. 
Let $p^{\tau}: C \longrightarrow C^{\tau}$ be the natural projection. 
Define the following smooth map  $$\iota^{\tau}: p^{\tau}(\tilde{U}) \longrightarrow C$$ 
by 
$\iota^{\tau}(c=(c_A)_{A \in \tau})=(d_A)_{A \subseteq V, A \neq \emptyset}$, where $d_A=c_A$ when $A \in \tau$, and $d_A$ is the restriction of $c_{\hat{A}}$ to $A$
otherwise (so that the restriction conditions $(C1)$, $(C2)$, $(C3)$, $(C4)$
are satisfied for $(A,\hat{A})$). Note that such a restriction is well-defined and smooth
from $p_{\hat{A}}(\tilde{U})$ to $C(A;M)$ since $A$ is not contained
a daughter of $\hat{A}$. See the charts of Examples~\ref{exa2} and \ref{exa3}. In particular, $\iota^{\tau}$ is smooth. The proofs of the following assertions are left to the reader.

\begin{itemize}
\item $p^{\tau}(\tilde{C}_V(M)) \subseteq C^{\tau}_V(M)$.
\item $p^{\tau} \circ \iota^{\tau}_{|p^{\tau}(\tilde{U})}=\mbox{Identity}(p^{\tau}(\tilde{U}))$
\item $\iota^{\tau} \circ p^{\tau}_{|\tilde{U} \cap \tilde{C}_V(M)}=\mbox{Identity}(\tilde{U} \cap \tilde{C}_V(M))$.
\item $\iota^{\tau}(p^{\tau}(\tilde{U}) \cap C^{\tau}_V(M)) \subset \tilde{C}_V(M)$.
\end{itemize}

We shall prove the following lemma.

\begin{lemma}
\label{lemdifconfinfini}
There exist
\begin{enumerate}
\item $\varepsilon>0$,
$O= [0,\varepsilon[^k$,
\item an open neighborhood $W$ of $w^0$ in $\tilde{W}$,
\item an open neighborhood $U^{\tau}$ of $p^{\tau}(c^0)$ in $p^{\tau}(\tilde{U})$,
\item a smooth map $(p^{\tau} \circ \xi): O \times W \longrightarrow U^{\tau}$ such that 
\begin{itemize}
\item $(p^{\tau} \circ \xi)(0;w^0)=p^{\tau}(c^0)$,
\item $(p^{\tau} \circ \xi)(O \times W) \subset {C}^{\tau}_V(M)$, and 
\item $p_V \circ \xi(\omega \in O,w)$ is an injective map from $V$ to $(M \setminus \infty)$
if and only if $\omega \in ]0,\infty[^k$,
\end{itemize}
\item a smooth map $r^{\tau}: p^{\tau}(\tilde{U}) \longrightarrow \RR^k \times \tilde{W}$ such that 
\begin{itemize}
\item $r^{\tau} \circ p^{\tau} \circ \xi$ is the identity of $O\times W$, 
\item $r^{\tau}(U^{\tau} \cap {C}^{\tau}_V(M)) \subseteq O \times W$, and 
\item the restriction of $p^{\tau} \circ \xi \circ r^{\tau}$ to $U^{\tau} \cap C^{\tau}_V(M)$ is the identity of $U^{\tau} \cap C^{\tau}_V(M)$.
\end{itemize}
\end{enumerate}
\end{lemma}

This lemma implies Proposition~\ref{propdifconf} in this case because $\xi=\iota^{\tau} \circ (p^{\tau} \circ \xi)$, $U=(p^{\tau})^{-1}(U^{\tau})$ and
$r=r^{\tau} \circ p^{\tau}$ have the desired properties under its conclusions.


\medskip

\noindent{\bf Construction of $p^{\tau} \circ \xi$, $O$ and $W$.\/}

\medskip
\noindent Set $$P=((\nu_i)_{i \in \{1, \dots, \sigma\}};(\mu_A)_{A \in \tau_d};(s_i)_{i \in \{1, \dots, \sigma\}};(w_A)_{A \in \tau_d}) \in \RR^{\tau_s} \times \RR^{\tau_d} \times \tilde{W},$$
$$P^0=((0)_{i \in \{1, \dots, \sigma\}};(0)_{A \in \tau_d};(s^0_i)_{i \in \{1, \dots, \sigma\}}; (w^0_A)_{A \in \tau_d}) )=(0;w^0).$$
When $A \in \tau_d$, define
$$\tilde{w}_A =\tilde{w}_A(P)=\sum_{C \in \tau; C \subseteq A}\left(\prod_{D \in \tau; C \subseteq D \subset A} \mu_D \right) w_C \in S((\RR^3)^{A \setminus b(A)}) \subset S((\RR^3)^V).$$
$$\tilde{w}_A=w_A + \sum_{C \in D(A)}\mu_C \tilde{w}_C.$$
Set
$$\begin{array}{ll}
\tilde{s}_{\sigma} =s_{\sigma} + \mu_{V(\sigma)} \tilde{w}_{V(\sigma)} & \mbox{if}\; V(\sigma) \in \tau_d\\
\tilde{s}_{\sigma} =s_{\sigma} + \sum_{C \in D(V(\sigma))}\mu_C \tilde{w}_C & \mbox{otherwise.}
\end{array}$$
For $i= \sigma-1, \sigma-2, \dots, 1$, inductively define
$$\tilde{s}_i =s_i + \nu_{i+1} \tilde{s}_{i+1} +\sum_{C \in D(V(i)); C \neq V(i+1)}\mu_C \tilde{w}_C.$$
Define
$$\lambda_{V(r)}=\lambda_r=\prod_{i=1}^r \nu_i$$
so that $(\lambda_i \tilde{s}_i)_{|V(i+1)}=\lambda_{i+1} \tilde{s}_{i+1}$.

The $\tilde{w}_A$, $\lambda_r$, and $\tilde{s}_i$ are smooth functions defined on $\RR^k \times \tilde{W}$, such that $\tilde{w}_A(P^0)=w_A^0$ and $\tilde{s}_i(P^0)=s_i^0$.

In particular, since the norms of these vectors are $1$ for $P^0$, we can choose neighborhoods $O$ of $0$ in $[0,\infty[^k$
and $W$ of $w^0$ in $\tilde{W}$, so that for any $P$ in $O \times W$
\begin{itemize}
\item the norms of the $\tilde{w}_A(P)$ and $\tilde{s}_i(P)$ do not vanish,
\item $\frac{\tilde{w}_A(P)}{\parallel \tilde{w}_A(P) \parallel} \in O(A;b(A);\tau)$, and,
\item $\frac{\tilde{s}_i(P)}{\parallel \tilde{s}_i(P) \parallel} \in O(V(i);\tau)$.
\end{itemize}

We choose $O$ and $W$ so that these properties are satisfied for any $A \in \tau_d$, and for any $i=1,2,\dots,\sigma$.

When $A \in \tau$, let ${V(i(A))}$ be the smallest element of $\tau_s$
such that $A \subseteq V(i(A))$.
Define $\tilde{s}_A \in (\RR^3)^A$ as the restriction of $\tilde{s}_{i(A)}$ to $A$.

In order to define $p^{\tau} \circ \xi$, we define its projections $\xi_A(P)$ onto the factors 
$C(A;M)$ for $A \in \tau$.\\
When ${A} \in \tau \setminus \tau_d$, 
$\xi_A(P)=\psi(A;\phi_{\infty})(\lambda_{i(A)} \parallel \tilde{s}_A \parallel;
\frac{\tilde{s}_A}{\parallel \tilde{s}_A \parallel}).$\\
When ${A} \in \tau_d$,
$\xi_A(P)=\psi(A;\phi_{\infty};b(A))(\ell_A;u_A;m_A;v_A)$
with
$${\ell}_A=
\lambda_{i(A)} \sqrt{\sharp A}\parallel \tilde{s}_{A}(b(A)) \parallel$$
$${u}_A=\frac{\tilde{s}_{A}(b(A))}{\parallel \tilde{s}_{A}(b(A)) \parallel}$$
$${m}_A=\left(\prod_{D \in \tau_d; A \subseteq D \subseteq V(i(A))} \mu_D \right)
\frac{\parallel \tilde{w}_A \parallel}{\sqrt{\sharp A}\parallel \tilde{s}_{A}(b(A)) \parallel}$$ 
$${v}_A=\frac{\tilde{w}_A}{\parallel \tilde{w}_A \parallel}.$$
Thus, we defined a smooth map $p^{\tau} \circ \xi$ from $O \times W$ to $p^{\tau} (\tilde{U})$.

\medskip

\noindent{\bf Checking that $p^{\tau} \circ \xi$ satisfies $(C1)$.  \/}

\medskip
Since the restriction of $\lambda_{i}  \tilde{s}_i$ to $V(i+1)$ is
$\lambda_{i+1}  \tilde{s}_{i+1}$, when $i \leq \sigma -1$, it is enough to check that for any $A$, $\Pi_A(\xi_A(P))$ is equal to 
$\phi_{\infty} \circ \left(\lambda_{i(A)}  (\tilde{s}_{i(A)})_{|A}\right)$.
When ${A} \notin \tau_d$
it is obvious.
Let us now consider the case when ${A} \in \tau_d$. \\
$$\tilde{s}_{A}=\mbox{constant map} + \left(\prod_{D \in \tau_d; A \subseteq D \subseteq V(i(A))} \mu_D \right) \tilde{w}_{{A}}.$$
Therefore
$$\phi_{\infty}^{-1} \circ \left(\Pi_A(\xi_A(P))\right)$$
$$=\lambda_{i(A)}\left( \left(\tilde{s}_{A}(b(A))\right)^A + 
\left(\prod_{D \in \tau_d; A \subseteq D \subseteq V(i(A))} \mu_D \right)
\tilde{w}_A
\right)$$
$$=\lambda_{i(A)}\tilde{s}_{A}.$$

\medskip

\noindent{\bf Checking that $p^{\tau} \circ \xi$ satisfies $(C3)$.  \/}

\medskip
Here, we need to check that when $A \subset B$, (and when $\lambda_{i(B)}=0$)
$(\tilde{s}_{i(B)})_{|A}$ is a $(\geq 0)$ multiple of $(\tilde{s}_{i(A)})_{|A}$.
Since $(\tilde{s}_{i(B)})_{|V(i(A))}=\left(\prod_{j=i(B)+1}^{i(A)}\nu_j \right) \tilde{s}_{i(A)}$, we are done.

\medskip

\noindent{\bf Checking that $p^{\tau} \circ \xi$ satisfies $(C2)$ and $(C4)$.  \/}

\medskip
These conditions must be checked for some $A \subset B$, when the restriction of $\tilde{s}_{i(B)}$ to $B$ is constant.
In this case, since $\frac{\tilde{s}_{i(B)}}{\norm{\tilde{s}_{i(B)}}} \in O(V(i(B));\tau)$, ${B} \in \tau_d$, and 
therefore
${A} \in \tau_d$.
These conditions say that, up to translation, $(\tilde{w}_{{B}})_{|A}$ is a $(\geq 0)$ multiple of $(\tilde{w}_{{A}})$
(when $\left(\prod_{D \in \tau_d; B \subseteq D \subseteq V(i(B))} \mu_D \right)=0$). They are realised with $\left(\prod_{D \in \tau_d; A \subseteq D \subset {B}} \mu_D \right)=0$ as a factor.

\medskip

\noindent {\bf We have proved that $p^{\tau} \circ \xi(O \times W) \subset {C}^{\tau}_V(M)$}
and it is easy to see that $p^{\tau} \circ \xi(P^0=(0;w^0))=p^{\tau}(c^0)$.

\medskip

\noindent{\bf Checking that $p_V \circ \xi((\mu_A,\nu_i) \in O,w \in W)$ is injective and does not reach $\infty$
if and only if all the $\mu_A$ and the $\nu_i$ are non zero. \/}

\medskip
Remember from the proof that $p^{\tau} \circ \xi$ satisfies $(C1)$, that the 
restriction $\Pi_A(\xi_A(P))$ of $\Pi_V(\xi_V(P))$ is
$\phi_{\infty} \circ \left(\lambda_{i(A)}  (\tilde{s}_{i(A)})_{|A}\right)$ for $A \in \tau$.
Now, $\Pi_V(\xi_V(P))$ is injective and does not reach $\infty$
if and only if $\lambda_{1}\tilde{s}_{1}$ is injective and does not reach $0$.

In particular, if $\Pi_V(\xi_V(P))$ is injective and does not reach $\infty$,
all the restrictions $\lambda_{i(A)}  (\tilde{s}_{i(A)})_{|A}$ are injective and do not reach $0$ and this easily implies that the $\mu_A$ and the $\nu_i$ are non zero.

Conversely, assume that the $\mu_A$ and the $\nu_i$ are non zero, and let us prove that $\lambda_{1}\tilde{s}_{1}$ is injective and does not reach $0$. Since $\nu_1=\lambda_1$, it is enough to prove that $\tilde{s}_{1}$ is injective  and does not reach $0$.
Let $a$ and $b$ be in $V$, and let $A$ be the smallest element of $\tau$ that contains both of them. If $A \in \tau_d$, $\tilde{w}_A(a) \neq \tilde{w}_A(b)$ and $\tilde{s}_{i(A)}(a) \neq \tilde{s}_{i(A)}(b)$. If $A \notin \tau_d$, $A=V(i(A))$, and $\tilde{s}_{i(A)}(a) \neq \tilde{s}_{i(A)}(b)$. Since $\tilde{s}_{1\,|V(i(A))}$ is a non zero multiple
of $\tilde{s}_{i(A)}$, it separates $a$ and $b$, and $\tilde{s}_{1}$ is injective. If $\tilde{s}_{1}(a)=0$, then $\tilde{s}_{i(\{a\})}(a)=0$, and this
is impossible, therefore $\tilde{s}_{1}$ does not reach $0$.

\medskip

\noindent{\bf Construction of $r^{\tau}$.\/}

\medskip

Let $c \in p^{\tau}(\tilde{U})$.
$$c=\left(\left(\psi(A;\phi_{\infty};b(A))(\ell_A;u_A;m_A;v_A)\right)_{A \in \tau_d};\left(\psi(A;\phi_{\infty})(\ell_A;S_A)\right)_{A \in (\tau \setminus \tau_d)}
\right).$$
We shall define
$$r^{\tau}(c)=((\nu_i)_{i \in \{1, \dots, \sigma\}};(\mu_A)_{A \in \tau_d};(s_i)_{i \in \{1, \dots, \sigma\}};(w_A)_{A \in \tau_d}).$$

\medskip

\noindent{\em Definition of $w_A$, for $A \in \tau_d$.}\\
Let $A \in \tau_d$. Define $w^1_A \in (\RR^3)^A$ by
$$w^1_A(a)=\left\{\begin{array}{ll} v_A(a) & \mbox{if}\; a \in \left(A \setminus
(\cup_{B \in D(A)}B) \right)\\
v_A(b(B)) & \mbox{if}\; a \in B\;\;\;\;\mbox{with} \;\;\;\;B \in D(A) \end{array}\right.$$
and set 
$$w_A=\frac{w^1_A}{\parallel w^1_A \parallel}.$$

\medskip

\noindent{\em Definition of $s_i$, for $i \in \{1, \dots, \sigma\}$.}\\
Let $V(i) \notin \tau_d$. Define $s^1_i \in (\RR^3)^{V(i)}$ by
$$s^1_i(a)=\left\{\begin{array}{ll} 0 &  \mbox{if}\; a \in V(i)_0 \\
S_{V(i)}(a) & \mbox{if}\; a \in \left(V(i) \setminus
(\cup_{B \in D(V(i))}B) \right)\\
S_{V(i)}(b(B)) & \mbox{if}\; a \in B\;\mbox{where} \;B \in D(V(i))\;\mbox{and} \;B \neq  V(i)_0 \end{array}\right.$$
and set 
$$s_{V(i)}=s_i=\frac{s^1_i}{\parallel s^1_i \parallel}.$$
If $V(\sigma) \in \tau_d$, then
$$s_{\sigma}=\left(\frac{u_{V(\sigma)}}{\sqrt{\sharp V(\sigma)}}\right)^{V(\sigma)}.$$

\medskip

\noindent{\em Definition of $\mu_A$, for $A \in \tau_d$.}\\
For any $A \in \tau$ such that $(\sharp A \geq 2)$, choose 
$b^{\prime}(A) \neq b(A) \in A$ \index{N}{bprimeA@$b^{\prime}(A)$} 
to be either the element of a son of $A$ or a basepoint of a daughter of $A$ that does not contain $b(A)$. \\
\begin{itemize}

\item If $A \in \tau_d$, (if $A \notin \tau_s$,) and if $\hat{A} \in \tau_d$,
$$\mu_A=\frac{ \langle v_{\hat{A}}(b^{\prime}(A))-v_{\hat{A}}(b(A)), w_{A}(b^{\prime}(A)) \rangle }
{ \langle  w_{A}(b^{\prime}(A)), w_{A}(b^{\prime}(A)) \rangle }
\frac{\parallel w_{\hat{A}}(b^{\prime}(\hat{A})) \parallel}{\parallel v_{\hat{A}}(b^{\prime}(\hat{A})) \parallel}.$$
\item  If $A \in \tau_d$, if $A \notin \tau_s$, and if $\hat{A} \notin \tau_d$,
$$\mu_A=\frac{ \langle S_{\hat{A}}(b^{\prime}(A))-S_{\hat{A}}(b(A)), w_{A}(b^{\prime}(A)) \rangle }
{ \langle  w_{A}(b^{\prime}(A)), w_{A}(b^{\prime}(A)) \rangle }
\frac{\parallel s_{\hat{A}}(b^{\prime}(\hat{A})) \parallel}{\parallel S_{\hat{A}}(b^{\prime}(\hat{A})) \parallel}.$$
\item If $V(\sigma) \in \tau_d$, then
$\mu_{V(\sigma)}=\frac{m_{V(\sigma)}}{\parallel \tilde{w}_{V(\sigma)}\parallel}$
with
$$\tilde{w}_{V(\sigma)}=\sum_{C \in \tau; C \subseteq V(\sigma)}\left(\prod_{D \in \tau; C \subseteq D \subset V(\sigma)} \mu_D \right) w_C.$$
\end{itemize}

\medskip

\noindent{\em Definition of $\nu_i$, for $i \in \{1, \dots, \sigma\}$.}\\
Set  $S_i=S_{V(i)}$ when $V(i) \notin \tau_d$, and set $b^{\prime}_i=b^{\prime}(V(i))$ when $\sharp V(i) > 1$.
\begin{itemize}
\item
When $i \geq 2$, 
\begin{itemize}
\item If $i=\sigma$, and, if $V(\sigma) \in \tau_d$ or if $\sharp V(\sigma)=1$, then
$$\nu_{i}=\frac{ \langle S_{i-1}(b_{\sigma}),s_{i}(b_{\sigma}) \rangle }{ \langle s_{i}(b_{\sigma}),s_{i}(b_{\sigma}) \rangle }
\frac{\parallel s_{i-1}(b^{\prime}_{i-1}) \parallel}{\parallel S_{i-1}(b^{\prime}_{i-1})\parallel}.$$
\item Otherwise,
$$\nu_{i}=\frac{ \langle S_{i-1}(b^{\prime}_{i}),s_{i}(b^{\prime}_{i}) \rangle }{ \langle s_{i}(b^{\prime}_{i}),s_{i}(b^{\prime}_{i}) \rangle }
\frac{\parallel s_{i-1}(b^{\prime}_{i-1}) \parallel}{\parallel S_{i-1}(b^{\prime}_{i-1})\parallel}$$
\end{itemize}
\item \begin{itemize}
\item If $V \in \tau_d$, or if $\sharp V=1$, $\nu_1=\ell_V$.
\item If $V \notin \tau_d$ and if  $\sharp V > 1$,
$$\nu_1=\ell_V\frac{ \langle S_V(b^{\prime}_{1}),s_1(b^{\prime}_{1}) \rangle }{ \langle s_1(b^{\prime}_{1}),s_1(b^{\prime}_{1}) \rangle }.$$

\end{itemize}
\end{itemize}

Then it is clear that $r^{\tau}$ is smooth from $p^{\tau}(\tilde{U})$ to $\RR^k \times \tilde{W}$.

\medskip

\noindent{\bf Checking that $r^{\tau} \circ p^{\tau} \circ \xi$ is the identity  of $O \times W$.\/}

\medskip
We compute $$r^{\tau} \circ p^{\tau} \circ \xi(P=((\nu_i)_{i \in \{1, \dots, \sigma\}};(\mu_A)_{A \in \tau_d};(s_i)_{i \in \{1, \dots, \sigma\}};(w_A)_{A \in \tau_d}))$$
$$=((\nu^{\prime}_i)_{i \in \{1, \dots, \sigma\}};(\mu^{\prime}_A)_{A \in \tau_d};(s^{\prime}_i)_{i \in \{1, \dots, \sigma\}};(w^{\prime}_A)_{A \in \tau_d}).$$
It is clear that  $w^{\prime}_A=w_A$ for any $A \in \tau_d$ and that $s^{\prime}_i=s_i$ if $V(i) \notin \tau_d$.\\
If $V(\sigma) \in \tau_d$, then $s_{\sigma}$ is constant and $\frac{\tilde{s}_{\sigma}(b_{\sigma})}{\parallel \tilde{s}_{\sigma}(b_{\sigma})\parallel}= \frac{s_{\sigma}(b_{\sigma})}{\parallel s_{\sigma}(b_{\sigma})\parallel}$.\\ Thus
$u_{V(\sigma)}= \frac{s_{\sigma}(b_{\sigma})}{\parallel s_{\sigma}(b_{\sigma})\parallel}$ and $s^{\prime}_{\sigma}=\left(\frac{u_{V(\sigma)}}{\sqrt{\sharp V(\sigma)}}\right)^{V(\sigma)}
=s_{\sigma}$.\\

\noindent{\em Checking that $\mu^{\prime}_A=\mu_A$, for $A \in \tau_d$.}\\
\begin{itemize}
\item If $A \in \tau_d$, if $A \notin \tau_s$, and if $\hat{A} \in \tau_d$,
then $$v_{\hat{A}}
=\frac{\tilde{w}_{\hat{A}}}{\parallel \tilde{w}_{\hat{A}} \parallel}
=\frac{\parallel {v}_{\hat{A}}(b^{\prime}(\hat{A})) \parallel}{\parallel \tilde{w}_{\hat{A}}(b^{\prime}(\hat{A})) \parallel}\tilde{w}_{\hat{A}}
=\frac{\parallel {v}_{\hat{A}}(b^{\prime}(\hat{A})) \parallel}{\parallel {w}_{\hat{A}}(b^{\prime}(\hat{A})) \parallel}\tilde{w}_{\hat{A}}$$
and 
$$\tilde{w}_{\hat{A}}(b^{\prime}(A))-\tilde{w}_{\hat{A}}(b(A))
=\mu_A \left( w_{A}(b^{\prime}(A))-w_A(b(A)) \right)=\mu_A w_{A}(b^{\prime}(A)).$$
Therefore, 
$\mu^{\prime}_A=\mu_A.$ 
\item If $A \in \tau_d$, if $A \notin \tau_s$, and if $\hat{A} \notin \tau_d$,
$$S_{\hat{A}}=\frac{\tilde{s}_{\hat{A}}}{\parallel \tilde{s}_{\hat{A}} \parallel}=\frac{\parallel {S}_{\hat{A}}(b^{\prime}(\hat{A})) \parallel}{\parallel \tilde{s}_{\hat{A}}(b^{\prime}(\hat{A})) \parallel}\tilde{s}_{\hat{A}}
=\frac{\parallel {S}_{\hat{A}}(b^{\prime}(\hat{A})) \parallel}{\parallel {s}_{\hat{A}}(b^{\prime}(\hat{A})) \parallel}\tilde{s}_{\hat{A}}$$
and $$\tilde{s}_{\hat{A}}(b^{\prime}(A))-\tilde{s}_{\hat{A}}(b(A))
=\mu_A \left( w_{A}(b^{\prime}(A))-w_A(b(A)) \right).$$
Therefore, 
$\mu^{\prime}_A=\mu_A.$
\item If $V(\sigma) \in \tau_d$, then
$m_{V(\sigma)}=\mu_{V(\sigma)}
\frac{\parallel \tilde{w}_{V(\sigma)} \parallel}
{\sqrt{{\sharp V(\sigma)}} \norm{\tilde{s}_{\sigma}(b_{\sigma})}}$,\\
where $\tilde{s}_{\sigma}(b_{\sigma})={s}_{\sigma}(b_{\sigma})
=\frac{{s}_{\sigma}(b_{\sigma})}{\sqrt{{\sharp V(\sigma)}} \norm{{s}_{\sigma}(b_{\sigma})}}$.\\
Thus, $m_{V(\sigma)}=\mu_{V(\sigma)}
\parallel \tilde{w}_{V(\sigma)} \parallel$, and $\mu^{\prime}_{V(\sigma)}=\mu_{V(\sigma)}$.
\end{itemize}

\medskip

\noindent{\em Proving that $\nu^{\prime}_i=\nu_i$, for $i \in \{1, \dots, \sigma\}$.}\\
\begin{itemize}
\item
When $i \geq 2$,\\
$S_{V({i-1})}=S_{i-1}=\frac{\tilde{s}_{i-1}}{\parallel \tilde{s}_{i-1} \parallel}= \frac{\parallel S_{i-1}(b^{\prime}_{i-1}) \parallel}{\parallel s_{i-1}(b^{\prime}_{i-1})\parallel} \tilde{s}_{i-1} $, $\tilde{s}_{i-1}(b_{\sigma})=\nu_{i}s_{i}(b_{\sigma})$, and\\
$\tilde{s}_{i-1}(b^{\prime}_{i})=\nu_{i}s_{i}(b^{\prime}_{i})$ if $V(i) \notin \tau_d$ and $\sharp V(i)  >  1$.
Therefore, $\nu^{\prime}_{i}=\nu_{i}$.
\item 
If $V \in \tau_d$, or if $\sharp V=1$,\\
$\nu^{\prime}_1=\ell_V$
where $\ell_V= \nu_1 \sqrt{\sharp V}\parallel \tilde{s}_V(b(V)) \parallel$, and $\tilde{s}_V(b(V))=s_V(b(V))=\frac{1}{\sqrt{\sharp V}}u_V$.
Therefore, $\nu_1=\ell_V$, and we are done.
\item If $V \notin \tau_d$ and if  $\sharp V > 1$,\\
$\nu^{\prime}_1=\ell_V\frac{ \langle S_V(b^{\prime}_{1}),s_1(b^{\prime}_{1}) \rangle }{ \langle s_1(b^{\prime}_{1}),s_1(b^{\prime}_{1}) \rangle }$, $\ell_V=\nu_1 \parallel \tilde{s}_1 \parallel$,
where $S_V=\frac{\tilde{s}_1}{\parallel \tilde{s}_1 \parallel}=
\frac{ \langle S_V(b^{\prime}_{1}),\tilde{s}_1(b^{\prime}_{1}) \rangle }{ \langle \tilde{s}_1(b^{\prime}_{1}),\tilde{s}_1(b^{\prime}_{1}) \rangle } \tilde{s}_1$, and $\tilde{s}_1(b^{\prime}_{1})=s_1(b^{\prime}_{1})$.
Thus, $\nu^{\prime}_1=\nu_1.$\\
\end{itemize}

\medskip

Thus, $r^{\tau} \circ p^{\tau} \circ \xi$ is the identity of $O \times W$.

\medskip

\noindent{\bf Checking that $r^{\tau}(p^{\tau}(\tilde{U}) \cap {C}^{\tau}_V(M)) \subseteq [0,\infty[^k \times \tilde{W}$.\/}

\medskip

Let $c=(c_A)_{A; A \in \tau} \in p^{\tau}(\tilde{U}) \cap {C}^{\tau}_V(M)$
with
$$c_A=\left\{\begin{array}{ll}\psi(A;\phi_{\infty};b(A))(\ell_A;u_A;m_A;v_A) & \mbox{if}
\;A \in \tau_d\\
\psi(A;\phi_{\infty})(\ell_A;S_A) & \mbox{if}
\; A \in (\tau \setminus \tau_d)\end{array}\right.$$

\begin{itemize}
\item If $A \in \tau_d$, (if $A \notin \tau_s$,) and if $\hat{A} \in \tau_d$,
then $\left(v_{\hat{A}}(b^{\prime}(A))-v_{\hat{A}}(b(A))\right)$ is a $(\geq 0)$
multiple of $\left(v_{A}(b^{\prime}(A))-v_{A}(b(A))\right)$ that is a positive multiple of $w_{A}(b^{\prime}(A))$, therefore,
$$\mu_A=\frac{ \parallel v_{\hat{A}}(b^{\prime}(A))-v_{\hat{A}}(b(A)) \parallel }
{\parallel  w_{A}(b^{\prime}(A))\parallel }
\frac{\parallel w_{\hat{A}}(b^{\prime}(\hat{A})) \parallel}{\parallel v_{\hat{A}}(b^{\prime}(\hat{A})) \parallel} \geq 0.$$
\item  If $A \in \tau_d$, if $A \notin \tau_s$, and if $\hat{A} \notin \tau_d$, similarly,
$$\mu_A=\frac{ \parallel S_{\hat{A}}(b^{\prime}(A))-S_{\hat{A}}(b(A)) \parallel }
{ \parallel  w_{A}(b^{\prime}(A)) \parallel }
\frac{\parallel s_{\hat{A}}(b^{\prime}(\hat{A})) \parallel}{\parallel S_{\hat{A}}(b^{\prime}(\hat{A})) \parallel} \geq 0.$$
\item If $V(\sigma) \in \tau_d$, then
$\mu_{V(\sigma)}=\frac{m_{V(\sigma)}}{\parallel \tilde{w}_{V(\sigma)}\parallel} \geq 0.$
\item
When $i \geq 2$, 
\begin{itemize}
\item If $i=\sigma$, and, if $V(\sigma) \in \tau_d$ or if $\sharp V(\sigma)=1$, then
$$\nu_{i}=\frac{\parallel S_{i-1}(b_{\sigma})\parallel }{ \parallel s_{i}(b_{\sigma})\parallel }
\frac{\parallel s_{i-1}(b^{\prime}_{i-1}) \parallel}{\parallel S_{i-1}(b^{\prime}_{i-1})\parallel} \geq 0.$$
\item Otherwise,
$$\nu_{i}=\frac{ \parallel S_{i-1}(b^{\prime}_{i})\parallel }{ \parallel s_{i}(b^{\prime}_{i})\parallel }
\frac{\parallel s_{i-1}(b^{\prime}_{i-1}) \parallel}{\parallel S_{i-1}(b^{\prime}_{i-1})\parallel} \geq 0$$
\end{itemize}
\item \begin{itemize}
\item If $V \in \tau_d$, or if $\sharp V=1$, $\nu_1=\ell_V \geq 0$.
\item If $V \notin \tau_d$ and if  $\sharp V > 1$,
$$\nu_1=\ell_V\frac{ \parallel S_V(b^{\prime}_{1})\parallel }{ \parallel s_1(b^{\prime}_{1})\parallel} \geq 0.$$

\end{itemize}
\end{itemize}

Therefore all the $\mu_A$ and the $\nu_i$ are positive.
Also, note that $r^{\tau}(p^{\tau}(c^0))=P^0$.

Now,  choose
$(\varepsilon  > 0)$ such that $[0,\varepsilon[^k \subset O$,
reduce $O$ into
$[0,\varepsilon[^k$ and set 
$$U^{\tau}=(r^{\tau})^{-1}(]-\varepsilon,\varepsilon[^k \times W).$$

Then $p^{\tau} \circ \xi(O \times W) \subseteq U^{\tau}$, 
$r^{\tau}(U^{\tau} \cap {C}^{\tau}_V(M)) \subseteq O \times W$.

\medskip

\noindent{\bf Checking that $p^{\tau} \circ \xi \circ r^{\tau}_{|U^{\tau} \cap {C}^{\tau}_V(M)}$ is the identity.\/}

Keep the above notation introduced to check that $r^{\tau}(p^{\tau}(\tilde{U}) \cap {C}^{\tau}_V(M)) \subseteq [0,\infty[^k \times \tilde{W}$
and assume that $c \in U^{\tau} \cap {C}^{\tau}_V(M)$.

\medskip

\noindent{\em Introducing more notation to prove that 
$\xi_A \circ r^{\tau}(c)=c_A$, for any $A \in \tau$.} \\
$$r^{\tau}(c)=((\nu_i)_{i \in \{1, \dots, \sigma\}};(\mu_A)_{A \in \tau_d};(s_i)_{i \in \{1, \dots, \sigma\}};(w_A)_{A \in \tau_d}).$$
Define $\tilde{w}_A$, for $A \in \tau_d$, $\tilde{s}_i$ and $\lambda_i$, for
$i \in \{1, \dots, \sigma\}$ as in the construction of $p^{\tau} \circ \xi$. Then
$$\xi_A(r^{\tau}(c))=\left\{\begin{array}{ll}\psi(A;\phi_{\infty};b(A))(\tilde{\ell}_A;\tilde{u}_A;\tilde{m}_A;\frac{\tilde{w}_{A}}{\parallel \tilde{w}_{A} \parallel}) & \mbox{if}
\;A \in \tau_d\\
\psi(A;\phi_{\infty})(\lambda_{i(A)} \parallel \tilde{s}_{i(A)} \parallel;
\frac{\tilde{s}_{i(A)}}{\parallel \tilde{s}_{i(A)} \parallel}) & \mbox{if}
\; A \in (\tau \setminus \tau_d)\end{array}\right.$$
where, for $A \in \tau_d$,
$$\tilde{\ell}_A=
\lambda_{i(A)} \sqrt{\sharp A}\parallel \tilde{s}_{i(A)}(b(A)) \parallel,$$
$$\tilde{u}_A=\frac{\tilde{s}_{i(A)}(b(A))}{\parallel \tilde{s}_{i(A)}(b(A)) \parallel},$$
and 
$$\tilde{m}_A=\left(\prod_{D \in \tau_d; A \subseteq D \subseteq V(i(A))} \mu_D \right)
\frac{\parallel \tilde{w}_{A} \parallel}{\sqrt{\sharp A}\parallel \tilde{s}_{i(A)}(b(A)) \parallel}.$$ 

\medskip

\noindent{\em Proving that $v_A=\frac{\tilde{w}_{A}}{\parallel \tilde{w}_{A} \parallel}$ when $A \in \tau_d$.}\\
Since $\norm{v_A}=1$, it suffices to prove that 
$$\tilde{w}_A=\frac{1}{\norm{w^1_A}}v_A$$
where because $w_A=\frac{1}{\norm{w^1_A}}w^1_A$
$$\norm{w^1_A}=\frac{\norm{w^1_A(b^{\prime}(A))}}{\norm{w_A(b^{\prime}(A))}}=\frac{\norm{v_A(b^{\prime}(A))}}{\norm{w_A(b^{\prime}(A))}}.$$
\begin{itemize}
\item When $A$ has no daughters, since $\tilde{w}_A=w_A=w^1_A=v_A$ and $\norm{w^1_A}=1$, we are done.
\item Assume that $\tilde{w}_C=\frac{1}{\norm{w^1_C}}v_C$ for any $C \in D(A)$.
Let $C \in D(A)$. \\
Since $c \in {C}^{\tau}_V(M)$, $\left((v_A)_{|C}-(v_A(b(C)))^C\right)$ is a positive multiple of $v_C$, while
$(v_C)_{|\{b(C),b^{\prime}(C)\}}$ is a positive multiple of $w_C$ that vanishes at $b(C)$. Therefore, 
$$\mu_C=\frac{1}{\norm{w^1_A}}\frac{\norm{v_A(b^{\prime}(C))-v_A(b(C))}}{\norm{w_C(b^{\prime}(C))}},$$
$$\frac{\mu_C}{\norm{w^1_C}}=\frac{1}{\norm{w^1_A}}\frac{\norm{v_A(b^{\prime}(C))-v_A(b(C))}}{\norm{v_C(b^{\prime}(C))}},$$
and 
$$\norm{w^1_A}\tilde{w}_A=w^1_A + \sum_{C \in D(A)}\frac{\norm{v_A(b^{\prime}(C))-v_A(b(C))}}{\norm{v_C(b^{\prime}(C))}} v_C=v_A.$$
\end{itemize}
This proves that $\tilde{w}_A=\frac{1}{\norm{w^1_A}}v_A$ by induction.

\medskip

\noindent{\em Proving that $S_A=\frac{\tilde{s}_{i(A)}}{\parallel \tilde{s}_{i(A)} \parallel}$ when $A \in (\tau \setminus \tau_d)$.}\\
If $\sharp A =1$, then $A=V(\sigma)$, $\tilde{s}_{i(A)}=s_{i(A)}=S_A$, and we are done.\\
Assume $A=V(i) \notin \tau_d$ and $\sharp A  > 1$.\\
We need to prove that, $S_i=S_{V(i)}=\frac{\tilde{s}_{i}}{\parallel \tilde{s}_{i} \parallel}$.
Again, it is enough to prove that $\tilde{s}_{i}=\frac{1}{\norm{s^1_i}}S_i$
where, since $s_i=\frac{1}{\norm{s^1_i}}s^1_i$, 

$$\norm{s^1_i}=\frac{\norm{s^1_i(b^{\prime}_i)}}{\norm{s_i(b^{\prime}_i)}}
=\frac{\norm{S_i(b^{\prime}_i)}}{\norm{s_i(b^{\prime}_i)}}$$
Let $C \in D(V(i))$, $C \neq V(i)_0$.\\
Then $\left(S_i(b^{\prime}(C))-S_i(b(C))\right)$ is a $(\geq 0)$ multiple of
$\left(v_C(b^{\prime}(C))=\norm{w^1_C}w_C(b^{\prime}(C))\right)$ (see the previous proof),
$$\mu_C=\frac{1}{\norm{s^1_i}}\frac{\norm{S_i(b^{\prime}(C))-S_i(b(C))}}{\norm{w_C(b^{\prime}(C))}},$$
and,
$$\frac{\mu_C}{\norm{w^1_C}}=\frac{1}{\norm{s^1_i}}\frac{\norm{S_i(b^{\prime}(C))-S_i(b(C))}}{\norm{v_C(b^{\prime}(C))}}.$$
Therefore
$$\begin{array}{ll}\norm{s^1_i}(\tilde{s}_{i})_{|(V(i) \setminus V(i)_0)} &
=s^1_i + \sum_{C \in D(V(i)); C \neq V(i)_0}\frac{\mu_C}{\norm{w^1_C}}\norm{s^1_i}
v_C\\
 &
=s^1_i + \sum_{C \in D(V(i)); C \neq V(i)_0}\frac{\norm{S_i(b^{\prime}(C))-S_i(b(C))}}{\norm{v_C(b^{\prime}(C))}}
v_C\\
& =(S_i)_{|(V(i) \setminus V(i)_0)}\end{array}$$
\begin{itemize}
\item When $V(\sigma) \notin \tau_d$, and when $i=\sigma$,\\
$V(i)_0=\emptyset$ and we can conclude.
\item When $V(\sigma) \in \tau_d$, and when $i=\sigma-1 \geq 1$, \\
By definition of $\nu_{\sigma}$, $$\nu_{\sigma}=\frac{\norm{S_i(b_{\sigma})}}{\norm{s_{\sigma}(b_{\sigma})}}\frac{\norm{s_i(b^{\prime}_i)}}{\norm{S_i(b^{\prime}_i)}}=\sqrt{\sharp V(\sigma)}\frac{\norm{S_i(b_{\sigma})}}{\norm{s^1_i}}.$$
Therefore
$$\begin{array}{ll}\norm{s^1_i}(\tilde{s}_{i})_{|V(i)_0} &
=\norm{s^1_i}\nu_{\sigma} \tilde{s}_{\sigma} =\norm{s^1_i}\nu_{\sigma} \left({s}_{\sigma} + \mu_{V(\sigma)} \tilde{w}_{V(\sigma)} \right)\\
 &
=\sqrt{\sharp V(\sigma)}\norm{S_i(b_{\sigma})}\left({s}_{\sigma} + \mu_{V(\sigma)} \norm{\tilde{w}_{V(\sigma)}}v_{V(\sigma)} \right)\\
&
=\sqrt{\sharp V(\sigma)}\norm{S_i(b_{\sigma})}\left(\left(\frac{{u}_{V(\sigma)}}{\sqrt{\sharp V(\sigma)}}\right)^{V(\sigma)} + m_{V(\sigma)} v_{V(\sigma)} \right)\\
& =(S_i)_{|V(\sigma)}\end{array}$$
since the latter right-hand side is a $(\geq 0)$ mutiple of the previous right-hand side, and since the norms of their values at $b_{\sigma}$ are the same.
\item When $V(i+1) \notin \tau_d$, and when $\tilde{s}_{i+1}=\frac{1}{\norm{s^1_{i+1}}}S_{i+1}$,
$$(\tilde{s}_{i})_{|V(i+1)}= \nu_{i+1} \tilde{s}_{i+1}= \frac{\nu_{i+1}}{\norm{s^1_{i+1}}} S_{i+1}$$
\begin{equation}
\label{eq317}
\nu_{i+1}=\left\{\begin{array}{ll} \frac{\norm{s^1_{i+1}}}{\norm{s^1_{i}}}
\frac{\norm{S_i(b_{\sigma})}}{\norm{S_{i+1}(b_{\sigma})}} &
\mbox{if} \; i+1=\sigma \mbox{, and if}\; \sharp V(\sigma)=1,\\
\frac{\norm{s^1_{i+1}}}{\norm{s^1_{i}}}
\frac{\norm{S_i(b^{\prime}_{i+1})}}{\norm{S_{i+1}(b^{\prime}_{i+1})}} & \mbox{otherwise.}
\end{array}\right.\end{equation}
Thus, in the second case,
$$\norm{s^1_{i}}(\tilde{s}_{i})_{|V(i+1)}=\frac{\norm{S_i(b^{\prime}_{i+1})}}{\norm{S_{i+1}(b^{\prime}_{i+1})}}S_{i+1}=(S_i)_{|V(i+1)}.$$
In any case, $\norm{s^1_{i}}(\tilde{s}_{i})_{|V(i+1)}=(S_i)_{|V(i+1)}$.
\end{itemize}
We conclude that $\tilde{s}_{i}=\frac{1}{\norm{s^1_{i}}}S_{i}$ for any $i$ such that $V(i) \notin \tau_d$, with a decreasing induction on $i$.
\medskip

\noindent{\em Proving that $\ell_A=\lambda_A \norm{\tilde{s}_A}$ when $A \in (\tau \setminus \tau_d)$.}\\
If such an $A$ exists, $V \notin \tau_d$. Let $A=V(i) \notin \tau_d$, and let us prove that $\ell_{V(i)}= \lambda_i \norm{\tilde{s}_i}$\\ where
$S_i=\norm{{s}^1_i}\tilde{s}_i$, $\norm{\tilde{s}_i}=\frac{1}{\norm{{s}^1_i}}$, $\lambda_i= \prod_{j=1}^{i}\nu_j$ and $\nu_1=\ell_V \norm{{s}^1_1}$. (The latter equality is obvious if $\sharp V=1$, otherwise $V \notin \tau_d$ and $\nu_1=\ell_V \frac{\norm{S_V(b^{\prime}_1)}}{\norm{\tilde{s}_V(b^{\prime}_1)}}$.)\\
If $i=1$, we are done.\\
Otherwise, for any $j < i$, according to \ref{eq317},
$$\nu_{j+1}=\frac{\norm{s^1_{j+1}}}{\norm{s^1_{j}}}
\frac{\norm{(S_j)_{|V(j+1)}}}{\norm{S_{j+1}}}$$
Therefore,
$$\begin{array}{ll} \lambda_i \norm{\tilde{s}_i}&=\norm{\tilde{s}_i}\prod_{j=1}^{i}\nu_j\\
&=\ell_V\frac{\norm{{s}^1_1}}{\norm{{s}^1_i}}\prod_{j=1}^{i-1}\nu_{j+1}\\
&=\ell_V\prod_{j=1}^{i-1}\frac{\norm{(S_j)_{|V(j+1)}}}{\norm{S_{j+1}}}
\end{array}$$
where $\ell_{V(j)}(S_j)_{|V(j+1)}=\ell_{V(j+1)}S_{j+1}$.
It follows that $\ell_{V(i)}= \lambda_i \norm{\tilde{s}_i}$ by induction on $i$.

\medskip

\noindent{\em Proving that $\tilde{u}_{V(\sigma)}=u_{V(\sigma)}$, $\tilde{m}_{V(\sigma)}=m_{V(\sigma)}$ and $\tilde{\ell}_{V(\sigma)}=\ell_{V(\sigma)}$  when ${V(\sigma)} \in \tau_d$.}\\
Let  ${V(\sigma)} \in \tau_d$. We already know that $v_{V(\sigma)}=\frac{\tilde{w}_{V(\sigma)}}{\norm{\tilde{w}_{V(\sigma)}}}$.
In particular,
$$\begin{array}{ll}\tilde{s}_{\sigma}&=s_{\sigma} + \mu_{V(\sigma)} \tilde{w}_{V(\sigma)}\\
&=s_{\sigma} + \mu_{V(\sigma)} \norm{\tilde{w}_{V(\sigma)}}{v}_{V(\sigma)}\\
&=\left(\frac{u_{V(\sigma)}}{\sqrt{\sharp V(\sigma)}} \right)^{V(\sigma)} + m_{V(\sigma)}{v}_{V(\sigma)},\\ \end{array}$$
$$\tilde{u}_{V(\sigma)}=\frac{\tilde{s}_{\sigma}(b({V(\sigma)}))}{\parallel \tilde{s}_{\sigma}(b({V(\sigma)})) \parallel}=u_{V(\sigma)},$$
$$\sqrt{\sharp V(\sigma)}\parallel \tilde{s}_{\sigma}(b({V(\sigma)})) \parallel=1,$$
$$\tilde{m}_{V(\sigma)}= \mu_{V(\sigma)}\frac{\norm{\tilde{w}_{V(\sigma)}}}{\sqrt{\sharp V(\sigma)}\parallel \tilde{s}_{\sigma}(b({V(\sigma)})) \parallel}=m_{V(\sigma)}.$$
We are left with the proof that $\tilde{\ell}_{V(\sigma)}=\ell_{V(\sigma)}$.
\begin{itemize}
\item
When $\sigma=1$,\\
$V(\sigma)=V$, $\tilde{\ell}_{V}=\lambda_1=\nu_1=\ell_V$, and we are done.
\item 
When $\sigma > 1$,\\
Since $p^{\tau} \circ \xi \circ r^{\tau}(c) \in {C}^{\tau}_V(M)$,\\
$$\tilde{\ell}_{V(\sigma)} \left(\left(\frac{\tilde{u}_{V(\sigma)}}{\sqrt{\sharp V(\sigma)}} \right)^{V(\sigma)} + \tilde{m}_{V(\sigma)}
\frac{\tilde{w}_{V(\sigma)}}{\parallel \tilde{w}_{V(\sigma)}\parallel} \right)
=\lambda_{\sigma -1} (\tilde{s}_{\sigma -1})_{|V(\sigma)},$$
Thus, $\tilde{\ell}_{V(\sigma)}\tilde{s}_{\sigma}=\ell_{V(\sigma -1)}(S_{V(\sigma -1)})_{|V(\sigma)}.$\\
Since $c \in {C}^{\tau}_V(M)$,\\
$${\ell}_{V(\sigma)}\tilde{s}_{\sigma}=\ell_{V(\sigma -1)}(S_{V(\sigma -1)})_{|V(\sigma)}.$$
This implies that $\tilde{\ell}_{V(\sigma)}=\ell_{V(\sigma)}$.
\end{itemize}
\medskip

\noindent{\em Proving that $\ell_A=\tilde{\ell}_A$, $\tilde{m}_A=m_A$
and $\tilde{u}_A=u_A$ when $A \in \tau_d$, $A \neq V(\sigma)$.}\\
We already know that $\Pi_V(c_V)=\Pi_V(\xi_V \circ r^{\tau}(c))$ in $M^V$.
This map from $V$ to $M$ may be written as $\phi_{\infty} \circ f$
where $f \in (\RR^3)^V = \ell_V S_V$ if $V \notin \tau_d$, and 
$f=\ell_V \tilde{s}_{1}$ if $V \in \tau_d$.
Since both $c$ and $p^{\tau} \circ \xi \circ r^{\tau}(c)$ satisfy $(C1)$, then
$\ell_A=\norm{f(b(A))}\sqrt{\sharp A}$ and $\tilde{\ell}_A=\norm{f(b(A))}\sqrt{\sharp A}$.
Therefore, $\ell_A=\tilde{\ell}_A$.\\
Now, $f_{|A}$ is a $(\geq 0)$ multiple of 
$\left( (u_A/ \sqrt{\sharp A})^A +m_A v_A) \right)$ and 
$\left( (\tilde{u}_A/ \sqrt{\sharp A})^A +\tilde{m}_A v_A \right)$.
Thus, if $f_{|A}\neq 0$, using $$\left( (u_A/ \sqrt{\sharp A})^A +m_A v_A \right)(b(A))=(u_A/ \sqrt{\sharp A})\;\;\;\mbox{and}\;\;\;\norm{u_A}=1,$$ we easily
conclude that 
$$(u_A/ \sqrt{\sharp A})^A +m_A v_A =(\tilde{u}_A/ \sqrt{\sharp A})^A +\tilde{m}_A v_A.$$
Now, if $f_{|A} = 0$, since $f_{|V(i(A))}$ is a $(\geq 0)$ multiple of 
$\tilde{s}_{i(A)}$, and since $\tilde{s}_{i(A)}$ does not vanish on $A$, we deduce that $f_{|V(i(A))}=0$.
Then $(C3)$ implies that $(\tilde{s}_{i(A)})_{|A}$ is a $(\geq 0)$ multiple of 
$\left( (u_A/ \sqrt{\sharp A})^A +m_A v_A) \right)$ and 
$\left( (\tilde{u}_A/ \sqrt{\sharp A})^A +\tilde{m}_A v_A) \right)$, and we conclude as before that 
$$(u_A/ \sqrt{\sharp A})^A +m_A v_A =(\tilde{u}_A/ \sqrt{\sharp A})^A +\tilde{m}_A v_A.$$
This implies of course that $\tilde{m}_A=m_A$
and $\tilde{u}_A=u_A$. \eop

This finishes the proof of Lemma~\ref{lemdifconfinfini} and thus Proposition~\ref{propdifconf} is proved.

\eop

To prove Lemma~\ref{lempropconsadface} in this case, we again look at the above proof. Here, $k=1$ if and only if $\tau_s=\{V\}$ and $\tau_d=\emptyset$, that is if and only if $\Pi_{\infty}(c^0_V)$ is an injective map from $V$ to $(T_{\infty}M \setminus 0)$ (up to dilation),
that is if and only if $c^0 \in F(\infty;V)$.
Lemma~\ref{lempropconsadface} is now proved.\eop

\subsection{Proof of Proposition~\ref{propcdeuxcoinc}}
\label{subpropcdeuxcoinc}
First define a smooth map of the following form
$$\begin{array}{llll}H:&C_2(M) &\longrightarrow &C(V;M) \times C(\{1\};M) \times  C(\{2\};M)\\
&c &\mapsto &(p_2(c),r_1(c),r_2(c))\end{array}$$
whose image will be in $C_V(M)$.

Since $C_2(M)$ is a blow-up of $C(V;M)$ along $\infty \times C_1(M)$ and 
$C_1(M) \times \infty$, we get the canonical smooth projection
$$p_2:C_2(M) \longrightarrow C(V;M) (\hfl{\Pi_V} M^V).$$
Let $i \in \{1,2\}$.
Let $p_{\{i\}}: M^V \longrightarrow M^{\{i\}}$ be the canonical restriction, and let $$\tilde{r}_i=p_{\{i\}} \circ \Pi_V \circ p_2: C_2(M) \longrightarrow M.$$
Let $r_i$ denote the restriction 
of $\tilde{r}_i$ from $ (\Pi_V \circ p_2)^{-1}\left((M \setminus \infty)^2\right)$ to $(M \setminus \infty)$. 
We are now going to define a smooth extension of the $r_i$ to $C_2(M)$ so that $H(C_2(M)) \subseteq C_V(M)$.

Let $\Pi_1:C_1(M) \longrightarrow M$ be the canonical projection.
Since the blow-up of $M \times (M\setminus \infty)$ along $\infty \times (M \setminus \infty)$ is canonically diffeomorphic to the product of the blow-up 
of $M$ at $\infty$ by $(M \setminus \infty)$, it is easy to observe the following lemma.

\begin{lemma}
For any two disjoint open subsets $V_1$ and $V_2$ of $M$, there is a canonical 
diffeomorphism
$$(\Pi_V \circ p_2)^{-1}(V_1 \times V_2) \hfl{(r_1,r_2)} \Pi_1^{-1}(V_1) \times \Pi_1^{-1} (V_2)$$
where $r_1$ and $r_2$ coincide with the previous maps $r_1$ and $r_2$ wherever
it makes sense.
\end{lemma}
\eop

Let us now prove that $r_i$ extends to a smooth projection from $C_2(M)$ onto
$C_1(M)$.
This extension will be necessarily unique and canonical, it will be denoted by $r_i$.
By symmetry, we only consider the case $i=1$.
It remains to define $r_1$ on $(\Pi_V \circ p_2)^{-1}(\infty,\infty)$
and to prove that it is smooth there. 
The canonical smooth projection $\Pi_{M^2(\infty,\infty)}$ from 
$C_2(M)$ to $M^2(\infty,\infty)$ maps $(\Pi_V \circ p_2)^{-1}(\infty,\infty)$
to $\left(ST_{(\infty,\infty)}M^2\right)$.
In turn, $ST_{(\infty,\infty)}M^2 \setminus S(0 \times T_{\infty}M)$
projects onto $ST_{\infty}M$ as the first coordinate.
It is easy to see that this smoothly extends the definition of $r_1$ outside $\Pi_{M^2(\infty,\infty)}^{-1}\left(S(0 \times T_{\infty}M) \right)$. To conclude, we recall the structure of $C_2(M)$ near $\Pi_{M^2(\infty,\infty)}^{-1}\left(S(0 \times T_{\infty}M) \right)$.
According to Proposition~\ref{propblodifdeux}, since the normal bundle of $\infty \times M$ at $(\infty,\infty)$ in $M^2$ is $(T_{\infty}M \times \{0\})$, $\Pi_{M^2(\infty,\infty)}^{-1}\left(S(0 \times T_{\infty}M) \right)$ is the product
$ST_{\infty}M \times S(0 \times T_{\infty}M)$. Define $r_1$ as the projection
on the fiber $ST_{\infty}M \subset C_1(M)$ in this product.
From the chart $\phi_{\infty}^2: (\RR^3)^2 \longrightarrow M^2$, that induces the chart
$$\psi_1:[0,\infty[ \times S((\RR^3)^2)\cong S^5 \longrightarrow M^2(\infty,\infty)$$
such that $\Pi_V \circ \psi_1(\lambda;(x,y))= (\phi_{\infty}( \lambda x),\phi_{\infty}( \lambda y))$
that in turn, induces the chart near
$S( 0 \times T_{\infty}M )$ 
$$\psi_2:[0,\infty[ \times \RR^3 \times S^2  \longrightarrow M^2(\infty,\infty)$$
such that $\Pi_V \circ \psi_2(\lambda;(x,y))= (\phi_{\infty}( \lambda x),\phi_{\infty}( \lambda y))$, we get a chart
$$\psi_3:[0,\infty[ \times ([0,\infty[ \times S^2) \times S^2 \longrightarrow C_2(M)$$
such that $\Pi_V \circ p_2 \circ \psi_3(\lambda;\mu;x;y))= (\phi_{\infty}( \lambda \mu x),\phi_{\infty}( \lambda y))$. Using a similar chart for $C_1(M)$, $r_1$ will read
$$(\lambda;\mu;x;y) \mapsto (\lambda\mu;x)$$
and is smooth.

Now, our map $H$ is well-defined and smooth.
The elements of $H(C_2(M))$ satisfy $(C1)$ and $(C3)$ (of Lemmas~\ref{lemcond1} and \ref{lemcond3}). Thus, since $V$ has two elements, $H(C_2(M)) \subseteq C_V(M)$ and we have defined a smooth map $H$ from
$C_2(M)$ to $C_V(M)$, that extends the identity of $\breve{C}_2(M)$.

Let $p_V: C_V(M)  \longrightarrow C(V;M)$ be the canonical projection.
Define 
$$K: C_V(M)  \longrightarrow  C_2(M)$$
so that
$$K(c_V,c_1,c_2) = \left\{\begin{array}{ll}
				p_2^{-1}(c_V) &\;\mbox{if}\; c_V \notin (\infty \times C_1(M)) \cup (C_1(M) \times \infty)\\
(r_1,r_2)^{-1}(c_1,c_2)&\;\mbox{if}\;\Pi_V(c_V) \;\mbox{is not constant.}
\end{array} \right.
$$
The map $K$ is consistently defined outside $p_V^{-1}((\infty \times \partial C_1(M)) \cup (\partial C_1(M) \times \infty))$, and it is smooth there.
We shall extend $K$ by the canonical identifications on $p_V^{-1}((\infty \times \partial C_1(M)) \cup (\partial C_1(M) \times \infty))$ and
use the charts of $C_V(M)$, near $p_V^{-1}((\infty \times \partial C_1(M)) \cup (\partial C_1(M) \times \infty))$ to prove the smoothness.
For example, $p_V^{-1}((\infty \times \partial C_1(M))$ is the subset of 
$S(\{0\} \times T_{\infty}M) \times S(T_{\infty}M \times \{0\}) \times S(\{0\} \times T_{\infty}M)$ where the first and third coordinate coincide. There,
$\tau_d=\emptyset$, $\tau_s=\{V,\{1\}\}$, and the chart
$\xi$ of Subsection~\ref{subsecinf} reads:\\
$(\nu_1,\nu_2,s_1=s_V,s_2=s_{\{1\}}) \mapsto ( \psi(V;\phi_{\infty})(\nu_1\norm{s_V + \nu_2 s_{\{1\}}};\frac{s_V + \nu_2 s_{\{1\}}}{\norm{s_V + \nu_2 s_{\{1\}}}}),$\\
$\psi(\{1\};\phi_{\infty})(\nu_1\nu_2;s_{\{1\}}),
\psi(\{2\};\phi_{\infty})(\nu_1;{s_V(2)}))$.\\
Mapping $\xi(\nu_1,\nu_2,s_1=s_V,s_2=s_{\{1\}})$ to $\psi_3(\nu_1,\nu_2,s_{\{1\}},s_V(2))$ (where $\psi_3$ is defined above in this subsection)
smoothly extends $K$ to $p_V^{-1}((\infty \times \partial C_1(M))$. Similarly, $K$ smoothly extends to $p_V^{-1}(\partial C_1(M) \times \infty)$.
Then $K$ and $H$ are smooth maps that extend the identity  of $\breve{C}_2(M)$ that is dense in both spaces. Therefore  $K$ and $H$
are inverse of each other and they are diffeomorphisms.
\eop

\newpage
\begin{center}{\bf \large Terminology}\end{center}
\addcontentsline{toc}{section}{Terminology}
\begin{theindex}

  \item blow-up, 37, 38

  \indexspace

  \item dilation, 37

  \indexspace

  \item edge-orientation, 6

  \indexspace

  \item form
    \subitem admissible, 18
    \subitem antisymmetric, 5
    \subitem fundamental, 5

  \indexspace

  \item half-edge, 5

  \indexspace

  \item Jacobi diagram, 5
    \subitem automorphism of, 5, 16
    \subitem edge-oriented, 16
    \subitem labelled, 16
    \subitem orientation of, 5

  \indexspace

  \item linking number, 12

  \indexspace

  \item orientation of a finite set, 6

  \indexspace

  \item Pontryagin class, 8, 29
  \item Pontryagin number, 8, 9

  \indexspace

  \item trivialisation standard near $\infty$, 4

  \indexspace

  \item vertex-orientation, 6

\end{theindex}

\newpage
\begin{center}{\bf \large Notation}\end{center}
\addcontentsline{toc}{section}{Notation}
\begin{theindex}
\begin{multicols}{3}
  \item $\CA(\emptyset)$, 10
  \item $\CA_n(\emptyset)$, 5
  \item AS, 6
  \item $\sharp \mbox{Aut}(\Gamma)$, 5, 16

  \indexspace

  \item $b(A)$, 48, 55
  \item $b^{\prime}(A)$, 51, 59
  \item $B^3(r)$, 4

  \indexspace

  \item $C(A;b(A);\tau)$, 48
  \item $C(A;M)$, 37, 40
  \item $\breve{C}_n(M)$, 14
  \item $\breve{C}_{V(\Gamma)}(M)$, 6
  \item $C_1(M)$, 3, 37
  \item $\tilde{C}_V(M)$, 44
  \item $C_2(M)$, 3, 12, 14, 19, 45
  \item $C_V(M)$, 14, 37, 43, 44
  \item $C(V(i);\tau)$, 55

  \indexspace

  \item $\partial_1(C_{V}(M))$, 15
  \item $\partial_1(S_{V}(X))$, 16
  \item $\delta_n$, 17, 18

  \indexspace

  \item ${\cal E}_n$, 17
  \item $E(\Gamma)$, 5
  \item $E_1$, 9, 10, 27, 34

  \indexspace

  \item $F(B)$, 15
  \item $f(B)(X)$, 15
  \item $F(\infty;B)$, 14
  \item $F(V)$, 25

  \indexspace

  \item $\Gamma_B$, 20
  \item $G_M(\rho)$, 28, 32

  \indexspace

  \item $\HH$, 7
  \item $H(\Gamma)$, 5

  \indexspace

  \item $I_{\Gamma,F}$, 19
  \item $I_{\Gamma}(M;\Omega)$, 17
  \item $I_{\Gamma}(\omega_M)$, 7
  \item $I_{\Gamma}(\omega_T)$, 10
  \item IHX, 6, 23
  \item $\iota$, 5, 10
  \item $\overline{\iota}$, 5
  \item $i^2(m^{\CC}_r)$, 9, 30, 31

  \indexspace

  \item $\xi$, 10, 11
  \item $\xi_n$, 10, 18, 24
  \item $\xi_1$, 34, 36

  \indexspace

  \item $\lambda$, 3, 11
  \item $\lambda_W$, 11
  \item $L_M$, 12

  \indexspace

  \item $M^A(\infty^A)$, 37
  \item $m_r$, 30
  \item $m^{\CC}_r$, 30
  \item $M^2(\infty,\infty)$, 3

  \indexspace

  \item $O(A;b(A);\tau)$, 48
  \item $\Omega$, 16
  \item $\omega_T$, 10, 17, 36
  \item $O(V(i);\tau)$, 55

  \indexspace

  \item $p_A$, 37
  \item $p_e$, 6, 10
  \item $P(\Gamma)$, 16
  \item $\phi_{\infty}$, 39
  \item $\Pi_A$, 37
  \item $P_i(\Gamma)$, 16
  \item $\Pi_{\infty}(c^0_{A})$, 54
  \item $\Pi_{\infty,d}(c^0_{A})$, 54
  \item $p_M(\tau_M)$, 5
  \item $p_1$, 8, 9, 29, 32
  \item $\psi(A;\phi;b)$, 41
  \item $\psi(A;\phi_{\infty})$, 42
  \item $\psi(A;\phi_{\infty};b)$, 42
  \item $\psi(g)$, 9

  \indexspace

  \item $\rho$, 7, 10, 28, 30

  \indexspace

  \item $\breve{S}_n(X)$, 9, 14
  \item $\breve{S}_V(X)$, 9
  \item $S_i(T_{\infty}M^B)$, 14
  \item $S_2(E_1)$, 10, 17, 34
  \item $S(V)$, 3
  \item $S_V(\RR^3)$, 45
  \item $S_{V}(X)$, 14

  \indexspace

  \item $\tau(c^0)$, 48, 54
  \item $\tau_d$, 54, 55
  \item $\tau_s$, 54, 55

  \indexspace

  \item $V(\Gamma)$, 5

  \indexspace

  \item $X(Z)$, 37, 38

  \indexspace

  \item $Z(M)$, 10, 11
  \item $Z(M;\tau_M)$, 10, 17
  \item $z_n(E;\omega)$, 24, 25
  \item $Z_n(M)$, 7
  \item $Z_n(M;\tau_M)$, 7
  \item $Z_n(\omega_M)$, 18
  \item $z_n(\tau_M)$, 17
\end{multicols}
\end{theindex}


\end{document}